\documentclass[10pt,a4paper]{article}

 \usepackage{amssymb,amsmath,amsfonts,amssymb, amsthm}%,numbysec}
 \usepackage{graphics,graphicx,color}

 \usepackage{amsmath}    % remove % for AMS-LaTeX %\usepackage{amssymb}
 \advance\hoffset by 5mm

 \title{Global existence and regularity for the full\\
   coupled  Navier-Stokes and Q-tensor  system }
 \date{\today}
\author{Marius Paicu\footnote{Universit\'e Paris-Sud, Laboratoire de Math\'ematiques, 91405 Orsay Cedex, France. E-mail: marius.paicu@math.u-psud.fr }\,\,\,\,  and Arghir Zarnescu \footnote{Mathematical Institute, 24-29 St. Giles', Oxford, OX1 3LB, United Kingdom. E-mail: zarnescu@maths.ox.ac.uk} }

\newtheorem{remark}{Remark}
 \newtheorem{lemma} {Lemma}
 \newtheorem{proposition}{Proposition}
 \newtheorem{theorem} {Theorem}

\begin{document}

  \maketitle

\section{Introduction}

\par In this paper we continue the study, initiated in \cite{pz1}, concerning  the global existence of solutions for a system describing the evolution of a nematic liquid crystal fluid. This system is a prototypical example of a certain type of non-Newtonian complex fluids, in which the stress tensor of the fluid has, in addition to the usual Newtonian part, a  component due to the presence of particles embedded in the flow, namely the liquid crystal molecules. The evolution of the flow is influenced by the presence of these particles and on the other hand the evolution of the flow affects the direction and position of the liquid crystal molecules.  This situation is modelled through a forced Navier-Stokes system, describing the flow,  coupled with a parabolic-type system describing the evolution of the nematic crystal director fields ($Q$-tensors, that is traceless and symmetric $d$-by-$d$ matrices, $d=2,3$). 

\smallskip\par In our previous work, \cite{pz1}, we assumed that a certain parameter, $\xi$, is zero, which had the effect of cancelling certain terms. In the current work we do not make this assumption and study the full system, observing that the presence of these additional terms has a non-trivial effect, namely the quadruply exponential increase of the high norms (that will be detailed below). We also estimate differently certain terms already existent in the simplified system and improve the estimates in \cite{pz1}. 

\par The full coupled system has, as well as the simplified version in \cite{pz1}, a Lyapunov functional made of two parts: the free energy due to  the director fields and the kinetic energy of the fluid. This functional describes, from a physical point of view, the dissipation of the energy of the complex fluid.

\smallskip
\par      In the first part of the paper we use the apriori bounds on the solution (provided by the energy dissipation) to prove the existence of  global weak solutions in the natural energy space. In the second part,   we study the case where the fluid evolves in the two dimensional space and  prove the existence of a global regular solution issued from, an appropriately regular, initial data. 
  In the two dimensional space we also show that for an appropriately regular initial data the weak and the strong solutions  coincide. 
  
\medskip
  \par The main part of our  study concerns the high regularity of the solutions that start from regular (enough) initial data. We use, at a higher level of regularity,  the cancellations that  made possible the existence of the Lyapunov functional, to avoid estimating certain terms with maximal number of derivatives. Thus we show that for this type of complex fluids the existence of an energy dissipation is intrinsically related to the high regularity of the solutions. Moreover, the differential inequality relating the high Sobolev norms of the solution (inequality that allows us to obtain uniform bounds in high Sobolev spaces), is not completely classical and takes a form that is different from, for  example, the classical situation of global wellposedness for incompressible Euler equation in two dimensions. Indeed, in our proof we use the  logarithmic Sobolev embedding of $H^{1+\epsilon}$ in $L^\infty$ in conjunction with   the precise growth of the constant of the Sobolev embedding of $H^1$ in any $L^p$ (which is $C\sqrt{p})$, and an optimal choice of the Lebsegue index $p$ depending on the norm of the solution. This way we obtain a differential inequality with a double-logarithmical correction and this allows us to obtain a global in time control of high Sobolev norms of the solution.  

\medskip
 \par     There exist several competing theories that attempt to capture the complexity of nematic liquid crystals, and a comparative discussion and further references are available for instance in \cite{lin&liu}, \cite{mz}. In the present paper we use one of the most comprehensive description of nematics, the $Q$-tensor description, proposed by P.G. de Gennes \cite{dg}. There exist various specific models that all use the $Q$-tensor description and a comparative discussion of the main models is available for instance in \cite{sonnet}.
      In this paper we  use a model proposed by Beris and Edwards \cite{berised}, that one can find  in the physics literature for instance in \cite{equationsELreduction}, \cite{equations}. An important feature of this model is that if one assumes smooth solutions and one formally takes $Q(x)=s_+(n(x)\otimes n(x)-\frac{1}{d}Id)$, with $s_+$ a constant (depending on the parameters of the system, see for instance \cite{mz})  and $n:\mathbb{R}^d\to\mathbb{S}^{d-1}$ smooth (with $d$ the dimension of the space), then the equations reduce (see \cite{equationsELreduction}) to the generally accepted equations of Ericksen, Leslie and Parodi \cite{leslie&ericksen}. The system we study is related structurally to other models of complex fluids coupling a transport equation with a forced Navier-Stokes system  \cite{cftz}, \cite{constantin&masmoudi}, \cite{constantin&seregin}, \cite{chemin&masmoudi},  \cite{lions&masmoudi}, \cite{linandco}, \cite{masmoudiFENE}. In our case the Navier-Stokes equations are coupled with a parabolic type system, but we also have two more derivatives (than in the previously mentioned models)  in the forcing term of the Navier-Stokes equations.  The Ericksen-Leslie-Parodi system describing nematic liquid crystals, whose structure is closer to our system (but  that has one less derivative in the forcing term of the Navier-Stokes equations) was studied in \cite{lin}, \cite{linLE}, \cite{linarma}.

      In the following we  use a partial Einstein summation convention, that is we  assume summation over repeated {\it greek} indices, but not over the repeated  {\it latin} indices. We consider the equations as described in \cite{equationsELreduction},~\cite{equations} but assume that the fluid has constant  density in time.
      We denote

\begin{equation}
 S(\nabla u,Q)\stackrel{\rm{def}}{=}(\xi D+\Omega)(Q+\frac{1}{d}Id)+(Q+\frac{1}{d}Id)(\xi D-\Omega)-2\xi (Q+\frac{1}{d}Id)\textrm{tr}(Q\nabla u)
 \label{rel:defS}
\end{equation} where $D\stackrel{def}{=}\frac{1}{2}\left(\nabla u+(\nabla u)^T\right)$ and $\Omega\stackrel{def}{=}\frac{1}{2}
\left(\nabla u-(\nabla u)^T\right)$ are the symmetric part and the antisymmetric part, respectively, of the velocity gradient matrix $\nabla u$. The constand $d$ is the dimension of the space and $Q$ is a function on $\mathbb{R}^d$ with values into $S_0^{(d)}$(see the notations paragraph below). The term $S(\nabla, Q)$ appears in the equation of motion of the order-parameter, $Q$, and describes how the flow gradient rotates and stretches the order-parameter. The constant $\xi$ depends on the molecular details of a given liquid crystal and measures the ratio between the tumbling and the aligning effect that a shear flow would exert over the liquid crystal directors.

      We also denote:

\begin{equation}
 H\stackrel{\rm{def}}{=}-aQ+b[Q^2-\frac{\textrm{tr}(Q^2)}{3}Id]-cQ\textrm{tr}(Q^2)+L\Delta Q
\label{rel:defH}
\end{equation} where $L>0$. It will also be convenient to denote 

\begin{equation}
F\stackrel{\rm{def}}{=}H-L\Delta Q
\label{rel:defF}
\end{equation}

      With the notations above we have the coupled system:

\begin{equation}
\left\{\begin{array}{l}
      (\partial_t+u\cdot \nabla)Q-S(\nabla u,Q)=\Gamma H \\
         \partial_t u_\alpha+u_\beta   \partial_\beta u_\alpha=\nu   \partial_{\beta\beta} u_\alpha+   \partial_\alpha p+   \partial_\beta \tau_{\alpha\beta}+   \partial_\beta\sigma_{\alpha\beta}\\
          \partial_\gamma u_\gamma=0
     \end{array}\right.
     \label{system}
\end{equation} where $\Gamma>0,\nu>0$ and  we have the symmetric part of the additional stress tensor:

\begin{equation}
\tau_{\alpha\beta}\stackrel{\rm{def}}{=}-\xi \left(Q_{\alpha\gamma}+\frac{\delta_{\alpha\gamma}}{d}\right)H_{\gamma\beta}-\xi H_{\alpha\gamma}\left(Q_{\gamma\beta}+\frac{\delta_{\gamma\beta}}{d}\right)+2\xi (Q_{\alpha\beta}+\frac{\delta_{\alpha\beta}}{d})Q_{\gamma\delta}H_{\gamma\delta}-L\left(   \partial_\beta Q_{\gamma\delta}   \partial_\alpha Q_{\gamma\delta}+\frac{\delta_{\alpha\beta}}{d}Q_{\nu\varepsilon}Q_{\nu\varepsilon}\right)\end{equation} and an antisymmetric part:

\begin{equation}
\sigma_{\alpha\beta}\stackrel{\rm{def}}{=}Q_{\alpha\gamma}H_{\gamma\beta}-H_{\alpha\gamma}Q_{\gamma\beta}
\end{equation}

\

      We also need to assume from now on that 
\begin{equation}
c>0
\label{c+}
\end{equation} This assumption is necessary from a modelling point of view (see \cite{apala}, \cite{mz}) so that the energy $\mathcal{F}$ (see next section, relation (\ref{freeenergy})) is bounded from below, and it is also necessary for having global solutions (see Proposition ~\ref{prop:aprioriest} and its proof).

\smallskip     {\bf Notations and conventions}  Let
$S_0^{(d)}\subset \mathbb{M}^{d\times d}$ denote the space of Q-tensors in dimension $d$,  i.e.
 $$S_0^{(d)}\stackrel{\rm{def}}{=} \left\{Q \in \mathbb{M}^{d\times d};
Q_{ij}=Q_{ji},\textrm{tr}(Q) = 0, i,j=1,\dots,d \right\}$$
      We use the Frobenius norm of a matrix $
 \left| Q \right|\stackrel{\rm{def}}{=}\sqrt{\textrm{tr}Q^2} =\sqrt{ Q_{\alpha\beta}
Q_{\alpha\beta}}$ and define Sobolev spaces of $Q$-tensors in terms of this norm. For instance $H^1(\mathbb{R}^d,S_0^{(d)})\stackrel{\rm{def}}{=}\{Q:\mathbb{R}^d\to S_0^{(d)}, \int_{\mathbb{R}^d} |\nabla Q(x)|^2+|Q(x)|^2\,dx<\infty\}$. For $A,B\in S_0$ we denote $A\cdot B=\textrm{tr}(AB)$ and $|A|=\sqrt{\textrm{tr}(A^2)}$. We also denote  $|\nabla Q|^2(x)\stackrel{\rm{def}}{=}Q_{\alpha\beta,\gamma}(x)Q_{\alpha\beta,\gamma}(x)$ and $|\Delta Q|^2(x)\stackrel{\rm{def}}{=}\Delta Q_{\alpha\beta}(x)\Delta Q_{\alpha\beta}(x)$. We recall also that $\Omega_{\alpha\beta}\stackrel{\rm{def}}{=}\frac{1}{2}\left(   \partial_\beta u_\alpha-   \partial_\alpha u_\beta\right)$ and $u_{\alpha,\beta}\stackrel{\rm{def}}{=}   \partial_\beta u_\alpha$, $Q_{ij,k}\stackrel{\rm{def}}{=}   \partial_k Q_{ij}$.

\section{ The  energy decay and apriori estimates}
      Let us denote the free energy of the director fields:

\begin{equation}
\mathcal{F}(Q)=\int_{\mathbb{R}^d} \frac{L}{2}|\nabla Q|^2+\frac{a}{2}\textrm{tr}(Q^2)-\frac{b}{3}\textrm{tr}(Q^3)+\frac{c}{4}\textrm{tr}^2(Q^2)\,dx
\label{freeenergy}
\end{equation}

      In the absence of the flow, when $u=0$ in the equations (\ref{system}), the free energy is a Lyapunov functional of the system. If $u\not=0$ we still have a Lyapunov functional for (\ref{system}) but this time one that includes the kinetic energy of the system. More precisely we have:

\begin{proposition}
The system (\ref{system}) has a Lyapunov functional:
\begin{equation}
E(t)\stackrel{\rm{def}}{=}\frac{1}{2}\int_{\mathbb{R}^d}|u|^2(t,x)\,dx+\int_{\mathbb{R}^d}\frac{L}{2} |\nabla Q|^2(t,x)+\frac{a}{2}\textrm{tr}(Q^2(t,x))-\frac{b}{3}\textrm{tr}(Q^3(t,x))+\frac{c}{4}\textrm{tr}^2(Q^2(t,x))\,dx
\end{equation}
      If $d=2,3$ and $(Q,u)$ is a smooth solution of (\ref{system}) such that $Q\in L^\infty(0,T; H^1(\mathbb{R}^d))\cap L^2(0,T;H^2(\mathbb{R}^d))$  and $u\in L^\infty(0,T;L^2(\mathbb{R}^d))\cap L^2(0,T;H^1(\mathbb{R}^d))$ then, for all $t<T$, we have:
\begin{equation}
\frac{d}{dt}E(t)=-\nu\int_{\mathbb{R}^d}|\nabla u|^2\,dx-\Gamma\int_{\mathbb{R}^d} \textrm{tr}\left(L\Delta Q-aQ+b[Q^2-\frac{\textrm{tr}(Q^2)}{d}Id]-cQ\textrm{tr}(Q^2)\right)^2\,dx\le 0
\label{energydecay}
\end{equation}
\label{prop:Lyapunov}
\end{proposition}

\smallskip     {\bf Proof.} We multiply the first equation in (\ref{system}) to the right by $-H$, take the trace, integrate over $\mathbb{R}^d$ and by parts and sum with the second equation multiplied by $u$ and integrated over $\mathbb{R}^d$ and by parts (let us observe that because of our assumptions on $Q$ and $u$ we do not have boundary terms, when integrating by parts). We obtain:

\begin{eqnarray}
\frac{d}{dt}\int_{\mathbb{R}^d}\frac{1}{2}|u|^2+\frac{L}{2}|\nabla Q|^2+\frac{a}{2}\textrm{tr}(Q^2)-\frac{b}{3}\textrm{tr}(Q^3)+\frac{c}{4}\textrm{tr}^2(Q^2)\,dx\nonumber\\+\nu\int_{\mathbb{R}^d}|\nabla u|^2\,dx+\Gamma\int_{\mathbb{R}^d}\textrm{tr}\left(L\Delta Q-aQ+b[Q^2-\frac{\textrm{tr}(Q^2)}{d}Id]-cQ\textrm{tr}(Q^2)\right)^2\,dx\nonumber\\=\underbrace{\int_{\mathbb{R}^d}u\cdot\nabla Q_{\alpha\beta}\left(-aQ_{\alpha\beta}+b[Q_{\alpha\gamma}Q_{\gamma\beta}-\frac{\delta_{\alpha\beta}}{d}\textrm{tr}(Q^2)]-cQ_{\alpha\beta}\textrm{tr}(Q^2))\right)\,dx}_{\stackrel{\rm{def}}{=}\mathcal{I}}\nonumber\\+\underbrace{\int_{\mathbb{R}^d}\left(-\Omega_{\alpha\gamma} Q_{\gamma\beta}+Q_{\alpha\gamma}\Omega_{\gamma\beta}\right)\left(-aQ_{\alpha\beta}+b[Q_{\alpha\delta}Q_{\delta\beta}-\frac{\delta_{\alpha\beta}}{d}\textrm{tr}(Q^2)]-cQ_{\alpha\beta}\textrm{tr}(Q^2))\right)\,dx}_{\stackrel{\rm{def}}{=}\mathcal{II}}
\nonumber
\end{eqnarray}
\begin{eqnarray}-\xi\underbrace{\int_{\mathbb{R}^d} \big(Q_{\alpha\gamma}+\frac{\delta_{\alpha\gamma}}{d}\big)D_{\gamma\beta}H_{\alpha\beta}\,dx}_{\stackrel{\rm{def}}{=}\mathcal{J}_1}-\xi\underbrace{\int_{\mathbb{R}^d} D_{\alpha\gamma}\left(Q_{\gamma\beta}+\frac{\delta_{\gamma\beta}}{d}\right)H_{\alpha\beta}\,dx}_{\stackrel{\rm{def}}{=}\mathcal{J}_2}
+2\xi\underbrace{\int_{\mathbb{R}^d}\left(Q_{\alpha\beta}+\frac{\delta_{\alpha\beta}}{d}\right)H_{\alpha\beta}\textrm{tr}(Q\nabla u)\,dx}_{\stackrel{\rm{def}}{=}\mathcal{J}_3}\nonumber\\
+L\underbrace{\int_{\mathbb{R}^d}u_{\gamma} Q_{\alpha\beta,\gamma}\Delta Q_{\alpha\beta}\,dx}_{\stackrel{\rm{def}}{=}\mathcal{A}}\underbrace{
-\frac{L}{2}\int_{\mathbb{R}^d} u_{\alpha,\gamma}Q_{\gamma\beta}\Delta Q_{\alpha\beta}\,dx}_{\stackrel{\rm{def}}{=}\mathcal{B}}
\nonumber
\end{eqnarray}
\begin{eqnarray}
\underbrace{+\frac{L}{2}\int_{\mathbb{R}^d}u_{\gamma,\alpha}Q_{\gamma\beta}\Delta Q_{\alpha\beta}}_{\stackrel{\rm{def}}{=}\mathcal{C}}\,dx
\underbrace{+\frac{L}{2}\int_{\mathbb{R}^d}Q_{\alpha\gamma}u_{\gamma,\beta}\Delta Q_{\alpha\beta}\,dx}_{\mathcal{C}}\underbrace{-\frac{L}{2}\int_{\mathbb{R}^d}Q_{\alpha\gamma}u_{\beta,\gamma}\Delta Q_{\alpha\beta}\,dx}_{\mathcal{B}}\nonumber\\
+L\underbrace{\int_{\mathbb{R}^d}Q_{\gamma\delta,\alpha}Q_{\gamma\delta,\beta}u_{\alpha,\beta}\,dx}_{\stackrel{\rm{def}}{=}\mathcal{AA}}\underbrace{-L\int_{\mathbb{R}^d} Q_{\alpha\gamma}\Delta Q_{\gamma\beta}u_{\alpha,\beta}\,dx}_{\stackrel{\rm{def}}{=}\mathcal{CC}}\underbrace{+L\int_{\mathbb{R}^d}\Delta Q_{\alpha\gamma}Q_{\gamma\beta}u_{\alpha,\beta}\,dx}_{\stackrel{\rm{def}}{=}\mathcal{BB}}\nonumber\\
+\xi\underbrace{\int_{\mathbb{R}^d} \big(Q_{\alpha\gamma}+\frac{\delta_{\alpha\gamma}}{d}\big)H_{\gamma\beta}u_{\alpha,\beta}\,dx}_{\stackrel{\rm{def}}{=}\mathcal{JJ}_1}+\xi\underbrace{\int_{\mathbb{R}^d}H_{\alpha\gamma}\big(Q_{\gamma\beta}+\frac{\delta_{\gamma\beta}}{d}\big)u_{\alpha,\beta}\,dx}_{\stackrel{\rm{def}}{=}\mathcal{JJ}_2}-2\xi\underbrace{\int_{\mathbb{R}^d} \big(Q_{\alpha\beta}+\frac{\delta_{\alpha\beta}}{d}\big)u_{\alpha,\beta}\textrm{tr}(QH)\,dx}_{\stackrel{\rm{def}}{=}\mathcal{JJ}_3}\nonumber
\end{eqnarray}
\begin{eqnarray}
=\underbrace{-L\int_{\mathbb{R}^d} u_{\alpha,\gamma}Q_{\gamma\beta}\Delta Q_{\alpha\beta}\,dx}_{2\mathcal{B}}\underbrace{+L\int_{\mathbb{R}^d}u_{\gamma,\alpha}Q_{\gamma\beta}\Delta Q_{\alpha\beta}\,dx}_{2\mathcal{C}} \underbrace{-L\int_{\mathbb{R}^d} Q_{\alpha\gamma}\Delta Q_{\gamma\beta}u_{\alpha,\beta}\,dx}_{\mathcal{CC}}\underbrace{+L\int_{\mathbb{R}^d}\Delta Q_{\alpha\gamma}Q_{\gamma\beta}u_{\alpha,\beta}\,dx}_{\mathcal{BB}}=0
\label{Lyapunovcancellation}
\end{eqnarray} where $\mathcal{I}=0$ (since $\nabla\cdot u=0$), $\mathcal{II}=0$ (since $Q_{\alpha\beta}=Q_{\beta\alpha}$) and for the second equality we used

\begin{eqnarray}\underbrace{\int_{\mathbb{R}^d}u_\gamma Q_{\alpha\beta,\gamma}\Delta Q_{\alpha\beta}\,dx}_{\mathcal{A}}\underbrace{+\int_{\mathbb{R}^d}Q_{\gamma\delta,\alpha}Q_{\gamma\delta,\beta}u_{\alpha,\beta}\,dx}_{\mathcal{AA}}=\int_{\mathbb{R}^d}u_\gamma Q_{\alpha\beta,\gamma}\Delta Q_{\alpha\beta}\,dx\nonumber\\-\int_{\mathbb{R}^d}Q_{\gamma\delta,\alpha}Q_{\gamma\delta,\beta\beta}u_\alpha\,dx-
\int_{\mathbb{R}^d}Q_{\gamma\delta,\alpha\beta}Q_{\gamma\delta,\beta}u_\alpha\,dx=\int_{\mathbb{R}^d}\frac{1}{2}Q_{\gamma\delta,\beta}Q_{\gamma\delta,\beta}u_{\alpha,\alpha}\,dx=0\nonumber
\end{eqnarray} together with $Q_{\alpha\alpha}=H_{\alpha\alpha}=u_{\alpha,\alpha}=0$, $\mathcal{J}_3=\mathcal{JJ}_3$ and

\begin{eqnarray}
\mathcal{J}_1+\mathcal{J}_2= \int_{\mathbb{R}^d} \frac{1}{2}Q_{\alpha\gamma}u_{\gamma,\beta}H_{\alpha\beta}+\frac{1}{2}Q_{\alpha\gamma}u_{\beta,\gamma}H_{\alpha\beta}+\frac{1}{2}u_{\alpha,\gamma}Q_{\gamma\beta}H_{\alpha\beta}+\frac{1}{2}u_{\gamma,\alpha}Q_{\gamma\beta}H_{\alpha\beta}\,dx\nonumber\\+\frac{2}{d}\int_{\mathbb{R}^d}D_{\alpha\beta}H_{\alpha\beta}
=\int_{\mathbb{R}^d}\frac{1}{2}\big(Q_{\alpha\gamma}u_{\gamma,\beta}H_{\alpha\beta}+u_{\gamma,\alpha}Q_{\gamma\beta}H_{\alpha\beta}\big)+\frac{1}{2}\big(Q_{\alpha\gamma}u_{\beta,\gamma}H_{\alpha\beta}+u_{\alpha,\gamma}Q_{\gamma\beta}H_{\alpha\beta}\big)\,dx\nonumber\\+\frac{1}{d}\int_{\mathbb{R}^d}(u_{\alpha,\beta}+u_{\beta,\alpha})H_{\alpha\beta}\,dx
=\int_{\mathbb{R}^d} H_{\beta\alpha}Q_{\alpha\gamma}u_{\gamma,\beta}+ Q_{\gamma\alpha}H_{\alpha\beta}u_{\beta,\gamma}\,dx+\frac{2}{d}\int_{\mathbb{R}^d} u_{\alpha,\beta}H_{\alpha\beta}\,dx=\mathcal{JJ}_1+\mathcal{JJ}_2\nonumber
\end{eqnarray} 
\par Finally,  the last equality in (\ref{Lyapunovcancellation}) is a consequence of the straightforward identities $2\mathcal{B}+\mathcal{BB}=2\mathcal{C}+\mathcal{CC}=0.$ 
$\Box$

\bigskip
      In the following we assume that there exists a smooth solution of (\ref{system}) and  obtain estimates on the behaviour of various norms:

\begin{proposition} Let $(Q,u)$ be a smooth solution of (\ref{system}) in dimension $d=2$ or $d=3$, with restriction (\ref{c+}), and smooth initial data $(\bar Q(x),\bar u(x))$, that decays fast enough at infinity so that we can integrate by parts in space (for any $t\ge 0$) without boundary terms.
 \par     (i) For $(\bar Q,\bar u)\in H^1\times L^2$,we have
\begin{equation}
\|Q(t,\cdot)\|_{H^1}\le C_1+\bar C_1 e^{\bar C_1t}\|\bar Q\|_{H^1}, \forall t\ge 0
\label{apriorih1}
\end{equation} with $C_1,\bar C_1$ depending on $(a,b,c,d,\Gamma,L, \nu,\bar Q,\bar u)$.
 \par     (ii) For  $(\bar Q,\bar u)\in H^1\times L^2$, we have:
\begin{eqnarray}
\|u(t,\cdot)\|_{L^2}^2+2\nu\int_0^t\|\nabla u(s,\cdot)\|_{L^2}^2\,ds+L\|\nabla Q(t,\cdot)\|_{L^2}^2+\Gamma L^2\int_0^t \|\Delta Q(s,\cdot)\|_{L^2}^2\,ds\le C_2+\bar C_2e^{\bar C_2 t}
\end{eqnarray} with the constants $C_2,\bar C_2$ depending on $(a,b,c,d,L,\Gamma,\bar u,\bar Q, \nu)$.
\label{prop:aprioriest}
\end{proposition}

\smallskip     {\bf Proof.} We multiply the first equation in (\ref{system}) by $Q$, take the trace, integrate over $\mathbb{R}^d$ and by parts and we obtain:

\begin{eqnarray}
\frac{1}{2} \frac{d}{dt}\int_{\mathbb{R}^d} |Q|^2(t,x)\,dx=\Gamma\Big(-L\int_{\mathbb{R}^d} |\nabla Q|^2\,dx-a\int_{\mathbb{R}^d}|Q(x)|^2\,dx+b\int_{\mathbb{R}^d}\textrm{tr}(Q^3)\,dx-c\int_{\mathbb{R}^d}|Q|^4\,dx\Big)\nonumber\\
+\underbrace{\int_{\mathbb{R}^d}\textrm{tr}(\Omega Q^2-Q\Omega Q)\,dx}_{\stackrel{\rm{def}}{=}\mathcal{I}}\nonumber\\
+\underbrace{\xi\int_{\mathbb{R}^d} D_{\alpha\gamma}(Q_{\gamma\beta}+\frac{\delta_{\gamma\beta}}{d})Q_{\alpha\beta}+(Q_{\alpha\gamma}+\frac{\delta_{\alpha\gamma}}{d})D_{\gamma\beta}Q_{\alpha\beta}-2(Q_{\alpha\beta}+\frac{\delta_{\alpha\beta}}{d})Q_{\alpha\beta}\textrm{tr}(Q\nabla u)\,dx}_{\stackrel{\rm{def}}{=}\mathcal{II}}\nonumber
\end{eqnarray}

      Recalling that $Q$ is symmetric we have $\mathcal{I}=0$. Also:

\begin{displaymath}
|\mathcal{II}|=|2\xi||\int_{\mathbb{R}^d} \frac{1}{d}D_{\alpha\beta}Q_{\alpha\beta}+D_{\alpha\gamma}Q_{\gamma\beta}Q_{\beta\alpha}-Q_{\alpha\beta}Q_{\alpha\beta}\textrm{tr}(Q\nabla u)\,dx|\le C(\xi,d)\int_{\mathbb{R}^d}\varepsilon |\nabla u|^2+\frac{1}{\varepsilon}(|Q|^2+|Q|^6)\,dx
\end{displaymath} hence we get:

\begin{equation}
   \frac{d}{dt} \int_{\mathbb{R}^d}|Q|^2\,dx\le C(\xi,d)\varepsilon\int_{\mathbb{R}^d} |\nabla u|^2\,dx+C(\xi,\Gamma,L,a,b,c,d)\frac{1}{\varepsilon}\int_{\mathbb{R}^d}|Q|^2+|Q|^6\,dx\nonumber 
\end{equation}

      Adding the last relation multiplied by $A^2$ (with $A\in\mathbb{R}$ a constant to be chosen) and (\ref{energydecay}) we get

\begin{eqnarray}
   \frac{d}{dt}\int_{\mathbb{R}^d} \frac{L}{2}|\nabla Q|^2+A^2|Q|^2+\frac{a}{2}|Q|^2-\frac{b}{3}\textrm{tr}(Q^3)+\frac{c}{4}|Q|^4\,dx\le (C(\xi,d)\varepsilon A^2-\nu)\int_{\mathbb{R}^d}|\nabla u|^2\nonumber\\
   +\frac{C(\xi,\Gamma,L,a,b,c,d)A^2}{\varepsilon} \int_{\mathbb{R}^d}|Q|^2+|Q|^6\,dx
\label{lowerbdenergy1}
\end{eqnarray}

      Let us observe that for $Q$  a traceless, symmetric, $3\times 3$ matrix we have:

\begin{equation}
\textrm{tr}(Q^3)\le \frac{3\delta}{8}\textrm{tr}^2(Q^2)+\frac{1}{\delta} \textrm{tr}(Q^2),\forall \delta>0
\label{est:Q3}
\end{equation}

      Indeed, if $Q$ has the eigenvalues $x,y,-x-y$ then $\textrm{tr}(Q^3)=-3xy(x+y)$, $\textrm{tr}(Q^2)=2(x^2+y^2+xy)$ and the inequality (\ref{est:Q3}) follows. Then, choosing $\delta$ appropriately small and  for $A$ large enough we have:

\begin{equation}
\int_{\mathbb{R}^d}\frac{L}{2}|\nabla Q(t,x)|^2+\frac{A^2}{2}|Q(t,x)|^2\,dx\le \int_{\mathbb{R}^d}\frac{L}{2}|\nabla Q(t,x)|^2+A^2|Q(t,x)|^2+\frac{a}{2}|Q(t,x)|^2-\frac{b}{3}\textrm{tr}(Q^3(t,x))+\frac{c}{4}|Q(t,x)|^4\,dx
\label{lowerbdenergy2}
\end{equation} (note that we can choose $\delta>0$ appropriately small so that we have the previous inequality precisely because  of our assumption  (\ref{c+}), namely $c>0$)
\par In the case $d=2$ we have $\textrm{tr}(Q^3)=0$ (as $Q$ is traceless and symmetric) but we still need the assumption $c>0$ in order to have the estimate (\ref{lowerbdenergy2}).

  \par    Using together (\ref{lowerbdenergy1}) and (\ref{lowerbdenergy2}) and choosing $\varepsilon>0$ appropriately small  so that 
$C(\xi,d)\varepsilon A^2-\nu<0$ we get:

\begin{eqnarray}
\int_{\mathbb{R}^d} \frac{L}{2}|\nabla Q(t,x)|^2+\frac{A^2}{2}|Q(t,x)|^2\,dx\le C_1+C_2\int_0^t\int_{\mathbb{R}^d} |Q(s,x)|^2+|Q(s,x)|^6\,ds dx\nonumber\\
\le C_1+C_3\int_0^t\int_{\mathbb{R}^d}  \frac{L}{2}|\nabla Q(s,x)|^2+\frac{A^2}{2}|Q(s,x)|^2\,ds\,dx
\end{eqnarray} where $C_1$ depends on the initial data $A,\varepsilon, L,a,b,c,d$ and $C_2,C_3$ depend on $A,\varepsilon, L,a,b,c,d,\Gamma$.  Thus we obtain the claimed estimate (\ref{apriorih1}).

\smallskip      (ii) Relation (\ref{energydecay}) implies

\begin{eqnarray}
 \frac{L}{2}\|\nabla Q(t,\cdot)\|_{L^2}^2+\frac{1}{2}\|u(t,\cdot)\|_{L^2}^2+\nu\int_0^t \|\nabla u(s,\cdot)\|_{L^2}^2\,ds+\Gamma L^2 \int_0^t \|\Delta Q(s,\cdot)\|_{L^2}^2ds\nonumber\\
 \le C\int_{\mathbb{R}^d} \textrm{tr}(Q^2(t,x))+\textrm{tr}^2(Q^2(t,x))\,dx+C\int_{\mathbb{R}^d} \textrm{tr}(Q^2(0,x))+\textrm{tr}^2(Q^2(0,x))\,dx+ \frac{L}{2}\|\nabla Q(0,\cdot)\|_{L^2}+\frac{1}{2}\|u(0,\cdot)\|_{L^2}^2\nonumber\\
 - \Gamma \int_0^t\int_{\mathbb{R}^d} \textrm{tr}\Big(L\Delta Q\big(aQ-bQ^2+cQ\textrm{tr}(Q^2)\big)\Big)\,dx\,ds-\Gamma \int_0^t\int_{\mathbb{R}^d} \textrm{tr}\Big(\big(aQ-bQ^2+cQ\textrm{tr}(Q^2)\big)L\Delta Q\Big)\,dx\,ds\nonumber\\
 +\Gamma\int_0^t\int_{\mathbb{R}^d}\textrm{tr}\Big(aQ-bQ^2+cQ\textrm{tr}(Q^2)\Big)^2\,dx\,ds
 \end{eqnarray}

      In the last inequality we  use Holder inequality to estimate $\Delta Q$ in $L^2$  and absorb it in the left hand side  while the terms without gradients are estimated using (\ref{apriorih1}) and interpolation between the $L^2$ and $L^6$ norms.$\Box$

\section{Weak solutions}

      A pair $(Q,u)$ is called a weak solution  of the system (\ref{system}), subject to initial data
\begin{equation}
Q(0,x)=\bar Q(x)\in H^1(\mathbb{R}^d),\, u(0,x)=\bar u(x)\in L^2(\mathbb{R}^d), \nabla\cdot \bar u=0\,\textrm{ in }\mathcal{D}'(\mathbb{R}^d)
\label{initialdata}
\end{equation}
if   $Q\in L^\infty_{loc}(\mathbb{R}_+;H^1)\cap L^2_{loc}(\mathbb{R}_+;H^2)$, $u\in L^\infty_{loc}(\mathbb{R}_+;L^2)\cap L^2_{loc}(\mathbb{R}_+;H^1)$  and for every compactly supported $\varphi\in C^\infty([0,\infty)\times \mathbb{R}^d; S_0^{(d)})$, $\psi\in C^\infty ([0,\infty)\times\mathbb{R}^d;\mathbb{R}^d)$ with $\nabla\cdot\psi=0$ we have
\begin{eqnarray}\int_0^\infty\int_{\mathbb{R}^d}(-Q\cdot    \partial_t\varphi-\Gamma L\Delta Q\cdot \varphi)- Q\cdot u\nabla_x \varphi dx\,dt\nonumber\\-\int_0^\infty\int_{\mathbb{R}^d}(\xi D+\Omega)(Q+\frac{1}{d}Id)\cdot\varphi+(Q+\frac{1}{d}Id)(\xi D-\Omega)\cdot\varphi-2\xi (Q+\frac{1}{d}Id)\textrm{tr}(Q\nabla u)\cdot\varphi\,dx\,dt\nonumber\\
=\int_{\mathbb{R}^d}\bar Q(x)\cdot \varphi(0,x)\,dx+\Gamma\int_0^\infty\int_{\mathbb{R}^d}\Big\{-aQ+b[Q^2-\frac{\textrm{tr}(Q^2)}{d}Id]-cQ\textrm{tr}(Q^2)\Big\}\cdot\varphi\,\,dx\,dt
\label{weaksol1}
\end{eqnarray} 
and
\begin{eqnarray}\int_0^\infty\int_{\mathbb{R}^d}-u   \partial_t\psi-u_\alpha u_\beta   \partial_\alpha\psi_\beta+\nu\nabla u\nabla\psi\,\, dt\,dx-\int_{\mathbb{R}^d}\bar u(x)\psi(0,x)\,dx\nonumber\\
=L\int_0^\infty\int_{\mathbb{R}^d} 
 Q_{\gamma\delta,\alpha} Q_{\gamma\delta,\beta}\psi_{\alpha,\beta} -Q_{\alpha\gamma}\Delta Q_{\gamma\beta}\psi_{\alpha,\beta}+\Delta Q_{\alpha\gamma}Q_{\gamma\beta}\psi_{\alpha,\beta}\,\,dx\,dt\nonumber\\
+\xi\int_0^\infty\int_{\mathbb{R}^d}\left(Q_{\alpha\gamma}+\frac{\delta_{\alpha\gamma}}{d}\right)H_{\gamma\beta}\psi_{\alpha,\beta}+H_{\alpha\gamma}\left(Q_{\gamma\beta}+\frac{\delta_{\gamma\beta}}{d}\right)\psi_{\alpha,\beta}-2(Q_{\alpha\beta}+\frac{\delta_{\alpha\beta}}{d})Q_{\gamma\delta}H_{\gamma\delta}\psi_{\alpha,\beta}\,dx\,dt
\label{weaksol2}
\end{eqnarray} 
\begin{proposition} For $d=2,3$ there exists a weak solution $(Q,u)$ of the system (\ref{system}), with restriction (\ref{c+}), subject to initial conditions (\ref{initialdata}). The solution $(Q,u)$ is such that $Q\in L^\infty_{loc}(\mathbb{R}_+;H^1)\cap L^2_{loc}(\mathbb{R}_+;H^2)$ and $u\in L^\infty_{loc}(\mathbb{R}_+;L^2)\cap L^2_{loc}(\mathbb{R}_+;H^1)$.
\label{prop:weak}
\end{proposition} 

\smallskip     {\bf Proof.}  We define the mollifying operator

$$\widehat{J_nf}(\xi)=1_{[\frac{1}{n},n]}(|\xi|)\hat f(\xi)$$ and consider the system:

\begin{equation}
\left\{\begin{array}{l}
         \partial_t Q^{(n)}+J_n \Big(\mathcal{P}J_n u^n \nabla J_nQ^{(n)}\Big)   -J_n\Big((\xi \mathcal{P}J_n D^n+\mathcal{P}J_n\Omega^n)(J_n Q^{(n)}+\frac{1}{d}Id)\Big)+J_n\Big((J_nQ^{(n)}+\frac{1}{d}Id)(\xi\mathcal{P}J_n D^n-\mathcal{P}J_n\Omega^n)\Big)\\-2\xi J_n \Big((J_n Q^{(n)}+\frac{1}{d}Id)\textrm{tr} J_n\big(J_n Q^{(n)}\nabla\mathcal{P}J_n u^n\big)\Big)
=\Gamma L\Delta J_n Q^{(n)}+\\+\Gamma\Big(-aJ_n Q^{(n)}+b[J_n(J_nQ^{(n)} J_n Q^{(n)})-\frac{\textrm{tr}(J_n(J_nQ^{(n)}J_n Q^{(n)}))}{d}Id]-cJ_n\bigg(J_nQ^{(n)}\textrm{tr}(J_n(J_nQ^{(n)}J_n Q^{(n)}))\Big)\bigg) \\
         \partial_t u^n +\mathcal{P}J_n(\mathcal{P}J_nu^n\nabla\mathcal{P}J_n u^n)=-\xi\mathcal{P}J_n\nabla\cdot \Big(\left( J_n Q^{(n)}+\frac{1}{d}Id\right)J_n \tilde H^{(n)}\Big)-\xi\mathcal{P} J_n\nabla\cdot\Big( J_n \tilde H^{(n)}\left(J_n Q^{(n)}+\frac{1}{d}Id\right)\Big)\nonumber\\
 +2\xi \mathcal{P}J_n\nabla\cdot\Big(\big(J_nQ^{(n)}+\frac{1}{d}\big)J_n\big(J_nQ^{(n)}J_n\tilde H^{(n)}\big)\Big)- L\mathcal{P}J_n(\nabla\cdot\left(\textrm{tr}(\nabla J_n Q^{(n)}\nabla J_n Q^{(n)})-\frac{1}{d}|\nabla J_n Q^{(n)}|^2 Id\right))\nonumber\\+L\mathcal{P}(\nabla\cdot J_n\left(J_nQ^{(n)} \Delta J_nQ^{(n)}-\Delta J_n Q^{(n)} J_n Q^{(n)}\right))+\nu\Delta \mathcal{P}J_nu^n\\
     \end{array}\right.
      \label{approxsystem+}
\end{equation} where $\mathcal{P}$ denotes the Leray projector onto divergence-free vector fields and $\tilde H^{(n)}\stackrel{def}{=}LJ_n\Delta Q^{(n)} -aJ_n Q^{(n)}+b[J_n(J_nQ^{(n)} J_n Q^{(n)})-\frac{\textrm{tr}(J_n(J_nQ^{(n)}J_n Q^{(n)}))}{d}Id]-cJ_n\bigg(J_nQ^{(n)}\textrm{tr}(J_n(J_nQ^{(n)}J_n Q^{(n)}))\bigg)$.

      The system above can be regarded as an ordinary differential equation in $L^2$ verifying the conditions of the Cauchy-Lipschitz theorem. Thus it admits a unique maximal solution $(Q^{(n)}, u^n)\in C^1([0,T_n); L^2(\mathbb{R}^d;\mathbb{R}^{d\times d})\times L^2(\mathbb{R}^d,\mathbb{R}^d))$. As we have $(\mathcal{P}J_n)^2=\mathcal{P}J_n$ and $J_n^2=J_n$ the pair $(J_n Q^{(n)}, \mathcal{P}J_n u^n)$ is also a solution of (\ref{approxsystem+}). By uniqueness we have $(J_n Q^{(n)}, \mathcal{P}J_n u^n)=(Q^{(n)}, u^n)$ hence $(Q^{(n)}, u^n)\in C^1([0,T_n),H^\infty)$ and $(Q^{(n)}, u^n)$ satisfy the system:

\begin{equation}
\left\{\begin{array}{l}
         \partial_t Q^{(n)}+J_n\big(u^n \nabla Q^{(n)}\big)-J_n\Big((\xi D^n+\Omega^n)(Q^{(n)}+\frac{1}{d}Id)\Big)+J_n\Big((Q^{(n)}+\frac{1}{d}Id)(\xi D^n-\Omega^n)\Big)\\-2\xi J_n \Big((Q^{(n)}+\frac{1}{d}Id)\textrm{tr} J_n\big(Q^{(n)}\nabla u^n\big)\Big)=\Gamma L\Delta Q^{(n)}\\
      +\Gamma\Big(-a Q^{(n)}+b[J_n(Q^{(n)}  Q^{(n)})-\frac{\textrm{tr}(J_n(Q^{(n)} Q^{(n)}))}{d}Id]-cJ_n\bigg(Q^{(n)}\textrm{tr}(J_n(Q^{(n)} Q^{(n)}))\Big)\bigg) \\
         \partial_t u^n +\mathcal{P}J_n(u^n\nabla u^n)=-\xi\mathcal{P}J_n\nabla\cdot \Big(\left(Q^{(n)}+\frac{1}{d}Id\right) \bar H^{(n)}\Big)-\xi\mathcal{P} J_n\nabla\cdot\Big( \bar H^{(n)}\left(Q^{(n)}+\frac{1}{d}Id\right)\Big)\\
 +2\xi\mathcal{P} J_n\nabla\cdot\Big(\big(Q^{(n)}+\frac{1}{d}\big)J_n\big(Q^{(n)}\bar  H^{(n)}\big)\Big)- L\mathcal{P}J_n(\nabla\cdot\left(\textrm{tr}(\nabla  Q^{(n)}\nabla  Q^{(n)})-\frac{1}{d}|\nabla Q^{(n)}|^2 Id\right))\\+L\mathcal{P}(\nabla\cdot J_n\left(Q^{(n)} \Delta Q^{(n)}-\Delta Q^{(n)} Q^{(n)} \right))+\nu\Delta u^n
      \end{array}\right.
      \label{approxsystem++}
\end{equation}  where $\bar H^{(n)}\stackrel{def}{=}L\Delta Q^{(n)} -a Q^{(n)}+b[J_n(Q^{(n)}Q^{(n)})-\frac{\textrm{tr}(J_n(Q^{(n)} Q^{(n)}))}{d}Id]-cJ_n\bigg(Q^{(n)}\textrm{tr}(J_n(Q^{(n)} Q^{(n)}))\bigg)$.

      We can argue as in the proof of  the apriori estimates and the same estimates hold for the approximating system (\ref{approxsystem++}). These estimates allow us to conclude that $T_n=\infty$ and we also get the following apriori bounds:

\begin{eqnarray}
\sup_n \|Q^{(n)}\|_{L^2(0,T;H^2)\cap L^\infty(0,T;H^1)}<\infty\nonumber\\
\sup_n \|u^n\|_{L^\infty(0,T;L^2)\cap L^2(0,T;H^1)}<\infty
\label{weaksolapriori}
\end{eqnarray} for any $T<\infty$.

      The pair $(Q^{(n)}, u^n)$ is also a weak solution of the approximating system (\ref{approxsystem++}) hence for every compactly supported $\varphi\in C^\infty([0,\infty)\times \mathbb{R}^d; S_0^{(d)})$, $\psi\in C^\infty ([0,\infty)+\times\mathbb{R}^d;\mathbb{R}^d)$ with $\nabla\cdot\psi=0$ we have:

\begin{eqnarray}\int_0^\infty\int_{\mathbb{R}^d}(-Q^{(n)}\cdot    \partial_t\varphi-\Gamma L \Delta Q^{(n)}\cdot\varphi) - J_n\big(Q^{(n)}\cdot u^n\big)\nabla_x\varphi-J_n\Big((\xi D^n+\Omega^n)(Q^{(n)}+\frac{1}{d}Id)\Big)\cdot\varphi \,\,dx\,dt\nonumber\\
-\int_0^\infty\int_{\mathbb{R}^d}J_n\Big((Q^{(n)}+\frac{1}{d}Id)(\xi D^n-\Omega^n)\Big)\cdot\varphi-2\xi J_n\bigg((Q^{(n)}+\frac{1}{d}Id)\textrm{tr}J_n\big(Q^{(n)}\nabla u^n\big)\bigg)\cdot\varphi\,dx\,dt\nonumber\\
=\int_{\mathbb{R}^d}\bar Q(x)\cdot \varphi(0,x)\,dx+\Gamma\int_0^\infty\int_{\mathbb{R}^d}\{-aQ^{(n)}+b[J_n\left(Q^{(n)}\right)^2-\frac{\textrm{tr}\big(J_n\left(Q^{(n)}\right)^2\big)}{d}Id]-cJ_n\bigg(Q^{(n)}\textrm{tr}(J_n(Q^{(n)})^2)\}\bigg)\cdot\varphi\,\,dx\,dt
\label{weaksol1+}
\end{eqnarray}
and
\begin{eqnarray}\int_0^\infty\int_{\mathbb{R}^d}-u^n   \partial_t\psi-J_n(u^n_\alpha u^n_\beta)   \partial_\alpha\psi_\beta+\nu\nabla u^n\nabla\psi\,\,dx\, dt-\int_{\mathbb{R}^d}\bar u(x)\psi(0,x)\,dx\nonumber\\
=L\int_0^\infty\int_{\mathbb{R}^d}\Big\{ J_n\left(Q^{(n)}_{\gamma\delta,\alpha} Q^{(n)}_{\gamma\delta,\beta}\right)\psi_{\alpha,\beta}-J_n\left(Q^{(n)}_{\alpha\gamma}\Delta Q^{(n)}_{\gamma\beta}-\nu\Delta Q^{(n)}_{\alpha\gamma}Q^{(n)}_{\gamma\beta}\right)\psi_{\alpha,\beta}\Big\}\,\,dx\,dt\nonumber\\
+\xi\int_0^\infty\int_{\mathbb{R}^d}\Big\{J_n\Big(\left(Q_{\alpha\gamma}^{(n)}+\frac{\delta_{\alpha\gamma}}{d}\right)\bar H_{\gamma\beta}^{(n)}\Big)\psi_{\alpha,\beta}+J_n\Big(\bar H_{\alpha\gamma}^{(n)}\left(Q_{\gamma\beta}^{(n)}+\frac{\delta_{\gamma\beta}}{d}\right)\Big)\psi_{\alpha,\beta}\Big\}\,dx\,dt\nonumber\\
-2\xi\int_0^\infty\int_{\mathbb{R}^d}\Big\{J_n\Big(\left(Q_{\alpha\beta}^{(n)}+\frac{\delta_{\alpha\beta}}{d}\right)J_n\big(Q_{\gamma\delta}^{(n)}\bar H_{\gamma\delta}^{(n)}\big)\Big)\psi_{\alpha,\beta}\Big\}\,dx\,dt
\label{weaksol2+}
\end{eqnarray}

      We consider the solutions of (\ref{approxsystem++}) and taking into account the bounds (\ref{weaksolapriori}) we get, by classical compactness and weak convergence arguments, that there exists a $Q\in L^\infty_{loc}(\mathbb{R}_+;H^1)\cap L^2_{loc}(\mathbb{R}_+;H^2)$ and a $u\in L^\infty_{loc}(\mathbb{R}_+;L^2)\cap L^2_{loc}(\mathbb{R}_+;H^1)$ so that, on a subsequence, we have:

\begin{eqnarray}
Q^{(n)}\rightharpoonup Q\textrm{ in } L^2(0,T;H^2)\,\textrm{ and }Q^{(n)}\to Q\textrm{ in }L^2(0,T;H_{loc}^{2-\varepsilon}),\forall \varepsilon>0\nonumber\\
Q^{(n)}(t)\rightharpoonup Q(t)\textrm{ in }H^1\textrm{ for all }t\in\mathbb{R}_+\nonumber\\
u^n\rightharpoonup u\textrm{ in }L^2(0,T;H^1)\,\textrm{ and }u^n\to u\textrm{ in }L^2(0,T;H_{loc}^{1-\varepsilon}),\forall\varepsilon>0\nonumber\\
u^n(t)\rightharpoonup u(t)\,\textrm{ in }L^2\textrm{ for all }t\in\mathbb{R}_+
\label{convergences}
\end{eqnarray}

 These  convergences allow us to the pass to the limit in the weak solutions 
(\ref{weaksol1+}),(\ref{weaksol2+}) to obtain a weak solution of (\ref{system}), namely  (\ref{weaksol1}),(\ref{weaksol2}). Of all the terms there are only two types of terms that are slightly difficult to treat in passing  to the limit. A first type  is  a term in (\ref{weaksol2+}), namely 

$$L\int_0^\infty\int_{\mathbb{R}^d}J_n\left(Q^{(n)}_{\alpha\gamma}\Delta Q^{(n)}_{\gamma\beta}-\Delta Q^{(n)}_{\alpha\gamma}Q^{(n)}_{\gamma\beta}\right)\psi_{\alpha,\beta}\,\,dx\,dt=L\int_0^\infty\int_{\mathbb{R}^d}\left(Q^{(n)}_{\alpha\gamma}\Delta Q^{(n)}_{\gamma\beta}-\Delta Q^{(n)}_{\alpha\gamma}Q^{(n)}_{\gamma\beta}\right)\cdot J_n\psi_{\alpha,\beta}\,\,dx\,dt.$$

 Recalling that $\psi$ is compactly supported we have that there exists a time $T>0$ so that $\psi(t,x)=J_n\psi(t,x)=0,\forall t>T, x\in\mathbb{R}^d,n\in\mathbb{N}$. Taking into account that $\psi$ is compactly supported and the convergences (\ref{convergences}) one can easily pass to the limit the terms $   \partial_\beta J_n\psi_\alpha Q^{(n)}_{\alpha\gamma}$ and $   \partial_\beta J_n\psi_\alpha Q^{(n)}_{\gamma\beta}$ strongly in $L^2(0,T;L^2)$. Indeed we have:

\begin{equation}
   \partial_\beta J_n\psi_\alpha Q^{(n)}_{\alpha\gamma}-   \partial_\beta\psi_\alpha Q_{\alpha\gamma}=
\underbrace{\Big(   \partial_\beta J_n\psi_\alpha-   \partial_\beta\psi_\alpha\Big) Q^{(n)}_{\alpha\gamma}}_{\mathcal{I}}+
\underbrace{   \partial_\beta\psi_\alpha\Big(Q^{(n)}_{\alpha\gamma}-Q_{\alpha\gamma}\Big)}_{\mathcal{II}}
\label{lim:decomptwo}
\end{equation} and the first term, $\mathcal{I}$, converges to $0$, strongly in $L^2(0,T;L^2)$ because $\psi$ is smooth and compactly supported, hence $   \partial_\beta J_n\psi-   \partial_\beta\psi$ converges to zero in any $L^q(0,T;L^p)$  and $Q^{(n)}$ is bounded in $L^\infty$ in time and $L^p$ in space ($1<p<\infty$ if $d=2$ and $2\le p\le 6$ if $d=3$, due to the bounds (\ref{weaksolapriori})). On the other hand  the second term $\mathcal{II}$ converges strongly to zero in $L^2(0,T;L^2)$ because of (\ref{convergences}) and the fact that $\psi$ is compactly supported.
 
       Relations (\ref{convergences}) give that   $\Delta Q^{(n)}_{\gamma\beta}$, $\Delta Q^{(n)}_{\alpha\gamma}$ converges weakly in $L^2(0,T;L^2)$. Thus  we get convergence to the limit term 
 
 \begin{eqnarray}
 L\int_0^\infty \int_{\mathbb{R}^d} (\Delta Q_{\gamma\beta})(   \partial_\beta \psi_\alpha Q_{\alpha\gamma})dxdt-L\int_0^\infty\int_{\mathbb{R}^d}(\Delta Q_{\alpha\gamma})(   \partial_\beta\psi_\alpha Q_{\gamma\beta})dxdt\nonumber\\
=L\int_0^T \int_{\mathbb{R}^d} (\Delta Q_{\gamma\beta})(   \partial_\beta \psi_\alpha Q_{\alpha\gamma})dxdt-L\int_0^T\int_{\mathbb{R}^d}(\Delta Q_{\alpha\gamma})(   \partial_\beta\psi_\alpha Q_{\gamma\beta})dxdt.
 \end{eqnarray}
 
       Another type of term that could cause difficulties in passing to the limit is a part of the term in last line of (\ref{weaksol2+}) namely
 
 \begin{eqnarray}
 \int_0^\infty\int_{\mathbb{R}^d}\Big\{J_n\Big(Q_{\alpha\beta}^{(n)}J_n\big(Q_{\gamma\delta}^{(n)}\Delta Q_{\gamma\delta}^{(n)}\big)\Big)\psi_{\alpha,\beta}\Big\}\,dx\,dt=\int_0^\infty\int_{\mathbb{R}^d}Q_{\gamma\delta}^{(n)}\Delta Q_{\gamma\delta}^{(n)} J_n\Big(J_n\psi_{\alpha,\beta}Q_{\alpha\beta}^{(n)}\Big)\,dx\,dt
 \end{eqnarray}

       In order to treat this term we claim first that 
 
 \begin{equation}
 \|Q^{(n)}_{\gamma\delta} J_n\Big(J_n\psi_{\alpha,\beta}Q^{(n)}_{\alpha\beta}\Big)-Q^{(n)}_{\gamma\delta}J_n\psi_{\alpha,\beta}Q^{(n)}_{\alpha\beta}\|_{L^2(0,T;L^2)}\to 0
 \label{lim:simplu}
 \end{equation}
 
       Indeed we have 
 
 \begin{eqnarray}
 \|Q^{(n)}_{\gamma\delta} J_n\Big(J_n\psi_{\alpha,\beta}Q^{(n)}_{\alpha\beta}\Big)-Q^{(n)}_{\gamma\delta}J_n\psi_{\alpha,\beta}Q^{(n)}_{\alpha\beta}\|_{L^2(0,T;L^2)}\le \|Q^{(n)}\|_{L^\infty(0,T;L^4)}\|J_n-Id\|_{L^4\to L^4}\|J_n(\nabla \psi)Q^{(n)}\|_{L^2(0,T;L^4)}\nonumber\\
 \le \|Q^{(n)}\|_{L^\infty(0,T;L^4)}\|J_n-Id\|_{L^4\to L^4}\|J_n(\nabla \psi)\|_{L^\infty(0,T;L^{12})}\|Q^{(n)}\|_{L^2(0,T;L^6)}\to 0
 \end{eqnarray} where we denoted $\|J_n-Id\|_{L^4\to L^4}$ the norm of the operator $J_n-Id$ acting on $L^4$ and used the fact that this norm converges to zero, together with the bounds (\ref{weaksolapriori}). Thus we have the claim (\ref{lim:simplu}).

       Using a decomposition as in (\ref{lim:decomptwo}) with $Q^{(n)}_{\gamma\delta}Q^{(n)}_{\alpha\beta}$ instead of $Q^{(n)}_{\alpha\gamma}$ we get that $Q^{(n)}_{\gamma\delta}J_n\psi_{\alpha,\beta}Q^{(n)}_{\alpha\beta}$ converges strongly, in $L^2(0,T;L^2)$ to $Q_{\gamma\delta}\psi_{\alpha,\beta}Q_{\alpha\beta}$. This, together with (\ref{lim:simplu}) ensures that 
 
 \begin{equation}
\| Q^{(n)}_{\gamma\delta} J_n\Big(J_n\psi_{\alpha,\beta}Q^{(n)}_{\alpha\beta}\Big)-Q_{\gamma\delta}\psi_{\alpha,\beta}Q_{\alpha\beta}\|_{L^2(0,T;L^2)}\to 0
 \end{equation}
 
  Relations (\ref{convergences}) give that   $\Delta Q^{(n)}_{\gamma\delta}$, converges weakly in $L^2(0,T;L^2)$. Thus  we get convergence to the limit term 
 
 \begin{eqnarray}
  \int_0^\infty\int_{\mathbb{R}^d}Q_{\gamma\delta}\Delta Q_{\gamma\delta}\psi_{\alpha,\beta}Q_{\alpha\beta}\,dx\,dt= \int_0^T\int_{\mathbb{R}^d}Q_{\gamma\delta}\Delta Q_{\gamma\delta}\psi_{\alpha,\beta}Q_{\alpha\beta}\,dx\,dt
 \end{eqnarray}

  $\Box$

\section{ Higher regularity in $2$D, using the dissipation principle}

      In this section we restrict ourselves to dimension two and show that starting from an initial data with some higher regularity, we can obtain more regular solutions. More precisely, we have:

\begin{theorem}
Let $s>0$ and $(\bar Q,\bar u)\in H^{s+1}(\mathbb{R}^2)\times H^s(\mathbb{R}^2)$. There exists a global a solution $(Q(t,x),u(t,x))$ of the system (\ref{system}), with restriction (\ref{c+}), subject  to initial conditions
$$Q(0,x)=\bar Q(x),\, u(0,x)=\bar u(x)$$ and $Q\in L^2_{loc}(\mathbb{R}_+; H^{s+2}(\mathbb{R}^2))\cap L^\infty_{loc}(\mathbb{R}_+;H^{s+1}(\mathbb{R}^2))$, $u\in L^2_{loc}(\mathbb{R}_+;H^{s+1}(\mathbb{R}^2)\cap L^\infty_{loc}(\mathbb{R}_+;H^s)$.
      Moreover, we have:
{    
\begin{equation}
L\|\nabla Q(t,\cdot)\|_{H^s(\mathbb{R}^2)}^2+\|u(t,\cdot)\|_{H^s(\mathbb{R}^2)}^2\le C \Big(e+\|\bar Q\|_{H^{s+1}(\mathbb{R}^2)}+\|\bar u\|_{H^s(\mathbb{R}^2)}\Big)^{e^{e^{e^{Ct}}}}
\label{rate:3exp}
\end{equation} } where the constant $C$ depends only on $\bar Q, \bar u$, $a,b,c$, $\Gamma$ and $L$. 
     If $\xi=0$ the increase in time of the norms above can be made to be only doubly exponential.
\label{theorem:reg}
\end{theorem}

 \bigskip  
      The proof of the theorem is mainly based on $H^s$ energy estimates and the following cancelation {    (that is also used implicitly in showing the dissipation of the energy in Proposition ~\ref{prop:Lyapunov}): }

\begin{lemma}
\label{anulare}
For any symmetric matrices $Q', Q\in\mathbb{R}^{d\times d}$ and $\Omega_{\alpha\beta}=\frac{1}{2}(u_{\alpha,\beta}-u_{\beta,\alpha})\in\mathbb{R}^{d\times d}$ (decaying fast enough at infinity so that we can integrate by parts, in the formula below, without boundary terms) we have:
 $$\int_{\mathbb{R}^d} \textrm{tr}\big((\Omega Q' -Q'\Omega)\Delta Q\big)\,dx-\int_{\mathbb{R}^d}   \partial_\beta(Q'_{\alpha\gamma}\Delta
Q_{\gamma\beta}-\Delta Q_{\alpha\gamma} Q'_{\gamma\beta})u_\alpha\,dx=0$$
\end{lemma}

\smallskip      {\bf Proof.} We note that

\begin{eqnarray}\int_{\mathbb{R}^d} \textrm{tr}\big((\Omega Q' -Q'\Omega)\Delta Q\big)\,dx=\int_{\mathbb{R}^d}\Omega_{\alpha\gamma}Q'_{\gamma\beta}\Delta Q_{\beta\alpha}-Q'_{\alpha\gamma}\Omega_{\gamma\beta}\Delta Q_{\beta\alpha}
=\int_{\mathbb{R}^d}\Omega_{\alpha\gamma}Q'_{\gamma\beta}\Delta Q_{\beta\alpha}+\Omega_{\beta\gamma}
Q'_{\gamma\alpha}\Delta Q_{\alpha\beta}\nonumber\\=2\int_{\mathbb{R}^d}\textrm{tr}\big(\Omega Q'\Delta Q\big)\,dx=\underbrace{\int_{\mathbb{R}^d} u_{\alpha,\beta}Q'_{\beta\gamma}\Delta Q_{\gamma\alpha}\,dx}_{\mathcal{I}_1}-\underbrace{\int_{\mathbb{R}^d}
u_{\beta,\alpha}Q'_{\beta\gamma}\Delta Q_{\gamma\alpha}\,dx}_{\mathcal{I}_2}
\end{eqnarray} and on the other hand 
$$-\int_{\mathbb{R}^d}    \partial_\beta (Q'_{\alpha\gamma}\Delta Q_{\gamma\beta}) u_\alpha={    \int_{\mathbb{R}^d} Q'_{\alpha\gamma}\Delta Q_{\gamma \beta}   \partial_\beta u_\alpha=}\int_{\mathbb{R}^d} Q'_{\beta\gamma}\Delta Q_{\gamma \alpha}   \partial_\alpha u_\beta=I_2$$ 
and also 
$$\int_{\mathbb{R}^d}    \partial_\beta(\Delta Q_{\alpha\gamma} Q'_{\gamma\beta}) u_\alpha=-\int_{\mathbb{R}^d}  Q'_{\beta \gamma}\Delta Q_{\gamma \alpha}   \partial_\beta u_\alpha=-I_1$$
which finishes the proof. $\Box$

\smallskip     \begin{remark} The main point in the proof of the theorem is to use the previous lemma to eliminate the highest derivatives in $u$ {     in the first equation of the system (\ref{system}) and the highest derivatives in $Q$ in the second equation of the system}.The proof  could have been done, alternatively, by differentiating the equations $k\ge 1$ times and using the previous lemma. However that would have required estimating some delicate commutators and  would have restricted the initial data to $(\bar Q,\bar u)\in H^2\times H^1$. The Littlewood-Paley approach that we use allows for $(\bar Q,\bar u)\in H^{s+1}\times H^s$ with $s>0$.
\end{remark}

\smallskip
      In order to prove the theorem we need to introduce some technical preliminaries:

\subsection{Littlewood-Paley theory}
We  define
$\mathcal{C}$ to be the ring of center
$0$, of small radius $1/2$ and great radius $2$. There exist two
nonnegative  radial
functions $\chi$ and $\varphi$ belonging respectively to~${\mathcal{D}} 
(B(0,1)) $ and to
${\mathcal{D}} (\mathcal{C}) $ so that
\begin{equation}
\label{lpfond1}
\chi(\xi) + \sum_{q\geq 0} \varphi (2^{-q}\xi) = 1,\forall \xi\in\mathbb{R}^d
\end{equation}
\begin{equation}
\label{lpfond2}
|p-q|\geq 2
\Rightarrow
{\rm Supp}\,\, \varphi(2^{-q}\cdot)\cap {\rm Supp}\,\, \varphi(2^{-p}\cdot)=\emptyset.
\end{equation}

\noindent
For instance, one can take $\chi \in \mathcal{D} (B(0,1))$ such that $
\chi  \equiv 1 $ on $B(0,1/2)$ and take
$$
\varphi(\xi) = \chi(\xi/2) -
\chi(\xi).
$$
Then, we are able to define the Littlewood-Paley decomposition. Let us denote
by~$\mathcal{F}$ the Fourier transform on~$\mathbb{R}^d$. Let
$h,\
\tilde h,\  \Delta_q, S_q$ ($q \in \mathbb{Z}$) be defined as follows:

$$\displaylines{
\label{defnotationdyadique}h = {\mathcal F}^{-1}\varphi\quad {\rm and}\quad \tilde h =
{\mathcal{F}}^{-1}\chi, \cr
\Delta_q u = \mathcal{F}^{-1}(\varphi(2^{-q}\xi)\mathcal{F} u) = 2^{qd}\int h(2^qy)u(x-y)dy,\cr
S_qu
=\mathcal{F}^{-1}(\chi(2^{-q}\xi)\mathcal{F} u) =2^{qd} \int \tilde h(2^qy)u(x-y)dy.\cr
}
$$

      We recall that for two appropriately smooth functions $a$ and $b$ we have Bony's paraproduct decomposition
\cite{Bony81}:

\begin{equation}
ab=T_a b+T_b a+R(a,b)
\end{equation} where 
\begin{displaymath}
T_a b=\sum_{\substack{q'}}S_{q'-1} a\Delta_{q'}b,\,\, T_b a=\sum_{q'}S_{q'-1}b\Delta_{q'}a\textrm{ and }R(a,b)=\sum_{\substack{q',\\ i\in\{0,\pm 1\}} }\Delta_{q'} a\Delta_{q'+i} b.
\end{displaymath}
 Then we have

\begin{equation}
\Delta_q(ab)=\Delta_q T_a b+\Delta_q T_b a+\Delta_q R(a,b)=\Delta_q T_a b+\Delta_q \tilde R(a,b)
\end{equation} where $\tilde R(a,b)=T_b a+R(a,b)=\Sigma_{q'} S_{q'+2}b\Delta_{q'}a$. Moreover:

\begin{eqnarray}
\Delta_q (ab)=\Sigma_{|q'-q|\le 5}\Delta_q (S_{q'-1}a\Delta_{q'}b)+\Sigma_{q'> q-5}\Delta_q(S_{q'+2}b\Delta_{q'}a)\nonumber\\=\Sigma_{|q'-q|\le 5}[\Delta_q,S_{q'-1}a]\Delta_{q'}b+\Sigma_{|q'-q|\le 5}S_{q'-1}a\Delta_q\Delta_{q'}b+\Sigma_{q'> q-5}\Delta_q(S_{q'+2}b\Delta_{q'}a)\nonumber\\
=\Sigma_{|q'-q|\le 5}[\Delta_{q},S_{q'-1} a]\Delta_{q'}b+\Sigma_{|q'-q|\le 5}(S_{q'-1} a-S_{q-1}a)\Delta_q\Delta_{q'}b\nonumber\\
+\Sigma_{q'> q-5}\Delta_q (S_{q'+2} b\Delta_{q'}a)+\underbrace{\Sigma_{|q'-q|\le 5}S_{q-1}a\Delta_q\Delta_{q'}b}_{=S_{q-1}a\Delta_q b}
\label{bonydecomp}
\end{eqnarray}

      In terms of this decomposition we can express the Sobolev norm  of an element $u$ in the space $H^s$ as:

\begin{displaymath}
\|u\|_{H^s}=\big(\|S_0 u\|_{L^2}^2+\sum_{q\in\mathbb{N}}2^{2qs}\|\Delta_q u\|_{L^2}^2\big)^{1/2}
\end{displaymath}
{    
      We will use the following well-known estimates:

\begin{lemma} (\cite{chemin},\cite{chemin&masmoudi})
      {\bf (i)} (Bernstein inequalities) $$2^{-q}\|\nabla S_q u\|_{L^p}\le C \|u\|_{L^p}, \forall 1\le p\le \infty$$
$$\| \Delta_q u\|_{L^p}\le C2^{-q}\|\Delta_q \nabla u\|_{L^p}\le C\|\Delta_q u\|_{L^p}, \forall 1\le p\le \infty$$

    {\bf (ii)} (Bernstein inequalities) $$\|\Delta_q u\|_{L^b}\le 2^{d(\frac{1}{a}-\frac{1}{b})q}\|\Delta_q u\|_{L^a},\textrm{ for } b\ge a\ge 1$$
$$\|S_q u\|_{L^b}\le 2^{d(\frac{1}{a}-\frac{1}{b})q}\|S_q u\|_{L^a},\textrm{ for } b\ge a\ge 1$$

      {\bf (ii)}(commutator estimate) 

\begin{equation}
\|[\Delta_q, u]v\|_{L^p}\leq C 2^{-q}\|\nabla u\|_{L^r}\|v\|_{L^s}
\label{commutator}
\end{equation}  with $\frac{1}{p}=\frac{1}{r}+\frac{1}{s}$. The constant $C$ depends only on the function $\varphi$ used in defining $\Delta_q$ but not on $p,r,s$.
\label{lemma:bernstein&commutator}
\end{lemma}
}

Proof:  For the commutator estimate we begin by writing 
\begin{align*}
[\Delta_q , u]v(x)=\Delta_q(uv)(x)-u(x)\Delta_q v(x)&=2^{qd}\int h(2^q y)(u(x-y)-u(x))v(x-y) dy\\
&=-2^{qd}\int_{\mathbb{R}^d}\int_0^1 h(2^q y) y\nabla  u(x-\tau y) v(x-y) dy d\tau\\
&=-2^{-q} \int_{\mathbb{R}^d\times [0,1]}\tilde h_{2^q}(y)\nabla u(x-\tau y) v(x-y) dy d\tau,
\end{align*}
where $\tilde h(y)\stackrel{\rm{def}}{=}y h(y)\in \mathcal S(R^d)^d$ and $\tilde h_{\lambda}(y)\stackrel{\rm{def}}{=}\lambda^{d}\tilde h(\lambda y)$. Using the Cauchy-Schwartz inequality and a change of variables, we get
\begin{align*}
|[\Delta_q , u]v(x)|&\leq 2^{-q}\int_0^1\bigg(\int_{\mathbb{R}^d} |\tilde h_{2^q}(y)||\nabla u(x-\tau y)|^{\frac{r}{p}}dy\bigg)^{\frac pr}\bigg(\int_{\mathbb{R}^d}|\tilde h_{2^q}(y)|v(x-y)|^{\frac sp}dy\bigg)^{\frac ps}d\tau \\
&=2^{-q}\int_0^1 \bigg(\int_{\mathbb{R}^d} |\tilde h_{2^q\tau^{-1}}(y)||\nabla u(x- y)|^{\frac{r}{p}}dy\bigg)^{\frac pr}\bigg(\int_{\mathbb{R}^d}|\tilde h_{2^q}(y)|v(x-y)|^{\frac sp}dy\bigg)^{\frac ps}d\tau.
\end{align*}
Taking the $L^p$ norm in the $x$ variable, using the Cauchy-Schwartz inequality in the $x$ variable and  convolution estimates we obtain
\begin{align*}
\|[\Delta_q , u] v\|_{L^p}&\leq 2^{-q}\bigg(\int_0^1 \| |\tilde h_{2^q\tau^{-1}}|\star |\nabla u|^{\frac rp}\|_{L^p}^{\frac pr}d\tau\bigg) \||\tilde h_{2^q}|\star |v|^{\frac sp}\|_{L^p}^{\frac ps}\\
&\leq  2^{-q}\|\tilde h\|_{L^1} \|\nabla u\|_{L^r}\|v\|_{L^s},
\end{align*} so the constant in the inequality is $C=\|\tilde h\|_{L^1}$ and it does not depend on $p,r,s$.
\subsection{Proof of theorem ~\ref{theorem:reg}}

\smallskip      {\it Step 1. Estimates of the high frequencies}

\smallskip
\par We apply $\Delta_q$ to the first equation in (\ref{system}) and use the decomposition (\ref{bonydecomp}) to expand  $\Delta_q(\Omega_{\alpha\gamma} Q_{\gamma\beta})$, $\Delta_q (D_{\alpha\gamma}Q_{\gamma\beta})$ (and $\Delta_q(Q_{\alpha\gamma}\Omega_{\gamma\beta})$, $\Delta_q(Q_{\alpha\gamma}D_{\gamma\beta}$ ) as  $\Delta_q \Omega_{\alpha\gamma} S_{q-1} Q_{\gamma\beta}$, $S_{q-1}Q_{\gamma\beta}\Delta_q  D_{\alpha\gamma}$  (respectively $S_{q-1}Q_{\alpha\gamma}\Delta_q\Omega_{\gamma\beta}$, $S_{q-1}Q_{\alpha\gamma}\Delta_q D_{\gamma\beta}$) plus corrections. We also  expand $\Delta_q(Q_{\alpha\beta}\textrm{tr}(Q\nabla u))$ as $S_{q-1}Q_{\alpha\beta}S_{q-1}Q_{\gamma\delta}\Delta_q u_{\gamma,\delta}$  plus corrections (by applying the formula (\ref{bonydecomp}) twice) and we get:
\begin{eqnarray}
   \partial_t \Delta_q Q_{\alpha\beta}- \Gamma L\Delta\Delta_q Q_{\alpha\beta}-\Delta_q \Omega_{\alpha\gamma} S_{q-1} Q_{\gamma\beta}+S_{q-1}Q_{\alpha\gamma}\Delta_q\Omega_{\gamma\beta}-\xi S_{q-1}Q_{\gamma\beta}\Delta_q  D_{\alpha\gamma}-\xi S_{q-1}Q_{\alpha\gamma}\Delta_q D_{\gamma\beta}\nonumber\\-\xi\Delta_q D_{\alpha\beta}+2\xi S_{q-1}Q_{\alpha\beta}S_{q-1}Q_{\gamma\delta}\Delta_q u_{\gamma,\delta}+\xi\delta_{\alpha\beta} \Delta_q(\textrm{tr}(Q\nabla u))=\big(\mathcal{T}_Q\big)_{\alpha\beta}\nonumber
   \end{eqnarray} where $\mathcal{T}_Q$ denotes the sum of the correction terms mentioned before together with some other terms that are easy to estimate using the apriori bounds in Proposition ~\ref{prop:aprioriest}. These terms are described in the Appendix $A$.

\par  Multiplying the previous equation by $-L\Delta \Delta_q Q_{\alpha\beta}$ and integrating over $\mathbb{R}^2$ and by parts we obtain:

\begin{eqnarray}
\frac{L}{2}   \partial_t \|\nabla \Delta_q Q\|_{L^2}^2+\Gamma L^2\|\Delta\Delta_q Q\|_{L^2}^2+L\int \Delta_q\Omega_{\alpha\gamma} S_{q-1} Q_{\gamma\beta}\Delta\Delta_q Q_{\alpha\beta}-L\int S_{q-1} Q_{\alpha\gamma}\Delta_q \Omega_{\gamma\beta}\Delta\Delta_q Q_{\alpha\beta}\nonumber\\+L\xi \int S_{q-1}Q_{\gamma\beta}\Delta_q  D_{\alpha\gamma}\Delta\Delta_q Q_{\alpha\beta}+L\xi\int S_{q-1}Q_{\alpha\gamma}\Delta_q D_{\gamma\beta}\Delta\Delta_q Q_{\alpha\beta}\nonumber\\+L\xi\int \Delta_q D_{\alpha\beta}\Delta\Delta_q Q_{\alpha\beta}-2L\xi\int S_{q-1}Q_{\alpha\beta}\Delta\Delta_q Q_{\alpha\beta}S_{q-1}Q_{\gamma\delta}\Delta_q u_{\gamma,\delta}=\Big(-L\Delta\Delta_q Q_{\alpha\beta},\big(\mathcal{T}_Q\big)_{\alpha\beta}\Big)
\label{firsthalf}
\end{eqnarray} where the terms on the right hand side are described in the Appendix $A$. 

 We apply  $\Delta_q$ to the second equation in (\ref{system}) and  use the decomposition (\ref{bonydecomp}) to expand  
$\Delta_q( Q_{\alpha\gamma}\Delta Q_{\gamma\beta}$ (respectively $\Delta_q (\Delta Q_{\alpha\gamma} Q_{\gamma\beta})$) as
$S_{q-1}Q_{\alpha\gamma}\Delta_q\Delta Q_{\gamma\beta}$ (respectively  $\Delta_q\Delta Q_{\alpha\gamma} S_{q-1}Q_{\gamma\beta}$)  plus correction terms. We also  expand $\Delta_q(Q_{\alpha\beta}\textrm{tr}(Q\Delta Q))$ as $S_{q-1}Q_{\alpha\beta}S_{q-1}Q_{\gamma\delta}\Delta_q \Delta Q_{\gamma\delta}$  plus corrections (by applying the formula (\ref{bonydecomp}) twice) and we get:

 we get:

\begin{eqnarray}
   \partial_t\Delta_q u_\alpha-\nu\Delta\Delta_q u_\alpha=   \partial_\alpha\Delta_q  p+L   \partial_\beta\left(S_{q-1}Q_{\alpha\gamma}\Delta_q\Delta Q_{\gamma\beta}-\Delta_q\Delta Q_{\alpha\gamma} S_{q-1}Q_{\gamma\beta}\right)\nonumber\\
-L\xi   \partial_\beta\big( S_{q-1}Q_{\alpha\gamma}\Delta_q\Delta Q_{\gamma\beta}+\Delta_q\Delta Q_{\alpha\gamma} S_{q-1} Q_{\gamma\beta}-2 S_{q-1} Q_{\alpha\beta}S_{q-1}Q_{\gamma\delta}\Delta \Delta_q Q_{\gamma\delta})\bigg)\nonumber\\
-\xi   \partial_\beta\big(L\Delta_q \Delta Q_{\alpha\beta}-\delta_{\alpha\beta}\Delta_q(\textrm{tr}(QH)\big)+\big(\mathcal{T}_u\big)_\alpha\nonumber
\end{eqnarray}  where $\mathcal{T}_u$ denotes the sum of the correction terms mentioned before together with some other terms that are easy to estimate using the apriori bounds in Proposition ~\ref{prop:aprioriest}. These term are described in the Appendix $A$.

\par      We multiply the last equation by $\Delta_q u_\alpha$, integrate over $\mathbb{R}^2$ and by parts to obtain:

\begin{eqnarray}
\frac{1}{2}   \partial_t\|\Delta_q u\|_{L^2}^2+\nu\|\Delta_q\nabla u\|_{L^2}^2+L\int S_{q-1}Q_{\alpha\gamma}\Delta_q\Delta Q_{\gamma\beta}\Delta_q u_{\alpha,\beta}-L\int \Delta_q\Delta Q_{\alpha\gamma} S_{q-1}Q_{\gamma\beta}\Delta_q u_{\alpha,\beta}\nonumber\\-L\xi\bigg(\int S_{q-1}Q_{\alpha\gamma}\Delta_q\Delta Q_{\gamma\beta}\Delta_q u_{\alpha,\beta}+\int\Delta_q\Delta Q_{\alpha\gamma} S_{q-1} Q_{\gamma\beta}\Delta_q u_{\alpha,\beta}-2\int S_{q-1} Q_{\alpha\beta}S_{q-1}Q_{\gamma\delta}\Delta_q\Delta Q_{\gamma\delta}\Delta_q u_{\alpha,\beta}\bigg)\nonumber\\
-\xi L\int \Delta_q\Delta Q_{\alpha\beta}\Delta_q u_{\alpha,\beta}=\Big(\big(\mathcal{T}_u\big)_\alpha,\Delta_q u_\alpha\Big)
\label{secondhalf}
\end{eqnarray}

      Summing (\ref{firsthalf}) and (\ref{secondhalf}) and using Lemma ~\ref{anulare} we get:

\begin{eqnarray}
   \partial_t \left(\frac{L}{2}\|\nabla \Delta_q Q\|_{L^2}^2+\frac{1}{2}\|\Delta_q u\|_{L^2}^2\right)+\nu\|\Delta_q\nabla u\|_{L^2}^2+\Gamma L^2\|\Delta\Delta_q Q\|_{L^2}^2=\Big(-L\Delta\Delta_q Q_{\alpha\beta},\big(\mathcal{T}_Q\big)_{\alpha\beta}\Big)+\Big(\big(\mathcal{T}_u\big)_\alpha,\Delta_q u_\alpha\Big)\nonumber
\end{eqnarray} 

      We denote by $\varphi(t)\stackrel{def}{=}L\|\nabla Q\|_{H^s}^2+\|u\|_{H^s}^2$ with $\varphi_1(t)\stackrel{def}{=}L\|S_0 \nabla Q\|_{L^2}^2+\|S_0 u\|_{L^2}^2$ the low-frequency part of $\varphi$ and  $\varphi_2(t)\stackrel{def}{=}\varphi(t)-\varphi_1(t)$ the high-frequency part of $\varphi$.

      The last inequality leads to the following estimate, that holds for any $\varepsilon\in (0,\frac{1}{2})$ and  whose technical proof is postponed to  Appendix $B$:

\begin{eqnarray}
\frac{1}{2}\frac{d}{dt}\varphi_2+\sum_{q\in\mathbb{N}} 2^{2qs}\Big(\frac{\Gamma L^2}{2}\|\Delta\Delta_q Q\|_{L^2}^2+\frac{\nu}{2}\|\nabla\Delta_q u\|_{L^2}^2\Big)\nonumber\\ 
\leq C\Big(1+\|\nabla u\|_{L^2}^2+\|u\|_{L^2}^2\|\nabla u\|_{L^2}^2+\|\nabla Q\|_{L^2}^2\|\Delta Q\|_{L^2}^2\Big)(\|\nabla Q\|_{H^s}^2+\|u\|_{H^s}^2)\nonumber\\+C\big(\|Q\|_{L^2}+\|Q\|_{L^4}^2\big)^2\|\nabla Q\|_{H^s}^2\nonumber\\+\frac{\Gamma L^2}{50}\|\Delta Q\|_{H^s}^2+\frac{\nu}{50}\|\nabla u\|_{H^s}^2+\xi^2 C\Big(1+\|\nabla u\|_{L^2}^2+\|u\|_{L^2}^2\|\nabla u\|_{L^2}^2+\|\nabla Q\|_{L^2}^2\|\Delta Q\|_{L^2}^2\Big)(\|\nabla Q\|_{H^s}^2+\|u\|_{H^s}^2)\nonumber\\
+\xi^2 C\bigg((1+\|Q\|_{L^\infty}^2\|\nabla u\|_{L^2}^2)\|\nabla Q\|_{H^s}^2+\sum_{j=2}^5\|Q\|_{L^{2(j-1)}}^{2(j-1)}\|\nabla Q\|_{H^s}^2\bigg)\nonumber\\
+\xi^2 C\big(\|\nabla Q\|_{L^{\frac{2}{\varepsilon}}}\|Q\|_{L^\infty}\big)^{\frac{2}{1-\varepsilon}}\|u\|_{H^s}^2
\label{longest}
\end{eqnarray}  where the constant $C$ is independent of $\xi$.

\bigskip     {\it Step 2. Estimates of the low frequencies}
 \smallskip\par     This is much easier than the previous step. We apply the operator $S_0$ to the first equation in (\ref{system}), multiply by $-LS_0{    \Delta}Q_{\alpha\beta}$, take the trace, integrate over $\mathbb{R}^2$ and by parts  and we get:
\begin{eqnarray}
\frac{L}{2}   \partial_t \|S_0\nabla  Q\|_{L^2}^2+\Gamma L^2\|\Delta S_0 Q\|_{L^2}^2\le \|u\|_{L^4}\|\nabla Q\|_{L^4}\|\Delta S_0 Q\|_{L^2}+C\|S_0(Q\nabla u)\|_{L^2}\|\Delta S_0 Q\|_{L^2}\nonumber\\
+L\|S_0\big(-aQ+b[Q^2-\frac{\textrm{tr}(Q^2)}{3}Id]-cQ\textrm{tr}(Q^2)\big)\|_{L^2}\|S_0 \Delta Q\|_{L^2}\nonumber\\
+C\xi\|S_0(\nabla uQ)\|_{L^2}\|\Delta S_0 Q\|_{L^2}+C\xi\|S_0(\nabla u)\|_{L^2}\|\Delta S_0 Q\|_{L^2}+C\xi \|S_0(Q^2\nabla u)\|_{L^2}\|S_0\Delta Q\|_{L^2}\nonumber\\
\le C\|u\|_{L^2}^{\frac{1}{2}}\|\nabla u\|_{L^2}^{\frac{1}{2}}\|\nabla Q\|_{L^2}^{\frac{1}{2}}\|\Delta Q\|_{L^2}^{\frac{1}{2}}\|\Delta S_0 Q\|_{L^2}+C\|Q\nabla u\|_{L^1}\|\Delta S_0 Q\|_{L^2}
\nonumber\\+C\|-aQ+b[Q^2-\frac{\textrm{tr}(Q^2)}{3}Id]-cQ\textrm{tr}(Q^2)\|_{L^2}^2+\frac{\Gamma L^2}{100}\|\Delta S_0 Q\|_{L^2}\nonumber\\
+C\xi^2\big(\|\nabla u Q\|_{L^1}^2+\|\nabla u\|_{L^2}^2+\|Q^2\nabla u\|_{L^1}^2)\nonumber
\end{eqnarray} hence

\begin{eqnarray}
\frac{L}{2}   \partial_t \|S_0\nabla Q\|_{L^2}^2+\frac{\Gamma L^2}{2}\|\Delta S_0 Q\|_{L^2}^2\le C\|u\|_{L^2}^2\|\nabla u\|_{L^2}^2+C\|\nabla Q\|_{L^2}^2\|\Delta Q\|_{L^2}^2+\|Q\|_{L^2}^2\|\nabla u\|_{L^2}^2\nonumber\\
+C(\|Q\|_{L^2}^2+\|Q\|_{L^4}^4+\|Q\|_{L^6}^6)+C\xi^2(1+\|Q\|_{L^2}^2+\|Q\|_{L^4}^4)\|\nabla u\|_{L^2}^2
\label{eq:low1}
\end{eqnarray}

      We aply $S_0$ to the second equation in (\ref{system}), multiply by $S_0 u$ and integrate over $\mathbb{R}^2$ and by parts to obtain:
\begin{eqnarray}
\frac{1}{2}   \partial_t \|S_0 u\|_{L^2}^2+\nu\|\nabla S_0 u\|_{L^2}^2\le \|S_0(u\nabla u)\|_{L^2}\|S_0 u\|_{L^2}+C\|S_0(\nabla Q\Delta Q)\|_{L^2}\|S_0u\|_{L^2}+C\|S_0(Q\Delta Q)\|_{L^2}\|S_0\nabla u\|_{L^2}\nonumber\\
+C\xi\|S_0\big((Q+Q^2)(L\Delta Q-aQ+b[Q^2-\frac{\textrm{tr}(Q^2)}{3}Id]-c Q\textrm{tr}(Q^2))\big)\|_{L^2}\|S_0\nabla u\|_{L^2}\nonumber\\
+C\xi\|S_0\big(L\Delta Q-aQ+b[Q^2-\frac{\textrm{tr}(Q^2)}{3}Id]-c Q\textrm{tr}(Q^2)\big)\|_{L^2}\|S_0\nabla u\|_{L^2}\nonumber\\
\le C\|u\nabla u\|_{L^1}^2+C\|S_0 u\|_{L^2}^2+C\|\nabla Q\Delta Q\|_{L^1}^2+C\|S_0 u\|_{L^2}^2+C\|Q\Delta Q\|_{L^1}^2+\frac{\nu}{2}\|S_0\nabla u\|_{L^2}^2\nonumber\\
+C\xi^2\Big[ (1+\|Q\|_{L^2}+\|Q\|_{L^4}^4)\|\Delta Q\|_{L^2}^2+\sum_{j=2}^5\|Q\|_{L^j}^j\Big]\nonumber
\end{eqnarray} hence

\begin{eqnarray}
\frac{1}{2}   \partial_t\|S_0 u\|_{L^2}^2+\frac{\nu}{2}\|\nabla S_0 u\|_{L^2}^2\le C\|u\|_{L^2}\|\nabla u\|_{L^2}^2+C\|\nabla Q\|_{L^2}^2\|\Delta Q\|_{L^2}^2+ C\|Q\|_{L^2}^2\|\Delta Q\|_{L^2}^2+C\|S_0 u\|_{L^2}^2\nonumber\\
+C\xi^2\Big[ (1+\|Q\|_{L^2}+\|Q\|_{L^4}^4)\|\Delta Q\|_{L^2}^2+\sum_{j=2}^5\|Q\|_{L^j}^j\Big]
\label{eq:low2}
\end{eqnarray}

      Summing (\ref{eq:low1}) and (\ref{eq:low2}) we obtain:
\begin{equation}
   \partial_t \varphi_1+\frac{\nu}{2}\|\nabla S_0 u\|_{L^2}^2+\frac{\Gamma L^2}{2}\|\Delta S_0 Q\|_{L^2}^2\le C\varphi+m(t)+\xi^2 n(t)\label{est:low}
\end{equation} where 

\begin{equation}
m(t)\stackrel{def}{=}C\Big(\|u\|_{L^2}^2\|\nabla u\|_{L^2}^2+\|\nabla Q\|_{L^2}^2\|\Delta Q\|_{L^2}^2+\|Q\|_{L^2}^2\|\Delta Q\|_{L^2}^2+\|Q\|_{L^2}^2\|\nabla u\|_{L^2}^2+\|Q\|_{L^2}^2+\|Q\|_{L^4}^4+\|Q\|_{L^6}^6\Big)\nonumber
\end{equation} and 

\begin{equation}
n(t)\stackrel{def}{=}(1+\|Q\|_{L^2}^2+\|Q\|_{L^4}^4)\|\nabla u\|_{L^2}^2+ (1+\|Q\|_{L^2}+\|Q\|_{L^4}^4)\|\Delta Q\|_{L^2}^2+\sum_{j=2}^5\|Q\|_{L^j}^j\nonumber
\end{equation}

\medskip     {\it Step 3. The estimates of the high norms}

      Summing (\ref{longest}) and (\ref{est:low}) we obtain:

\begin{eqnarray}\frac{1}{2}\varphi'(t)\le C\underbrace{\Big(1+\|\nabla u\|_{L^2}^2+\|u\|_{L^2}^2\|\nabla u\|_{L^2}^2+\|\nabla Q\|_{L^2}^2\|\Delta Q\|_{L^2}^2+\|Q\|_{L^2}^2+\|Q\|_{L^4}^4\Big)}_{\stackrel{def}{=}u(t)}\varphi(t)+m(t)\nonumber\\
+C\xi^2 \underbrace{\Big(1+\|\nabla u\|_{L^2}^2+\|u\|_{L^2}^2\|\nabla u\|_{L^2}^2+\|\nabla Q\|_{L^2}^2\|\Delta Q\|_{L^2}^2+\sum_{j=1}^4\|Q\|_{L^{2j}}^{2j}\Big)}_{\stackrel{def}{=}v(t)}\varphi(t)\nonumber\\
+\xi^2 C\|Q\|_{L^\infty}^2\|\nabla u\|_{L^2}^2\varphi(t)+\xi^2 C\big(\|\nabla Q\|_{L^{\frac{2}{\varepsilon}}}\|Q\|_{L^\infty}\big)^{\frac{2}{1-\varepsilon}}\varphi(t)+\xi^2 n(t)\nonumber
\end{eqnarray} where $u(t),v(t),m(t)$ and $n(t)$ are, by Proposition ~\ref{prop:aprioriest}, apriori bounded in $L^2(0,T)$, and increasing exponentially in time.

   If $\xi=0$ the above estimates together with Gronwall's lemma show that $\varphi$ increases like $e^{e^{ct}}$ for an appropriate constant $c>0$.

      In the general case, when $\xi\not=0$  we start by recalling use a fundamental ingredient in the global existence, namely the logarithmic estimate (see \cite{brezisgallouet}), for $s>0$,
$$\|Q\|_{L^\infty}\leq \|Q\|_{H^1}\sqrt{\ln(e+\frac{\|\nabla Q\|_{H^s}^2}{\|Q\|_{H^1}})},$$
and be denoting $f(t)\stackrel{def}{=}\|Q\|_{H^1}^2$ and we obtain
\begin{eqnarray}\varphi'(t)\leq C\Big(u(t)+\xi v(t)\Big)\varphi(t)+m(t)\nonumber\\
+\xi C f(t)\|\nabla u\|_{L^2}^2\ln(e+\frac{\varphi(t)}{\sqrt{f(t)}})\varphi(t)+\xi C\big(\|\nabla Q\|_{L^{\frac{2}{\varepsilon}}}\|Q\|_{L^\infty}\big)^{\frac{2}{1-\varepsilon}}\varphi(t)+\xi^2 n(t)\nonumber
\end{eqnarray}

      Observing that the function $h(x)\stackrel{def}{=}x\ln(e+\frac{\varphi}{\sqrt{x}})$  is increasing  the last relation implies:
\begin{eqnarray}
\varphi'(t)\le C\Big(u(t)+\xi v(t)\Big)\varphi(t)+m(t)\nonumber\\
+\xi C(1+f(t))\|\nabla u\|_{L^2}^2\varphi(t)\Big(\ln(e+\varphi(t))\Big)+\xi C\big(\|\nabla Q\|_{L^{\frac{2}{\varepsilon}}}\|Q\|_{L^\infty}\big)^{\frac{2}{1-\varepsilon}}\varphi(t)+\xi^2 n(t)
\label{eq:varphieq1}
\end{eqnarray}

      On the other hand, by using the interpolation inequality (see \cite{cheminxu}, and also \cite{mp}, Lemma $10$):
\begin{equation}
\|g\|_{L^{2p}}\le C\sqrt{p}\|g\|_{L^2}^{\frac{1}{p}}\|\nabla g \|_{L^2}^{1-\frac{1}{p}}
\label{trickinterpolation} 
\end{equation}
 we get:

$$\|\nabla Q\|_{L^{\frac{2}{\varepsilon}}}^{\frac{2}{1-\varepsilon}}\le\Big(\frac{1}{\varepsilon}\Big)^{\frac{1}{1-\varepsilon}}\|\nabla Q\|^{\frac{2\varepsilon}{1-\varepsilon}}\|\Delta Q\|_{L^2}^2\le\Big(\frac{1}{\varepsilon}\Big)^{\frac{1}{1-\varepsilon}} (1+\|\nabla Q\|_{L^2}^2)\|\Delta Q\|_{L^2}^2$$ where for the last inequality we assumed $0<\varepsilon<\frac{1}{2}$.

      Then (\ref{eq:varphieq1}) becomes:
\begin{eqnarray}
\varphi'(t)\le C\Big(u(t)+\xi v(t)\Big)\varphi(t)+m(t)\nonumber\\
+\xi C(1+f(t))\|\nabla u\|_{L^2}^2\varphi(t)\Big(\ln(e+\varphi(t))\Big)+\xi^2 n(t)\nonumber\\
\xi C(1+f(t))\|\Delta Q\|_{L^2}^2\big[(1+f(t))\ln(e+\varphi(t))\big]^{\frac{1}{1-\varepsilon}}\left(\frac{1}{\varepsilon}\right)^{\frac{1}{1-\varepsilon}}\varphi(t)
\label{varphieq2} 
\end{eqnarray}
      Observing that the constants in the interpolation inequality (\ref{trickinterpolation}) and in the commutator estimate (\ref{commutator}) do not depend on the space $L^p$ that we work with  and denoting $N\stackrel{def}{=}\ln(e+\varphi)$ we choose

$$\varepsilon\stackrel{def}{=}(1+\ln N)^{-1}$$ and observing that $[N(1+\ln N)]^{1+\frac{1}{\ln N}}\le CN(1+\ln N)$ for some  constant $C$ independent of $N$,  the last inequality becomes:

\begin{eqnarray}
\varphi'(t)\le C\Big(u(t)+\xi v(t)\Big)\varphi(t)+m(t)\nonumber\\
+\xi C(1+f(t))\|\nabla u\|_{L^2}^2\varphi(t)\Big(\ln(e+\varphi(t))\Big)+\xi^2 n(t)\nonumber\\
\xi C(1+f(t))^3\|\Delta Q\|_{L^2}^2\varphi(t)\ln\Big (e+\varphi(t)\Big)\bigg(1+\ln(e+\ln(\varphi(t)+e)\big)\bigg)
\label{varphieq1+} 
\end{eqnarray}

$\Box$

\section{Weak-Strong uniqueness in 2D}
In this section we consider a global weak solution and a strong one, starting from the same  initial data $(\bar Q,\bar u)\in H^{s+1}(\mathbb{R}^2)\times H^s(\mathbb{R}^2)$ with $s>0$ and we show that they are the same. More precisely:

\begin{proposition} Let $(\bar Q,\bar u)\in H^{s+1}(\mathbb{R}^2)\times H^s(\mathbb{R}^2)$ with $s>0$. By Proposition ~\ref{prop:weak} there exists a weak solution $(Q_1, u_1)$ of the system (\ref{system}), subject to restriction (\ref{c+}) and starting from initial data $(\bar Q,\bar u)$, such that 

 \begin{equation}Q_1\in L^\infty_{loc}(\mathbb{R}_+;H^1(\mathbb{R}^2))\cap L^2_{loc}(\mathbb{R}_+;  H^2(\mathbb{R}^2))\textrm{ and } u_1\in L^\infty_{loc}(\mathbb{R}_+;L^2(\mathbb{R}^2))\cap L^2_{loc}(\mathbb{R}_+; H^1(\mathbb{R}^2))
 \label{sol1}
 \end{equation}

       Theorem ~\ref{theorem:reg} gives the existence of a strong solution $(Q_2,u_2)$ such that

 \begin{equation}Q_2\in L^\infty_{loc}(\mathbb{R}_+; H^{s+1}(\mathbb{R}^2)\cap L^2_{loc}(\mathbb{R}_+; H^{s+2}(\mathbb{R}^2))\textrm{ and }u_2\in L^\infty(\mathbb{R}_+; H^s(\mathbb{R}^2))\cap L^2(\mathbb{R}_+; H^{s+1}(\mathbb{R}^2))
 \label{sol2}
 \end{equation} with $s>0$ and the same initial data $(\bar Q,\bar u)\in H^{s+1}(\mathbb{R}^2)\times H^s(\mathbb{R}^2)$. Then  $(Q_1,u_1)=(Q_2, u_2)$. 
\end{proposition}

\smallskip     {\bf Proof.}
We denote by $\delta Q=Q_1-Q_2$ and $\delta u=u_1-u_2$ which verify the following system
\begin{equation}
\left\{\begin{array}{l}
      (   \partial_t+  \delta u\nabla )\delta Q-\delta \Omega \delta Q+\delta Q\delta \Omega+\delta u\nabla Q_2+u_2\nabla\delta Q+Q_2\delta\Omega+\delta Q\Omega_2-\delta\Omega Q_2-\Omega_2\delta Q\\
 -\xi\big[\delta D\delta Q  +\delta Q\delta D+\delta D-2(\delta Q+\frac{1}{2}Id)\textrm{tr}(\delta Q\nabla\delta u)\big]\\
 -\xi[\delta D Q_2+D_2\delta Q+\delta Q D_2+ Q_2 \delta D-\textrm{tr}(\delta Q\nabla u_2)Id-\textrm{tr}(Q_2\nabla\delta u)Id]\\
 -2\xi\big[ \delta Q\textrm{tr}(\delta Q\nabla u_2)+\delta Q\textrm{tr}(Q_2\nabla \delta u)+\delta Q\textrm{tr}(Q_2\nabla u_2)+Q_2\textrm{tr}(\delta Q\nabla\delta u)+Q_2\textrm{tr}(\delta Q\nabla u_2)+Q_2\textrm{tr}(Q_2\nabla\delta u)\big]\\
   =\Gamma\Big(L\Delta \delta Q-a\delta Q+b[\delta QQ_1+Q_2\delta Q-\frac{\textrm{tr}\big(\delta QQ_1+Q_2\delta Q\big)}{2}Id]-c\delta Q\textrm{tr}(Q_1^2)-cQ_2\big[\textrm{tr}(Q_1\delta Q+\delta QQ_2)\big]\Big) \\
         \partial_t \delta u +\mathcal{P}(\delta u\nabla \delta u)=\nu\Delta \delta u-L\mathcal{P}\big(\nabla\cdot(\nabla \delta Q\nabla \delta Q-\frac{1}{2}|\nabla \delta Q|^2)\big)+L\mathcal{P}\big(\nabla\cdot(\delta Q\Delta \delta Q-\Delta \delta Q\delta Q)\big)\\
-\xi\nabla\cdot \big[\delta Q\delta H+\delta H\delta Q+\delta H-2(\delta Q+\frac{1}{2}Id)\textrm{tr}(\delta Q\delta H)\big]\\
     -\mathcal{P}(u_2\nabla\delta u+\delta u\nabla u_2)-L \mathcal{P}\bigg(\nabla\cdot\Big((\nabla\delta Q\nabla Q_2+\nabla Q_2\nabla\delta Q)-\frac{1}{2}\textrm{tr}\big(\nabla\delta Q\nabla Q_2+Q_2\nabla\delta Q\big)Id\Big)\bigg)\\
  -\xi \nabla\cdot\big[\delta QH_2+Q_2\delta H+\delta HQ_2+H_2\delta Q-\textrm{tr}(\delta QH_2)Id-(Q_2\delta H)Id\big]\\
-2\xi\nabla\cdot\big[\delta Q\textrm{tr}(\delta QH_2)+\delta Q\textrm{tr}(Q_2\delta H)+Q_2\textrm{tr}(\delta Q\delta H)+ Q_2\textrm{tr}( Q_2\delta H )+Q_2\textrm{tr}(\delta Q H_2 )+\delta Q\textrm{tr}( Q_2 H_2 )\big]\nonumber\\
 +L\mathcal{P}\big(\nabla\cdot(\delta Q\Delta Q_2+Q_2\Delta \delta Q-\Delta\delta Q Q_2-\Delta Q_2\delta Q)\big)\\
     \end{array}\right.
      \label{approxsystemdif}
      \end{equation}
      
\smallskip  
      We proceed similarly as in the proof of Proposition ~\ref{prop:Lyapunov}, namely we multiply  the first equation in (\ref{approxsystemdif}) to the right by $-L\Delta\delta Q+\delta Q$, integrate over $\mathbb{R}^2$ and by parts, take the trace and sum with the second equation in (\ref{approxsystemdif}) multiplied by $\delta u$ and integrated over $\mathbb{R}^2$ and by parts.
 Taking into account the cancellations analogous to the ones in (\ref{Lyapunovcancellation}) we obtain:
 \begin{align}
 \frac{d}{dt}\int_{\mathbb{R}^2}\frac{L}{2}|\nabla\delta Q(x)|^2+\frac{1}{2}|\delta Q(x)|^2+\frac{1}{2}|\delta u(x)|^2\,dx+\int_{\mathbb{R}^2}\nu|\nabla\delta u(x)|^2+\Gamma L^2 |\Delta\delta Q(x)|^2\,dx&\nonumber\\
 =L\int_{\mathbb{R}^2}\textrm{tr}\Big(\big[\delta u\nabla Q_2+u_2\nabla\delta Q+\delta Q\Omega_2-\Omega_2\delta Q\big]\Delta\delta Q\Big)\,dx+\underbrace{L\int_{\mathbb{R}^2}\textrm{tr}\Big(\big[Q_2\delta\Omega-\delta\Omega Q_2\big]\Delta\delta Q\Big)\,dx}_{\mathcal{A}}&\nonumber\\
 -\underbrace{\xi L\int_{\mathbb{R}^2}[\delta D Q_2+ Q_2\delta D]\Delta\delta Q\,dx}_{\mathcal{B}}-\xi L\int_{\mathbb{R}^2}[D_2\delta Q+\delta Q D_2]\Delta \delta Q\,dx&\nonumber\\
 -\underbrace{2\xi L\int_{\mathbb{R}^2}\big[ \delta Q\textrm{tr}(Q_2\nabla \delta u)+Q_2\textrm{tr}(\delta Q\nabla\delta u)+Q_2\textrm{tr}(Q_2\nabla\delta u)\big]\Delta\delta Q\,dx}_{\mathcal{C}}&\nonumber\\
  -2\xi L\int_{\mathbb{R}^2}\big[ \delta Q\textrm{tr}(\delta Q\nabla u_2)+
  \delta Q\textrm{tr}(Q_2\nabla u_2)+Q_2\textrm{tr}(\delta Q\nabla u_2)\big]\Delta\delta Q\,dx&\nonumber\\
  +\xi\int_{\mathbb{R}^2}[\delta D Q_2+D_2\delta Q+\delta Q D_2+ Q_2\delta D] \delta Q\,dx&\nonumber\\
 +2\xi \int_{\mathbb{R}^2}\big[ \delta Q\textrm{tr}(\delta Q\nabla u_2)+\delta Q\textrm{tr}(Q_2\nabla \delta u)+Q_2\textrm{tr}(\delta Q\nabla\delta u)+
  \delta Q\textrm{tr}(Q_2\nabla u_2)+Q_2\textrm{tr}(\delta Q\nabla u_2)+Q_2\textrm{tr}(Q_2\nabla\delta u)\big]\delta Q\,dx&\nonumber
  \end{align}
  \begin{align}
 -a\Gamma L\int_{\mathbb{R}^2}|\nabla\delta Q(x)|^2\,dx -b\Gamma L\int_{\mathbb{R}^2}\textrm{tr}\Big(\big(\delta Q(x) Q_1(x) +Q_2(x)\delta Q(x)\big)\Delta\delta Q(x)\Big)\,dx&\nonumber\\
+c\Gamma L\int_{\mathbb{R}^2}\textrm{tr}\big(\delta Q\Delta\delta Q\big)\textrm{tr}(Q_1^2)\,dx+c\Gamma L\int_{\mathbb{R}^2}\textrm{tr}(Q_2\Delta\delta Q)\textrm{tr}(Q_1\delta Q+\delta Q Q_2)\,dx&\nonumber\\
-\int_{\mathbb{R}^2}\textrm{tr}\big(\delta u\nabla Q_2\delta Q\big)\,dx-\int_{\mathbb{R}^2}\textrm{tr}\big(Q_2\delta\Omega\delta Q\big)\,dx-\underbrace{\int_{\mathbb{R}^2}\textrm{tr}\big(\delta Q\Omega_2\delta Q\big)dx}_{\mathcal{I}}&\nonumber\\
+\int_{\mathbb{R}^2}\textrm{tr}\big(\delta\Omega Q_2\delta Q\big)\,dx+\underbrace{\int_{\mathbb{R}^2}\textrm{tr}\big(\Omega_2(\delta Q)^2\big)\,dx}_{\mathcal{II}}-\Gamma L\int_{\mathbb{R}^2}|\nabla Q|^2\,dx&\nonumber\\
-a\Gamma\int_{\mathbb{R}^2}|\delta Q|^2\,dx+b\Gamma \int_{\mathbb{R}^2}\textrm{tr}\big(\delta Q Q_1\delta Q+Q_2(\delta Q)^2\big)\,dx&\nonumber\\
-c\Gamma\int_{\mathbb{R}^2}\textrm{tr}(Q_1)^2|\delta Q|^2\,dx-c\Gamma\int_{\mathbb{R}^2}\textrm{tr}(Q_2\delta Q)\textrm{tr}(Q_1\delta Q+\delta QQ_2)\,dx&\nonumber\\
-{    \int_{\mathbb{R}^2}\big(u_2\nabla \delta u+\delta u\nabla u_2)\delta u\,dx}+L\int_{\mathbb{R}^2}\big(\nabla\delta Q\nabla Q_2+\nabla Q_2\nabla\delta Q\big)\cdot\nabla \delta u\,dx&\nonumber\\
+\xi\int_{\mathbb{R}^2}\big[\delta Q\delta F+\delta F\delta Q+\delta F-2\delta Q\textrm{tr}(\delta Q\delta F)\big]\cdot\nabla\delta u\,dx\nonumber\\
\underbrace{L\xi\int_{\mathbb{R}^2}\big[Q_2\delta \Delta Q+\delta \Delta Q Q_2]\cdot\nabla\delta u\,dx}_{\mathcal{BB}}+\xi\int_{\mathbb{R}^2}\big[Q_2\delta F+\delta F Q_2]\cdot\nabla\delta u\,dx&\nonumber\\
+\xi\int_{\mathbb{R}^2}\big[\delta QH_2+H_2\delta Q\big]\cdot\nabla\delta u\,dx&\nonumber\\
  +\underbrace{2\xi L\int_{\mathbb{R}^2}\big[\delta Q\textrm{tr}(Q_2\delta \Delta Q)+Q_2\textrm{tr}(\delta Q\delta \Delta Q)+Q_2\textrm{tr}( Q_2\delta \Delta Q )\big]\cdot\nabla \delta u\,dx}_{\mathcal{CC}}&\nonumber\\
  +2\xi \int_{\mathbb{R}^2}\big[\delta Q\textrm{tr}(Q_2\delta F)+Q_2\textrm{tr}(\delta Q\delta F)+Q_2\textrm{tr}( Q_2\delta F )\big]\cdot\nabla \delta u\,dx&\nonumber\\
+2\xi\int_{\mathbb{R}^2}\big[\delta Q\textrm{tr}(\delta QH_2)+Q_2\textrm{tr}(\delta Q H_2 )+\delta Q\textrm{tr}( Q_2 H_2 )\big]\cdot\nabla \delta u\,dx&\nonumber\\
-L\int_{\mathbb{R}^2} \Big(\delta Q\Delta Q_2-\Delta Q_2\delta Q\Big)\cdot\nabla\delta u\,dx-\underbrace{L\int_{\mathbb{R}^2} \Big(Q_2\Delta\delta Q-\Delta\delta Q Q_2\Big)\cdot\nabla\delta u\,dx}_{\mathcal{AA}}&
 \end{align}

      Let us observe that Lemma ~\ref{anulare} implies    $\mathcal{A}-\mathcal{AA}=0$ and that one can easily show $\mathcal{B}-\mathcal{BB}=0$ and $\mathcal{C}-\mathcal{CC}=0$. Also $\mathcal{I}+\mathcal{II}=0$ and then  we easily obtain
\begin{eqnarray}
\frac 12\frac{d}{dt}(L\|\nabla\delta Q\|^2_{L^2}+\|\delta Q\|_{L^2}^2+\|\delta u\|_{L^2}^2)+\Gamma L^2\|\Delta \delta Q\|_{L^2}^2+\nu\|\nabla \delta u\|_{L^2}^2\leq L\|\Delta\delta Q\|_{L^2}\|\delta u\|_{L^4}\|\nabla Q_2\|_{L^4}\nonumber\\
+L\|u_2\|_{L^4}\|\nabla\delta Q\|_{L^4}\|\Delta\delta Q\|_{L^2}+2(1+ |\xi|)L\|\delta Q\|_{L^{\frac{2}{s}}}\|\Omega_2\|_{L^{\frac{2}{1-s}}}\|\Delta\delta Q\|_{L^2}\nonumber\\
+2|\xi|L\|\delta Q\|_{L^{\frac{4}{s}}}^2\|\nabla u_2\|_{L^{\frac{2}{1-s}}}\|\Delta \delta Q\|_{L^2}+
4|\xi| L\|Q_2\|_{L^\infty}\|\delta Q\|_{L^{\frac{2}{s}}}\|\nabla u_2\|_{L^{\frac{2}{1-s}}}\|\Delta \delta Q\|_{L^2}+|\xi|\|\nabla u_2\|_{L^2}\|\delta Q\|_{L^4}^2\nonumber\\
+2|\xi|\|\delta Q\|_{L^{\frac{4}{s}}}^2\|\nabla u_2\|_{L^{\frac{2}{1-s}}}\|\delta Q\|_{L^2}+2|\xi|\|Q_2\|_{L^\infty}\|\delta Q\|_{L^4}^2\|\nabla\delta u\|_{L^2}+2|\xi|\|u_2\|_{L^4}\|\delta Q\|_{L^8}^2\|\nabla\delta u\|_{L^2}\nonumber\\
+4|\xi|\|Q_2\|_{L^\infty}\|\delta Q\|_{L^4}^2\|\nabla u_2\|_{L^2}+2|\xi|\|Q_2\|_{L^\infty}^2\|\delta Q\|_{L^2}\|\nabla\delta u\|_{L^2}
\nonumber\\
|a|\Gamma L\|\nabla\delta Q\|_{L^2}^2+|b|\Gamma L\|\Delta\delta Q\|_{L^2}\|\delta Q\|_{L^4} \|Q_1\|_{L^4} +|b|\Gamma L\|Q_2\|_{L^\infty}\|\delta Q \|_{L^2}\|\Delta\delta Q\|_{L^2}\nonumber\\
+c\Gamma L\|\delta Q\|_{L^4}\|\Delta\delta Q\|_{L^2}\|Q_1\|_{L^8}^2+c\Gamma L\|Q_2\|_{L^\infty}\|\Delta\delta Q\|_{L^2}\big(\|Q_1\|_{L^4}+\|Q_2\|_{L^4}\big)\|\delta Q\|_{L^4}\nonumber\\
+\|\nabla Q_2\|_{L^4}\|\delta u\|_{L^4}\|\delta Q\|_{L^2}+2(1+ |\xi|)\|Q_2\|_{L^\infty}\|\nabla\delta u\|_{L^2}\|\delta Q\|_{L^2}\nonumber
\end{eqnarray}
\begin{eqnarray}
+\xi\bigg(\|\delta Q\|_{L^2}\|\nabla\delta u\|_{L^2}+(\|Q_1\|_{L^4}+\|Q_2\|_{L^4})\|\delta Q\|_{L^4}\|\nabla\delta u\|_{L^2}+(\|Q_1\|_{L^6}^2+\|Q_2\|_{L^6}^2)\|\delta Q\|_{L^6}\|\nabla\delta u\|_{L^2}\bigg)\nonumber\\
+\xi\bigg(\sum_{j=1}^3(\|Q_1\|_{L^{4j}}^j+\|Q_2\|_{L^{4j}}^j\bigg)\|\delta Q\|_{L^4}\|\nabla\delta u\|_{L^2}+\sum_{j=1}^3(\|Q_1\|_{L^{6j}}^j+\|Q_2\|_{L^{6j}}^j)\|\delta Q\|_{L^6}^2\|\nabla\delta u\|_{L^2}\bigg)\nonumber\\
+4|\xi|\|Q_2\|_{L^\infty}\big(|a|\|Q_2\|_{L^\infty}+|b|\|Q_2\|_{L^\infty}^2+|c|\|Q_2\|_{L^\infty}^3\big)\|\delta Q\|_{L^2}\|\nabla\delta u\|_{L^2}\nonumber\\
+|a|\Gamma\|\delta Q\|_{L^2}^2+\Gamma\big(|b|+c\|Q_2\|_{L^\infty}\big)\big(\|Q_1\|_{L^2}+\|Q_2\|_{L^2}\big)\|\delta Q\|_{L^4}^2\nonumber\\
+2|\xi| \|Q_2\|_{L^\infty}\bigg(|a|\|\delta Q\|_{L^2}+\big(|b|(\|Q_1\|_{L^4}+\|Q_2\|_{L^4})+|c|(\| Q_1\|_{L^8}^2+\| Q_2\|_{L^8}^2+\|Q_1\|_{L^8}\|Q_2\|_{L^8}\big)\|\delta Q\|_{L^4}\bigg)\|\nabla\delta u\|_{L^2}\nonumber\\
+2|\xi|L\|\delta Q\|_{L^{\frac{2}{s}}}\|\Delta Q_2\|_{L^{\frac{2}{1-s}}}\|\nabla\delta u\|_{L^2}+2|\xi|\|Q_2\|_{L^\infty}\bigg(|a|\|\delta Q\|_{L^2}+\big(|b|\|Q_2\|_{L^4}+|c|\| Q_2\|_{L^\infty}\|Q_2\|_{L^4}\big)\|\delta Q\|_{L^4}\bigg)\|\nabla\delta u\|_{L^2}\nonumber\\
+2|\xi|\big(\|Q_2\|_{L^\infty}^2\|\delta Q\|_{L^4}\|\nabla\delta u\|_{L^2}+2\|Q_2\|_{L^\infty}\|\delta Q\|_{L^8}^2\|\nabla\delta u\|_{L^2}\big)\times\nonumber\\
\times\bigg[|a|+|b|(\|Q_1\|_{L^4}+\|Q_2\|_{L^4})+c(\|Q_1\|_{L^8}^2+\|Q_2\|_{L^8}^2+\|Q_1\|_{L^8}\|Q_2\|_{L^8})\bigg]\nonumber\\
2|\xi|L\|\delta Q\|_{L^4}^2\big(|a|\|Q_2\|_{L^\infty}+|b|\|Q_2\|_{L^\infty}^2+|c|\|Q_2\|_{L^\infty}^3\big)+4|\xi|L\|Q_2\|_{L^\infty}\|\Delta Q_2\|_{L^{\frac{2}{1-s}}}\|\delta Q\|_{L^{\frac{2}{s}}}\|\nabla \delta u\|_{L^2}\nonumber\\
+\|\delta u\|_{L^4}^2\|\nabla u_2\|_{L^2}+2L\|\nabla Q_2\|_{L^4}\|\nabla\delta Q\|_{L^4}\|\nabla \delta u\|_{L^2}+2L\|\Delta Q_2\|_{L^{\frac{2}{1-s}}}\|\delta Q\|_{L^{\frac{2}{s}}}\|\nabla\delta u\|_{L^2}\nonumber
\end{eqnarray}

Using that $\|\delta Q\|_{L^{\frac{2}{s}}}\leq \frac{C}{\sqrt{s}}\|\delta Q\|_{L^2}^{s}\|\nabla \delta Q\|_{L^2}^{1-s}$ and $\|\Omega_2\|_{L^{\frac{2}{1-s}}}\leq C\|u_2\|_{H^{1+s}}$, we obtain the estimate by
\begin{eqnarray}
\le\frac{\nu}{2}\|\nabla \delta u\|_{L^2}^2+\frac{\Gamma L^2}{2}\|\Delta\delta Q\|_{L^2}^2+C\underbrace{\Big(\|\nabla u_2\|_{L^{\frac{2}{1-s}}}^2+\|\nabla Q_2\|_{L^{\frac {2}{1-s}}}^2\Big)}_{\mathcal{J}_1}\|\delta u\|_{L^2}^2\nonumber\\+C\underbrace{\Big(1+\|Q_2\|_{L^\infty}^4+\|\nabla u_2\|_{L^\infty}^2+{    \|\nabla Q_2\|_{H^s}^2}+\|Q_2\|_{H^s}^2+\|\Delta Q_2\|_{H^s}^2\Big)}_{\mathcal{J}_2}\|\delta Q\|_{L^2}^2+C\underbrace{\Big({    1+}\|u_2\|_{L^\infty}^2+\|\nabla Q_2\|_{L^\infty}^2\Big)}_{\mathcal{J}_3}\|\nabla \delta Q\|_{L^2}^2\nonumber\\+C\underbrace{\bigg(\|Q_1\|_{L^4}^2+\|Q_1\|_{L^8}^4+{    \|Q_2\|_{L^\infty}^2}\big(\|Q_1\|_{L^4}^2+\|Q_2\|_{L^4}^2\big)+\Gamma\big(|b|+c\|Q_2\|_{L^\infty}\big)\big(\|Q_1\|_{L^2}+\|Q_2\|_{L^2}\big)\bigg)}_{\mathcal{J}_4}\|\delta Q\|_{L^4}^2\nonumber\\
\underbrace{C\big(1+\|\nabla u_2\|_{L^2}^2+\|Q_2\|_{L^\infty}^6+\|Q_2\|_{L^\infty}^2\|\nabla u_2\|_{L^2}^2+\|Q_2\|_{L^\infty}^2\|\delta Q\|_{L^4}^2\big)}_{\mathcal{J}_5}\|\delta Q\|_{L^4}^2\nonumber\\
\underbrace{C\big(1+\|Q_2\|_{L^\infty}^2+\|Q_2\|_{L^\infty}^4\big)\bigg(1+\sum_{j=1}^3(\|Q_1\|_{L^{4j}}^{2j}+\|Q_2\|_{L^{4j}}^{2j})\bigg)}_{\mathcal{J}_6}\|\delta Q\|_{L^4}^2+\underbrace{C\Big(1+\sum_{j=1}^3(\|Q_1\|_{L^{6j}}^{2j}+\|Q_2\|_{L^{6j}}^{2j})\Big)^2}_{\mathcal{J}_7}\|\delta Q\|_{L^6}^2    \nonumber\\
\underbrace{C\bigg(1+\|u_2\|_{L^4}^2\|\delta Q\|_{L^8}^2+\|\delta Q\|_{L^8}^2\big(\|Q_2\|_{L^\infty}^2+\|Q_2\|_{L^\infty}^4\big)\Big(1+\|Q_1\|_{L^4}^2+\|Q_2\|_{L^4}^2+\|Q_1\|_{L^8}^4+\|Q_2\|_{L^8}^4\Big)\bigg)}_{\mathcal{J}_8}\|\delta Q\|_{L^8}^2\nonumber\\
\underbrace{C\big(\|Q_2\|_{L^\infty}^2\|\nabla u_2\|_{L^{\frac{2}{1-s}}}^2+\|\Delta Q_2\|_{L^{\frac{2}{1-s}}}^2+\|Q_2\|_{L^\infty}^2\|\Delta Q_2\|_{L^{\frac{2}{1-s}}}^2\big)}_{\mathcal{J}_9}\|\delta Q\|_{L^{\frac{2}{s}}}^2+\underbrace{C\|\nabla u_2\|_{L^{\frac{2}{1-s}}}^2(1+\|\delta Q\|_{L^{\frac{4}{s}}}^2)}_{\mathcal{J}_{10}}\|\delta Q\|_{L^{\frac{4}{s}}}^2
\end{eqnarray}

      We are in $2D$ so $\|\delta Q\|_{L^4}^2, \|\delta Q\|_{L^6}^2, \|\delta Q\|_{L^8}^2, \|\delta Q\|_{L^{\frac{4}{s}}}^2,\|\delta Q\|_{L^{\frac{2}{s}}}^2$ are controlled by $\|\delta Q\|_{L^2}^2+\|\nabla \delta Q\|_{L^2}^2$. The hypothesis, namely relations (\ref{sol1}) and (\ref{sol2}), ensure that the terms $\mathcal{J}_i,i=1,\dots,10$ are integrable in time (choosing $\varepsilon>0$ sufficiently small, depending on $s$)  thus using the last inequality and  Gronwall Lemma we obtain the uniqueness of the solution. $\Box$

\bigskip\noindent{\bf Acknowledgements} MP and AZ thank John M. Ball for stimulating discussions. {    MP  gratefully acknowledges the hospitality of Oxford University's OXPDE Center. }AZ acknowledges the support of  the EPSRC Science and Innovation award to the Oxford Center for Nonlinear PDE (EP/E035027/1).

\appendix
\section{The correction terms}

\par For the Q-tensor equations we have the following correction terms:

\smallskip
\begin{eqnarray}
\mathcal{T}_Q\stackrel{\rm{def}}{=}-\Delta_q( u_\gamma Q_{\alpha\beta,\gamma})+\Gamma\Delta_q [-aQ_{\alpha\beta}+b\left(Q_{\alpha\gamma}Q_{\gamma\beta}-\frac{\delta_{\alpha\beta}}{2}\textrm{tr}(Q^2)\right)-cQ_{\alpha\beta}\textrm{tr}(Q^2)]\nonumber\\+\sum_{|q'-q|\le 5} [\Delta_q;S_{q'-1}Q_{\gamma\beta}]\Delta_{q'}\Omega_{\alpha\gamma}+\sum_{|q'-q|\le 5}(S_{q'-1}Q_{\gamma\beta}-S_{q-1}Q_{\gamma\beta})\Delta_q\Delta_{q'}\Omega_{\alpha\gamma}+
\sum_{q'>q-5} \Delta_q\left(S_{q'+2}\Omega_{\alpha\gamma} \Delta_{q'}Q_{\gamma\beta}\right)\nonumber\\
-\sum_{|q'-q|\le 5} [\Delta_q;S_{q'-1}Q_{\alpha\gamma}]\Delta_{q'}\Omega_{\gamma\beta}-\Sigma_{|q'-q|\le 5}(S_{q'-1}Q_{\alpha\gamma}-S_{q-1}Q_{\alpha\gamma})\Delta_q\Delta_{q'}\Omega_{\gamma\beta}-
\sum_{q'>q-5} \Delta_q\left(S_{q'+2}\Omega_{\gamma\beta}\Delta_{q'}Q_{\alpha\gamma}\right)\nonumber\\
+\xi\sum_{|q'-q|\le 5} [\Delta_q, S_{q'-1}Q_{\gamma\beta}]\Delta_{q'}D_{\alpha\gamma}+\xi\sum_{|q'-q|\le 5} (S_{q'-1}Q_{\gamma\beta}-S_{q-1}Q_{\gamma\beta})\Delta_q\Delta_{q'}D_{\alpha\gamma}+\xi\sum_{q'>q-5}\Delta_q(S_{q'+2}D_{\alpha\gamma}\Delta_{q'}Q_{\gamma\beta})\nonumber\\
+\xi \sum_{|q'-q|\le 5} [\Delta_q, S_{q'-1}Q_{\alpha\gamma}]\Delta_{q'}D_{\gamma\beta}+\xi\sum_{|q'-q|\le 5}(S_{q'-1} Q_{\alpha\gamma}-S_{q-1}Q_{\alpha\gamma})\Delta_q\Delta_{q'}D_{\gamma\beta}+\xi\sum_{q'>q-5} \Delta_q(S_{q'+2}D_{\gamma\beta}\Delta_{q'}Q_{\alpha\gamma})\nonumber\\
-2\xi\sum_{|q'-q|\le 5}[\Delta_q, S_{q'-1}Q_{\alpha\beta}]\Delta_{q'}\textrm{tr}(Q\nabla u)-2\xi\sum_{|q'-q|\le 5} \big(S_{q'-1}Q_{\alpha\beta}-S_{q-1}Q_{\alpha\beta}\big)\Delta_q\Delta_{q'}\textrm{tr}(Q\nabla u)\nonumber\\
-2\xi\sum_{q'>q-5}\Delta_q\Big(S_{q'+2}\textrm{tr}(Q\nabla u)\Delta_{q'}Q_{\alpha\beta}\Big)\nonumber\\
-2\xi  S_{q-1}Q_{\alpha\beta}\bigg(\sum_{|q'-q|\le 5} [\Delta_q,S_{q'-1}Q_{\gamma\delta}]\Delta_{q'}u_{\gamma,\delta}+\sum_{|q'-q|\le 5} (S_{q'-1}Q_{\gamma\delta}-S_{q-1}Q_{\gamma\delta})\Delta_q\Delta_{q'}u_{\gamma,\delta}\bigg)\nonumber\\-2\xi S_{q-1}Q_{\alpha\beta}\sum_{q'>q-5}\Delta_q (S_{q'+2}u_{\gamma,\delta}\Delta_{q'}Q_{\gamma\delta}\bigg)\nonumber
\end{eqnarray} 

\par Then we get:
\begin{eqnarray}
\Big(-L\Delta\Delta_Q Q_{\alpha\beta},\big(\mathcal{T}_Q\big)_{\alpha\beta}\Big)
 =L\underbrace{\left(\Delta_q(u\nabla Q_{\alpha\beta}),\Delta\Delta_q Q_{\alpha\beta}\right)}_{\stackrel{def}{=}\mathcal{I}_1}-L\underbrace{\Sigma_{|q'-q|\le 5} \left([\Delta_q;S_{q'-1}Q_{\gamma\beta}]\Delta_{q'}\Omega_{\alpha\gamma},\Delta\Delta_q Q_{\alpha\beta}\right)}_{\stackrel{def}{=}\mathcal{I}_2}\nonumber\\-L\underbrace{\Sigma_{|q'-q|\le 5}\left((S_{q'-1}Q_{\gamma\beta}-S_{q-1}Q_{\gamma\beta})\Delta_q\Delta_{q'}\Omega_{\alpha\gamma},\Delta\Delta_q Q_{\alpha\beta}\right)}_{\stackrel{def}{=}\mathcal{I}_3}\nonumber
\end{eqnarray}
\begin{eqnarray}
-L\underbrace{\Sigma_{q'>q-5} \left( \Delta_q\left(S_{q'+2}\Omega_{\alpha\gamma} \Delta_{q'}Q_{\gamma\beta}\right),\Delta\Delta_q Q_{\alpha\beta}\right)}_{\stackrel{def}{=}\mathcal{I}_4}+L\underbrace{\Sigma_{|q'-q|\le 5}\left( [\Delta_q;S_{q'-1}Q_{\alpha\gamma}]\Delta_{q'}\Omega_{\gamma\beta},\Delta\Delta_q Q_{\alpha\beta}\right)}_{\stackrel{def}{=}\mathcal{I}_5}\nonumber\\
+L\underbrace{\Sigma_{|q'-q|\le 5}\left((S_{q'-1}Q_{\alpha\gamma}-S_{q-1}Q_{\alpha\gamma})\Delta_q\Delta_{q'}\Omega_{\gamma\beta},\Delta\Delta_q Q_{\alpha\beta}\right)}_{\stackrel{def}{=}\mathcal{I}_6}+L\underbrace{\Sigma_{q'>q-5}\left(  \Delta_q\left(S_{q'+2}\Omega_{\gamma\beta}\Delta_{q'}Q_{\alpha\gamma})\right),\Delta\Delta_q Q_{\alpha\beta}\right)}_{\stackrel{def}{=}\mathcal{I}_7}\nonumber
\end{eqnarray}
\begin{eqnarray}
-L\Gamma\underbrace{\Big(\Delta_q [-aQ_{\alpha\beta}+bQ_{\alpha\gamma}Q_{\gamma\beta}-cQ_{\alpha\beta}\textrm{tr}(Q^2)],\Delta\Delta_q Q_{\alpha\beta}\Big)}_{\stackrel{def}{=}\mathcal{I}_8}\nonumber\\
-L\xi\underbrace{\sum_{|q'-q|\le 5}\Big( [\Delta_q, S_{q'-1}Q_{\gamma\beta}]\Delta_{q'}D_{\alpha\gamma},\Delta\Delta_q Q_{\alpha\beta}\Big)}_{\stackrel{def}{=}\mathcal{I}_9}-L\xi\underbrace{\sum_{|q'-q|\le 5} \Big( (S_{q'-1}Q_{\gamma\beta}-S_{q-1}Q_{\gamma\beta})\Delta_q\Delta_{q'}D_{\alpha\gamma},\Delta\Delta_q Q_{\alpha\beta}\Big)}_{\stackrel{def}{=}\mathcal{I}_{10}}\nonumber\\-L\xi\underbrace{\sum_{q'>q-5}\Big(\Delta_q(S_{q'+2}D_{\alpha\gamma}\Delta_{q'}Q_{\gamma\beta}),\Delta\Delta_q Q_{\alpha\beta}\Big)}_{\stackrel{def}{=}\mathcal{I}_{11}}
-L\xi \underbrace{\sum_{|q'-q|\le 5} \Big([\Delta_q, S_{q'-1}Q_{\alpha\gamma}]\Delta_{q'}D_{\gamma\beta},\Delta_q\Delta Q_{\alpha\beta}\Big)}_{\stackrel{def}{=}\mathcal{I}_{12}}\nonumber
\end{eqnarray}
\begin{eqnarray}
-L\xi\underbrace{\sum_{|q'-q|\le 5}(S_{q'-1} Q_{\alpha\gamma}-S_{q-1}Q_{\alpha\gamma})\Delta_q\Delta_{q'}D_{\gamma\beta},\Delta_q\Delta Q_{\alpha\beta}\Big)}_{\stackrel{def}{=}\mathcal{I}_{13}}-L\xi\underbrace{\sum_{q'>q-5} \Big(\Delta_q(S_{q'+2}D_{\gamma\beta}\Delta_{q'}Q_{\alpha\gamma}),\Delta_q\Delta Q_{\alpha\beta}\Big)}_{\stackrel{def}{=}\mathcal{I}_{14}}\nonumber\\
+2L\xi\underbrace{\sum_{|q'-q|\le 5}\Big([\Delta_q, S_{q'-1}Q_{\alpha\beta}]\Delta_{q'}\textrm{tr}(Q\nabla u),\Delta_q\Delta Q_{\alpha\beta}\Big)}_{\stackrel{def}{=}\mathcal{I}_{15}}+2L\xi\underbrace{\sum_{|q'-q|\le 5}\Big( \big(S_{q'-1} Q_{\alpha\beta}-S_{q-1}Q_{\alpha\beta}\big)\Delta_q\Delta_{q'}\textrm{tr}(Q\nabla u),\Delta\Delta_q Q_{\alpha\beta}\Big)}_{\stackrel{def}{=}\mathcal{I}_{16}}\nonumber\\
+2L\xi\underbrace{\sum_{q'>q-5}\Big(\Delta_q(S_{q'+2}\textrm{tr}(Q\nabla u)\Delta_{q'}Q_{\alpha\beta}),\Delta_q \Delta Q_{\alpha\beta}\Big)}_{\stackrel{def}{=}\mathcal{I}_{17}}
+2L\xi\underbrace{\bigg( S_{q-1}Q_{\alpha\beta}\big(\sum_{|q'-q|\le 5} [\Delta_q,S_{q'-1}Q_{\gamma\delta}]\Delta_{q'}u_{\gamma,\delta}\big),\Delta_q\Delta Q_{\alpha\beta}\bigg)}_{\stackrel{def}{=}\mathcal{I}_{18}}\nonumber\\+2L\xi\underbrace{\bigg(S_{q-1}Q_{\alpha\beta}\big(\sum_{|q'-q|\le 5} (S_{q'-1}Q_{\gamma\delta}-S_{q-1}Q_{\gamma\delta})\Delta_q\Delta_{q'}u_{\gamma,\delta}\big),\Delta\Delta_q Q_{\alpha\beta}\bigg)}_{\stackrel{def}{=}\mathcal{I}_{19}}\nonumber\\+2L\xi\underbrace{\Big(S_{q-1}Q_{\alpha\beta}\sum_{q'>q-5}\Delta_q (S_{q'+2}u_{\gamma,\delta}\Delta_{q'}Q_{\gamma\delta}),\Delta_q\Delta Q_{\alpha\beta}\Big)}_{\stackrel{def}{=}\mathcal{I}_{20}}\nonumber
\end{eqnarray}

\par The correction terms for the Navier-Stokes part are:

\begin{eqnarray}
\big(\mathcal{T}_u\big)_\alpha\stackrel{\rm{def}}{=}-L   \partial_\beta\Delta_q\left(   \partial_\alpha Q_{\gamma\delta}   \partial_\beta Q_{\gamma\delta}-\frac{\delta_{\alpha\beta}}{3}   \partial_\lambda Q_{\gamma\delta}   \partial_\lambda Q_{\gamma\delta}\right)-\xi \Delta_q F_{\alpha\beta,\beta}\nonumber\\
-\Delta_q(u_\beta   \partial_\beta u_\alpha)-\xi   \partial_\beta\Big(\Delta_q\big(Q_{\alpha\gamma}F_{\gamma\beta}\big)+\Delta_q\big(F_{\alpha\gamma}Q_{\gamma\beta}\big)-2\Delta_q\big(Q_{\alpha\beta}\textrm{tr}(QF)\big)\Big)\nonumber\\
+L   \partial_\beta\Big(\Sigma_{|q'-q|\le 5}[\Delta_q;S_{q'-1} Q_{\alpha\gamma}]\Delta_{q'}\Delta Q_{\gamma\beta}+\Sigma_{|q'-q|\le 5}(S_{q'-1}Q_{\alpha\gamma}-S_{q-1}Q_{\alpha\gamma})\Delta_q\Delta_{q'} \Delta Q_{\gamma\beta}\Big)\nonumber\\
+L   \partial_\beta\Big(\Sigma_{q'> q-5}\Delta_q(S_{q'+2}\Delta Q_{\gamma\beta}\Delta_{q'} Q_{\alpha\gamma})
-\Sigma_{|q'-q|\le 5}[\Delta_q;S_{q'-1}  Q_{\gamma\beta}]\Delta_{q'}\Delta Q_{\alpha\gamma}\Big)\nonumber\\
-L   \partial_\beta\Big(\Sigma_{|q'-q|\le 5}(S_{q'-1} Q_{\gamma\beta}-S_{q-1} Q_{\gamma\beta})\Delta_q\Delta_{q'}\Delta Q_{\alpha\gamma}{    +}\Sigma_{q'> q-5}\Delta_q(S_{q'+2}\Delta Q_{\alpha\gamma}\Delta_{q'}  Q_{\gamma\beta})\Big)\nonumber
\end{eqnarray}
\begin{eqnarray}
-L\xi   \partial_\beta\bigg(\sum_{|q'-q|\le 5} [\Delta_q, S_{q'-1} Q_{\alpha\gamma}]\Delta_{q'}\Delta Q_{\gamma\beta}+\sum_{|q'-q|\le 5}(S_{q'-1} Q_{\alpha\gamma}-S_{q-1}Q_{\alpha\gamma})\Delta_q\Delta_{q'}Q_{\gamma\beta}\bigg)\nonumber\\
-L\xi   \partial_\beta\bigg(\sum_{q'>q-5} \Delta_q(S_{q'+2}\Delta Q_{\gamma\beta}\Delta_{q'}Q_{\alpha\gamma})+\sum_{|q'-q|\le 5} [\Delta_q, S_{q'-1} Q_{\gamma\beta} ]\Delta_{q'}\Delta Q_{\alpha\gamma}\bigg)\nonumber\\
-L\xi   \partial_\beta\bigg(\sum_{|q'-q|\le 5} (S_{q'-1}Q_{\gamma\beta}-S_{q-1}Q_{\gamma\beta})\Delta_q\Delta_{q'}\Delta Q_{\alpha\gamma}+\sum_{q'>q-5}\Delta_q (S_{q'+2}\Delta Q_{\alpha\gamma}\Delta_{q'} Q_{\gamma\beta})\bigg)\nonumber\\
+2L\xi   \partial_\beta\bigg( \sum_{|q'-q|\le 5} [\Delta_q, S_{q'-1} Q_{\alpha\beta}] \Delta_{q'}\textrm{tr}(Q\Delta Q)+ \sum_{|q'-q|\le 5} (S_{q'-1}Q_{\alpha\beta}-S_{q-1}Q_{\alpha\beta})\Delta_q\Delta_{q'}\textrm{tr}(Q\Delta Q)\bigg)\nonumber\\
+2L\xi   \partial_\beta\bigg(\sum_{q'>q-5} \Delta_q(S_{q'+2}\textrm{tr}(Q\Delta Q)\Delta_{q'}Q_{\alpha\beta})\bigg)\nonumber\\
+2L\xi   \partial_\beta\bigg(S_{q-1}Q_{\alpha\beta}\sum_{|q'-q|\le 5}[\Delta_q, S_{q'-1}Q_{\gamma\delta}]\Delta_{q'}\Delta Q_{\gamma\delta}+S_{q-1}Q_{\alpha\beta}\sum_{|q'-q|\le 5}(S_{q'-1}Q_{\gamma\delta}-S_{q-1}Q_{\gamma\delta})\Delta_q\Delta_{q'}\Delta Q_{\gamma\delta}\bigg)\nonumber\\+2L\xi   \partial_\beta\big(S_{q-1}Q_{\alpha\beta}\sum_{q'>q-5}\Delta_q(S_{q'+2}\Delta Q_{\gamma\delta}\Delta_{q'}Q_{\gamma\delta}\big)
\end{eqnarray}

\begin{eqnarray}
\Big(\big(\mathcal{T}_u\big)_\alpha,\Delta_q u_\alpha\Big)={    -}\underbrace{\left(\Delta_q(u_\beta   \partial_\beta u_\alpha),\Delta_q u_\alpha\right)}_{\stackrel{def}{=}\mathcal{J}_1}
+L\underbrace{\int \Delta_q\left(   \partial_\alpha Q_{\gamma\delta}   \partial_\beta Q_{\gamma\delta}-\frac{\delta_{\alpha\beta}}{3}   \partial_\lambda Q_{\gamma\delta}   \partial_\lambda Q_{\gamma\delta}\right)\Delta_q u_{\alpha,\beta}}_{\stackrel{def}{=}\mathcal{J}_2}\nonumber\\
-L\underbrace{\Sigma_{|q'-q|\le 5}\int[\Delta_q;S_{q'-1} Q_{\alpha\gamma}]\Delta_{q'}\Delta Q_{\gamma\beta}\Delta_q u_{\alpha,\beta}}_{\stackrel{def}{=}\mathcal{J}_3}
-L\underbrace{\int\Sigma_{|q'-q|\le 5}(S_{q'-1}Q_{\alpha\gamma}-S_{q-1}Q_{\alpha\gamma})\Delta_q\Delta_{q'}\Delta Q_{\gamma\beta}\Delta_q u_{\alpha,\beta}}_{\stackrel{def}{=}\mathcal{J}_4}\nonumber
\end{eqnarray}
\begin{eqnarray}-L\underbrace{\int\Sigma_{q'> q-5}\Delta_q(S_{q'+2} \Delta Q_{\gamma\beta}\Delta_{q'} Q_{\alpha\gamma})\Delta_q u_{\alpha,\beta}}_{\stackrel{def}{=}\mathcal{J}_5}+L\underbrace{\Sigma_{|q'-q|\le 5}\int [\Delta_q;S_{q'-1}  Q_{\gamma\beta}]\Delta_{q'} \Delta Q_{\alpha\gamma}\Delta_q u_{\alpha,\beta}}_{\stackrel{def}{=}\mathcal{J}_6}\nonumber\\
+L\underbrace{\int \Sigma_{|q'-q|\le 5}(S_{q'-1} Q_{\gamma\beta}-S_{q-1} Q_{\gamma\beta})\Delta_q\Delta_{q'}\Delta Q_{\alpha\gamma}\Delta_q u_{\alpha,\beta}}_{\stackrel{def}{=}\mathcal{J}_7}+L\underbrace{\int \Sigma_{q'>q-5}\Delta_q(S_{q'+2}\Delta Q_{\alpha\gamma}\Delta_{q'} Q_{\gamma\beta})\Delta_q u_{\alpha,\beta}}_{\stackrel{def}{=}\mathcal{J}_8}\nonumber\\
+L\xi\bigg(\underbrace{\big(\sum_{|q'-q|\le 5} [\Delta_q, S_{q'-1} Q_{\alpha\gamma}]\Delta_{q'}\Delta Q_{\gamma\beta}, \Delta_q u_{\alpha,\beta}\big)}_{\stackrel{def}{=}\mathcal{J}_9}+\underbrace{\big(\sum_{|q'-q|\le 5}(S_{q'-1} Q_{\alpha\gamma}-S_{q-1}Q_{\alpha\gamma})\Delta_q\Delta_{q'}\Delta Q_{\gamma\beta},\Delta_q u_{\alpha,\beta}\big)}_{\stackrel{def}{=}\mathcal{J}_{10}}\bigg)\nonumber
\end{eqnarray}
\begin{eqnarray}
+L\xi\bigg(\underbrace{\big(\sum_{q'>q-5} \Delta_q(S_{q'+2}\Delta Q_{\gamma\beta}\Delta_{q'}Q_{\alpha\gamma}), \Delta_q u_{\alpha,\beta}\big)}_{\stackrel{def}{=}\mathcal{J}_{11}}+\underbrace{\big(\sum_{|q'-q|\le 5} [\Delta_q, S_{q'-1} Q_{\gamma\beta} ]\Delta_{q'}\Delta Q_{\alpha\gamma},\Delta_q u_{\alpha,\beta}\big)}_{\stackrel{def}{=}\mathcal{J}_{12}}\bigg)\nonumber\\
+L\xi\bigg(\underbrace{\big(\sum_{|q'-q|\le 5} (S_{q'-1}Q_{\gamma\beta}-S_{q-1}Q_{\gamma\beta})\Delta_q\Delta_{q'}\Delta Q_{\alpha\gamma},\Delta_q u_{\alpha,\beta}\big)}_{\stackrel{def}{=}\mathcal{J}_{13}}+\underbrace{\big(\sum_{q'>q-5}\Delta_q (S_{q'+2}\Delta Q_{\alpha\gamma}\Delta_{q'} Q_{\gamma\beta}),\Delta_q u_{\alpha,\beta}\big)}_{\stackrel{def}{=}\mathcal{J}_{14}}\bigg)\nonumber
\end{eqnarray}
\begin{eqnarray}
-2L\xi\bigg(\underbrace{\big( \sum_{|q'-q|\le 5} [\Delta_q, S_{q'-1} Q_{\alpha\beta}] \Delta_{q'}\textrm{tr}(Q\Delta Q),\Delta_q u_{\alpha,\beta}\big)}_{\stackrel{def}{=}\mathcal{J}_{15}}+ \underbrace{\big(\sum_{|q'-q|\le 5} (S_{q'-1}Q_{\alpha\beta}-S_{q-1}Q_{\alpha\beta})\Delta_q\Delta_{q'}\textrm{tr}(Q\Delta Q),\Delta_q u_{\alpha,\beta}\big)}_{\stackrel{def}{=}\mathcal{J}_{16}}\bigg)\nonumber\\
-2L\xi\bigg(\underbrace{\big(\sum_{q'>q-5} \Delta_q(S_{q'+2}\textrm{tr}(Q\Delta Q)\Delta_{q'}Q_{\alpha\beta}),\Delta_q u_{\alpha,\beta}\big)}_{\stackrel{def}{=}\mathcal{J}_{17}}\bigg)
-2L\xi\underbrace{\bigg(S_{q-1}Q_{\alpha\beta}\sum_{|q'-q|\le 5}[\Delta_q, S_{q'-1}Q_{\gamma\delta}]\Delta_{q'}\Delta Q_{\gamma\delta},\Delta_q u_{\alpha,\beta}\bigg)}_{\stackrel{def}{=}\mathcal{J}_{18}}\nonumber\\-2L\xi\underbrace{\bigg(S_{q-1}Q_{\alpha\beta}\sum_{|q'-q|\le 5}(S_{q'-1}Q_{\gamma\delta}-S_{q-1}Q_{\gamma\delta})\Delta_q\Delta_{q'}\Delta Q_{\gamma\delta},\Delta_q u_{\alpha,\beta}\bigg)}_{\stackrel{def}{=}\mathcal{J}_{19}}\nonumber\\-2L\xi\underbrace{\bigg(S_{q-1}Q_{\alpha\beta}\sum_{q'>q-5}\Delta_q(S_{q'+2}\Delta Q_{\gamma\delta}\Delta_{q'}Q_{\gamma\delta}),\Delta_q u_{\alpha,\beta}\bigg)}_{\stackrel{def}{=}\mathcal{J}_{20}}\nonumber\\
-\xi\Big(\underbrace{\int\Delta_q\big(Q_{\alpha\gamma}F_{\gamma\beta}\big)\Delta_q u_{\alpha,\beta}+\int\Delta_q\big(F_{\alpha\gamma}Q_{\gamma\beta}\big)\Delta_q u_{\alpha,\beta}}_{\stackrel{def}{=}\mathcal{J}_{21}}-2\underbrace{\int\Delta_q\big(Q_{\alpha\beta}\textrm{tr}(QF)\big)\Delta_q u_{\alpha,\beta}\Big)}_{\stackrel{def}{=}\mathcal{J}_{22}}+\xi\underbrace{\int\Delta_q F_{\alpha\beta}\Delta_q u_{\alpha,\beta}}_{\stackrel{def}{=}\mathcal{J}_{23}}\nonumber
\end{eqnarray}

\section{Proof of estimate (\ref{longest})}
      In the following, $a_q(t)$  denotes a sequence in $l_{q}^2$  for all $t>0$ and $b_q(t)$ is a sequence in  $l^1_q$, $\forall t\ge 0$, sequences that can change from one line to the next. Moreover{     $\|\big(a_q(t)\big)_{q\in\mathbb{N}}\|_{l^2},\|\big(b_q(t)\big)_{q\in\mathbb{N}}\|_{l^1}\le C$} where the constant $C$ is independent of $t\ge 0$.
\begin{eqnarray}
|\mathcal{I}_1|=|\left(\Delta_q(u\nabla Q_{\alpha\beta}),\Delta_q\Delta Q_{\alpha\beta}\right)|\stackrel{(\ref{bonydecomp})}{=}|\underbrace{\int S_{q-1} u \Delta_q\nabla Q_{\alpha\beta}\Delta_q\Delta Q_{\alpha\beta}}_{\stackrel{\rm{def}}{=}\mathcal{I}_{1a}}
+\underbrace{\sum_{|q'-q|\le 5}\left([\Delta_q;S_{q'-1} u]\Delta_{q'}\nabla Q_{\alpha\beta},\Delta_q\Delta Q_{\alpha\beta}\right)}_{\stackrel{\rm{def}}{=}\mathcal{I}_{1b}}\nonumber\\
    +\underbrace{\sum_{|q'-q|\le 5}\left((S_{q'-1} u-S_{q-1} u)\Delta_q\Delta_{q'}\nabla Q_{\alpha\beta},\Delta\Delta_q Q_{\alpha\beta}\right)}_{\stackrel{\rm{def}}{=}\mathcal{I}_{1c}}+\underbrace{\sum_{q'\ge q-5}(\Delta_q(S_{q'+2}\nabla Q_{\alpha\beta}\Delta_{q'} u),\Delta_q\Delta Q_{\alpha\beta})}_{\stackrel{\rm{def}}{=}\mathcal{I}_{1d}}|\nonumber
\end{eqnarray}
We will use frequently  interpolation inequalities such as
$$\|f\|_{L^4(R^2)}\leq C\|f\|_{L^2(R^2)}^{\frac 12}\|\nabla f\|_{L^2(R^2)}^{\frac 12}.$$
We have
\begin{eqnarray}
|\mathcal{I}_{1a}|\le C \|u\|_{L^4}\|\Delta_q\nabla Q\|_{L^4}\|\Delta\Delta_q Q\|_{L^2}\le C2^{-2qs}b_q(t) \|u\|_{L^2}^{\frac 12}\|\nabla u\|_{L^2}^{\frac 12}\|\nabla Q\|_{H^s}^{\frac 12}\|\Delta Q\|_{H^s}^{\frac 32}\nonumber
\end{eqnarray}
On the other hand, using the commutator estimates and the  Bernstein inequality from Lemma ~\ref{lemma:bernstein&commutator}  we have
\begin{eqnarray}
    |\mathcal{I}_{1b}|\le \sum_{|q'-q|\le 5}\| [\Delta_q;S_{q'-1} u]\Delta_{q'}\nabla Q_{\alpha\beta}\|_{L^2}\|\Delta_q \Delta Q_{\alpha\beta}\|_{L^2}\le \sum_{|q'-q|\le 5} 2^{-q}\|\nabla S_{q'-1}u\|_{L^\infty}\|\nabla\Delta_{q'} Q_{\alpha\beta}\|_{L^2}\|\Delta_q\Delta Q_{\alpha\beta}\|_{L^2}\nonumber\\
    \leq C\sum_{|q'-q|\leq 5}2^{\frac{q'}{2}}\|S_{q'-1}u\|_{L^4}\|\nabla\Delta_{q'} Q_{\alpha\beta}\|_{L^2}\|\Delta_q\Delta Q_{\alpha\beta}\|_{L^2}
 \leq   C\|u\|_{L^4}b_q 2^{-2qs}\|\nabla Q_{\alpha\beta}\|_{H^{s+\frac 12}}\|\Delta Q_{\alpha\beta}\|_{H^s}\nonumber\\
    \le C\|u\|_{L^2}^{\frac 12}\|\nabla u\|_{L^2}^{\frac 12} 2^{-2qs}b_q(t)\|\nabla Q\|_{H^s}^{\frac 12}\|\Delta Q\|_{H^s}^{\frac 32}\nonumber
\end{eqnarray}

\begin{eqnarray}
|\mathcal{I}_{1c}|\le C\|u\|_{L^4}\|\Delta_q\nabla Q\|_{L^4}\|\Delta\Delta_q Q\|_{L^2}\le C2^{-2qs}b_q(t) \|u\|_{L^2}^{\frac 12}\|\nabla u\|_{L^2}^{\frac 12}\|\nabla Q\|_{H^s}^{\frac 12}\|\Delta Q\|_{H^s}^{\frac 32}\nonumber
\end{eqnarray}

\begin{eqnarray}
|\mathcal{I}_{1d}|\le\sum_{q'> q-5}|(\Delta_q(S_{q'+2}\nabla Q_{\alpha\beta}\Delta_{q'} u),\Delta_q\Delta Q_{\alpha\beta})|
\le \|\nabla Q\|_{L^4}\sum_{q'> q-5}2^{-(q'+q)s}2^{q's}\|\Delta_{q'}u\|_{L^4}2^{qs}\|\Delta_q\Delta Q\|_{L^2}\nonumber\\
    \le\|\nabla Q\|_{L^2}^{\frac 12}\|\Delta Q\|_{L^2}^{\frac 12} \sum_{q'> q-5}2^{-(q'+q)s}a_{q'}(t)\bar a_q\|u\|_{H^s}^{\frac 12}\|\nabla u\|_{L^2}^{\frac 12}\|\Delta Q\|_{H^s}\le
 C\|\nabla Q\|_{L^2}^{\frac 12}\|\Delta Q\|_{L^2}^{\frac 12} 2^{-2qs} b_q(t)\|u\|_{H^s}^{\frac{1}{2}}\|\nabla u\|_{H^s}^{\frac {1}{2}}\|\Delta Q\|_{H^s}\nonumber
\end{eqnarray}{     where $ b_q(t)=\Big(\sum_{q'>q-5} 2^{-(q'-q)s}a_{q'}(t)\Big)\bar a_{q}(t)$.} 

\begin{eqnarray}
|\mathcal{I}_2|=|\sum_{|q'-q|\le 5} \left([\Delta_q;S_{q'-1}Q_{\gamma\beta}]\Delta_{q'}\Omega_{\alpha\gamma},\Delta\Delta_q Q_{\alpha\beta}\right)|\le \sum_{|q'-q|\le 5}  2^{-q}\|S_{q'-1}\nabla Q_{\gamma\beta}\|_{L^\infty}\|\Delta_{q'} \Omega_{\alpha\gamma}\|_{L^2}\|\Delta\Delta_q Q_{\alpha\beta}\|_{L^2}\nonumber\\
\le \sum_{|q'-q|\le 5}C 2^{-q}2^{\frac {q'}{2}}\|S_{q'-1}\nabla Q_{\gamma\beta}\|_{L^4}2^{q'}\|\Delta_{q'}u\|_{L^2}\|\Delta\Delta_q Q_{\alpha\beta}\|_{L^2}
\le  C\sum_{|q'-q|\le 5}\|\nabla Q\|_{L^2}^{\frac 12}\|\Delta Q\|_{L^2}^{\frac 12}\|\Delta_{q'} u\|_{L^2}^{\frac 12}\|\Delta_{q'}\nabla u\|_{L^2}^{\frac 12}\|\Delta\Delta_q Q\|_{L^2}\nonumber\\
\le  C 2^{-2qs} b_q(t)\|\nabla Q\|_{L^2}^{\frac 12}\|\Delta Q\|^{\frac 12}_{L^2}\|u\|_{H^s}^{\frac 12}\|\nabla u\|_{H^s}^{\frac 12}\|\Delta Q\|_{H^s}\nonumber
\end{eqnarray}

\begin{eqnarray}
|\mathcal{I}_3|=|\sum_{|q'-q|\le 5}\left((S_{q'-1}Q_{\gamma\beta}-S_{q-1}Q_{\gamma\beta})\Delta_q\Delta_{q'}\Omega_{\alpha\gamma}, \Delta\Delta_q Q_{\alpha\beta}\right)|
\le \sum_{|q'-q|\le 5} \|\left(S_{q'-1}Q_{\gamma\beta}-S_{q-1}Q_{\gamma\beta}\right)\Delta_q\Delta_{q'}\Omega_{\alpha\gamma}\|_{L^2}\| \Delta\Delta_q Q_{\alpha\beta}\|_{L^2}\nonumber\\
\le C\sum_{|q'-q|\le 5}\|S_{q'-1}Q_{\gamma\beta}-S_{q-1}Q_{\gamma\beta}\|_{L^4}\|\Delta_q\Omega_{\alpha\gamma}\|_{L^4}
\| \Delta\Delta_q Q_{\alpha\beta}\|_{L^2}
\le \sum_{|q'-q|\le 5}2^{-q'}\|\tilde\Delta_{q'}\nabla Q_{\gamma\beta}\|_{L^4}2^{q}\|\Delta_q u\|_{L^4}\| \Delta\Delta_q Q_{\alpha\beta}\|_{L^2}\nonumber\\
\le C\|\nabla Q\|_{L^4}\|\Delta_q u\|_{L^2}^{\frac 12}\|\Delta_q \nabla u\|_{L^2}^{\frac 12}\|\Delta_q\Delta Q\|_{L^2}\le C2^{-2qs}b_q(t)\|\nabla Q\|_{L^2}^{\frac 12}\|\Delta Q\|_{L^2}^{\frac 12}\|u\|_{H^s}^{\frac 12}\|\nabla u\|_{H^s}^{\frac 12}\|\Delta Q\|_{H^s}\nonumber
\end{eqnarray} where $\tilde\Delta_q=\sum_{|i|\le 5}\Delta_q$.

\begin{eqnarray}
|\mathcal{I}_4|=|\sum_{q'> q-5}\left(\Delta_q\left(S_{q'+2}\Omega_{\alpha\gamma}\Delta_{q'}Q_{\gamma\beta}\right),\Delta_q \Delta Q_{\alpha\beta}\right)|\le \sum_{q'> q-5}\|\Delta_q(S_{q'+2}\Omega_{\alpha\gamma}\Delta_{q'}Q_{\gamma\beta})\|_{L^2}\|\Delta_q\Delta Q_{\alpha\beta}\|_{L^2}\nonumber\\
\le \sum_{q'> q-5}\|S_{q'+2}\Omega_{\alpha\gamma}\|_{L^4}\|\Delta_{q'}Q_{\gamma\beta}\|_{L^4}\|\Delta_q\Delta Q_{\alpha\beta}\|_{L^2}\le  \sum_{q'> q-5} C2^{q'}\|S_{q'+2} u\|_{L^4}\|\Delta_{q'}Q_{\gamma\beta}\|_{L^4} \|\Delta_q \Delta Q_{\alpha\beta}\|_{L^2}\nonumber\\
    \le \|u\|_{L^2}^{\frac 12}\|\nabla u\|_{L^2}^{\frac 12}C\sum_{q'> q-5} 2^{-q's-qs}2^{q'(s+1)}\|\Delta_{q'}Q_{\gamma\beta}\|_{L^2}^{\frac 12}\|\Delta_{q'} \nabla Q_{\gamma\beta}\|_{L^2}^{\frac 12}2^{qs}\|\Delta_q\Delta Q_{\alpha\beta}\|_{L^2}\nonumber\\
    \le C\|u\|_{L^2}^{\frac 12}\|\nabla u\|_{L^2}^{\frac 12}2^{-2qs} b_q(t)\|\nabla Q\|_{H^s}^{\frac 12}\|\Delta Q\|_{H^s}^{\frac 32}\nonumber
\end{eqnarray} where  $b_q(t)=\Big(\sum_{q'>q-5} 2^{-(q'-q)s}a_{q'}(t)\Big)\bar a_{q}(t)$.

\smallskip          The term $\mathcal{I}_k, k=5,6,7$ is estimated exactly as the term $\mathcal{I}_{k-3}$ that we have  already studied above.

      We claim that:
\begin{eqnarray}
|\mathcal{I}_8|\le 2^{-2qs}b_q(t) \bigg[C\bigg(1+\big(\sum_{j=2}^3 \|Q\|_{L^{2(j-1)}}^{j-1}\big)^2\bigg)\|\nabla Q\|_{H^s}^2+\frac{\Gamma L^2}{100}\|\Delta Q\|_{H^s}^2\bigg]
\label{i8}
\end{eqnarray}

      In order to prove the above estimate, we observe that the simplest terms are those of the form $(\Delta_q Q_{\alpha\beta},\Delta_q \Delta Q_{\alpha\beta})$ that can be easily estimated:

\begin{equation}
|(\Delta_q Q_{\alpha\beta},\Delta_q \Delta Q_{\alpha\beta})|\le 2^{-2qs}b_q(t)\|\nabla Q\|_{H^s}^2
\label{est:easiest}
\end{equation}

      For the rest of the terms we just consider a generic term from $\mathcal{I}_8$, namely $(\Delta_q(Q_{11}^j),\Delta_q\Delta Q_{\alpha\beta})$ where $2\le j\le 3$. We prove first the following:
      
 \begin{lemma}
 We have:
 \begin{equation}
 \|\Delta_q(Q_{11}^j)\|_{L^p}\le 2^{-qs}a_q(t)\|Q_{11}\|_{L^{p(j-1})}^{j-1}\|\nabla Q\|_{H^s} 
 \label{deltaqpowers}
 \end{equation} for $j\ge 2$.
 \label{lemma:deltaqpowers}
 \end{lemma}

\smallskip\par{\bf Proof.} We prove the statement by induction.

\par{\it Step $1$} We have:

$$\Delta_q(Q_{11}^2)=\sum_{q'>q-5}\Delta_q(S_{q'+2}Q_{11}\Delta_{q'}Q_{11})+\sum_{|q'-q|\le 5}\Delta_q (S_{q'-1}Q_{11}\Delta_{q'}Q_{11})$$ and

\begin{eqnarray}\|\sum_{q'>q-5}\Delta_q(S_{q'+2}Q_{11}\Delta_{q'}Q_{11})\|_{L^p}\le \|Q_{11}\|_{L^p}\sum_{q'>q-5}\|\Delta_{q'}Q_{11}\|_{L^\infty}\le \|Q_{11}\|_{L^p}\sum_{q'>q-5}2^{q'}\|\Delta_{q'}Q_{11}\|_{L^2}\nonumber\\
\le \|Q_{11}\|_{L^p}\sum_{q'>q-5}\|\Delta_{q'}\nabla Q_{11}\|_{L^2}\le \|Q_{11}\|_{L^p}2^{-qs}a_q(t)\|\nabla Q\|_{H^s}\nonumber
\end{eqnarray} where $a_q(t)=\sum_{q'>q-5}2^{-(q'-q)s}\bar a_{q'}(t)$.

\par On the other hand:

\begin{eqnarray}
\|\sum_{|q'-q|\le 5}\Delta_q (S_{q'-1}Q_{11}\Delta_{q'}Q_{11})\|_{L^p}\le \|Q_{11}\|_{L^p}\sum_{|q'-q|\le 5}\|\Delta_{q'} Q_{11}\|_{L^\infty}\nonumber\\ \le \|Q_{11}\|_{L^p}2^{-qs}\sum_{|q'-q|\le 5}2^{(q-q')s}2^{q's}\|\Delta_{q'}\nabla Q_{11}\|_{L^2}
\le \|Q_{11}\|_{L^p}2^{-qs}a_q(t)\|\nabla Q\|_{H^s}
\end{eqnarray}
\par The last two estimates prove Step $1$.

\smallskip\par{\it Step $2$} We assume the statement true for $j$ and we aim to prove it for $j+1$. We have
$$\Delta_q(Q_{11}^j Q_{11})=\sum_{q'>q-5}\Delta_q(S_{q'+2}Q_{11}^j\Delta_{q'}Q_{11})+\sum_{|q'-q|\le 5}\Delta_q(S_{q'-1}Q_{11}\Delta_q (Q_{11}^j))$$ and

\begin{eqnarray}
\sum_{q'>q-5}\|\Delta_q(S_{q'+2}Q_{11}^j\Delta_{q'}Q_{11}\|_{L^p}\le \sum_{q'>q-5}\|Q_{11}^j\|_{L^p}\|\Delta_{q'}Q_{11}\|_{L^\infty}\nonumber\\
\|Q_{11}^j\|_{L^p}\sum_{q'>q-5}2^{-q's}2^{q's}\|\nabla Q_{11}\|_{L^2}\le 2^{-qs}a_q(t)\|Q_{11}\|_{L^{pj}}^j \|\nabla Q\|_{H^s}
\nonumber
\end{eqnarray} where $a_q(t)=\sum_{q'>q-5}2^{-(q'-q)s}\bar a_q(t)$.

\par On the other hand, letting $r\stackrel{\rm{def}}{=}\frac{pj}{j-1}$  so that $\frac{1}{r}+\frac{1}{pj}=\frac{1}{p}$ we get:

\begin{eqnarray}
\sum_{|q'-q|\le 5} \|S_{q'-1}Q_{11}\Delta_{q'}Q_{11}^j\|_{L^p}\le \sum_{|q'-q|\le 5}\|Q_{11}\|_{L^{pj}}\|\Delta_{q'} Q_{11}^j\|_{L^r}\nonumber\\
\le \|Q_{11}\|_{L^{pj}}2^{-qs}a_q(t)\|Q\|_{L^{r(j-1)}}^{j-1}\|\nabla Q\|_{H^s}\le \|Q\|_{L^{pj}}^j2^{-qs}a_q(t)\|\nabla Q\|_{H^s}\nonumber
\end{eqnarray} where for the second inequality we used the inducation hypothesis.
\par The last two estimates show  Step $2$ and thus prove the lemma. $\Box$

      The lemma and estimate (\ref{est:easiest})  imply the claimed estimate (\ref{i8}).

\smallskip      The term $\mathcal{I}_k, k\in\{9,10,\dots,14\}$ is estimated exactly as the term $\mathcal{I}_{k-7}$ that we have already studied.

\smallskip       Using the commutator estimate (\ref{commutator}) with $p=2, r=\frac{2}{\varepsilon}, s=\frac{2}{1-\varepsilon}$ (where we restrict  $0<\varepsilon<\frac{1}{2}$) we get:
\begin{eqnarray}
|\mathcal{I}_{15}|=|\sum_{|q'-q|\le 5}\Big([\Delta_q, S_{q'-1}Q_{\alpha\beta}]\Delta_{q'}\textrm{tr}(Q\nabla u),\Delta_q\Delta Q\Big)|\le C\sum_{|q'-q|\le 5}2^{-q}\| S_{q'-1}\nabla Q\|_{L^{\frac{2}{\varepsilon}}}\|Q\|_{L^\infty} \|\Delta_{q'}\nabla u\|_{L^{\frac{2}{1-\varepsilon}}}\|\Delta\Delta_q Q\|_{L^2}\nonumber
\end{eqnarray}
using Bernstein inequality we have for $|q-q'|\leq 5$ and $\varepsilon\in(0,\frac 12)$,
$$2^{-q} \|\Delta_{q'}\nabla u\|_{L^{\frac{2}{1-\varepsilon}}}\leq C\|\Delta_{q'}u\|_{L^{\frac{2}{1-\varepsilon}}},$$
 and then, using the interpolation inequality (see \cite{cheminxu}, and also \cite{mp}, Lemma $10$):
$$\|f\|_{L^{2p}}\le C\sqrt{p}\|f\|_{L^2}^{\frac{1}{p}}\|\nabla f\|_{L^2}^{1-\frac{1}{p}}$$ with $p=\frac{1}{1-\varepsilon}\in[1,2]$, we obtain:

\begin{eqnarray}
|\mathcal{I}_{15}|\le C\sum_{|q'-q|\le 5} \|S_q\nabla Q\|_{L^{\frac{2}{\varepsilon}}}\|Q\|_{L^\infty}\|\Delta_q u\|_{L^2}^{1-\varepsilon}\|\Delta_q \nabla u\|_{L^2}^{\varepsilon}\|\Delta_q \Delta Q\|_{L^2},
\nonumber
\end{eqnarray}
where $C>0$ is constant independent of $\varepsilon\in (0,\frac 12)$.

      Using Young's inequality and assuming $0<\varepsilon<\frac{1}{2}$, $0<\eta<1$ we have $ab\le \frac{1-\varepsilon}{2}\frac{1}{\eta}^{\frac{2}{1-\varepsilon}}a^{\frac{2}{1-\varepsilon}}+
\frac{1+\varepsilon}{2}\eta^{\frac{2}{1+\varepsilon}}b^{\frac{2}{1+\varepsilon}}<\frac{1}{\eta^4} a^{\frac{2}{1-\varepsilon}}+\eta b^{\frac{2}{1+\varepsilon}}$ which implies, for appropriate $\eta$:

\begin{eqnarray}
|\mathcal{I}_{15}|\le C\sum_{|q'-q|\le 5}\Big(\|S_q\nabla Q\|_{L^{\frac{2}{\varepsilon}}}\|Q\|_{L^\infty}\Big)^{\frac{2}{1-\varepsilon}}\|\Delta_q u\|_{L^2}^2+\sum_{|q'-q|\le 5}\min\{\frac{\Gamma L^2}{100},\frac{\nu}{100}\}\|\Delta_q\nabla u\|_{L^2}^{\frac{2\varepsilon}{1+\varepsilon}}\|\Delta_q\Delta Q\|_{L^2}^{\frac{2}{1+\varepsilon}}
\nonumber
\end{eqnarray}

      We also use another form of Young's inequality, namely $ab\le \frac{\varepsilon a^{\frac{1+\varepsilon}{\varepsilon}}}{1+\varepsilon}+\frac{b^{1+\varepsilon}}{1+\varepsilon}<a^{\frac{1+\varepsilon}{\varepsilon}}+b^{1+\varepsilon}$ and obtain:

\begin{eqnarray}
|\mathcal{I}_{15}|\le C\sum_{|q'-q|\le 5}\bigg(\Big(\|S_q\nabla Q\|_{L^{\frac{2}{\varepsilon}}}\|Q\|_{L^\infty}\Big)^{\frac{2}{1-\varepsilon}}\|\Delta_q u\|_{L^2}^2+\frac{\nu}{100}\|\Delta_q\nabla u\|_{L^2}^2+\frac{\Gamma L^2}{100}\|\Delta_q\Delta Q\|_{L^2}^2\bigg)\nonumber\\ 
\le 2^{-2qs}\bigg(C\Big(\|S_q\nabla Q\|_{L^{\frac{2}{\varepsilon}}}\|Q\|_{L^\infty}\Big)^{\frac{2}{1-\varepsilon}}\|u\|_{H^s}^2+\frac{\nu}{100}\|\nabla u\|_{H^s}^2+\frac{\Gamma L^2}{100}\|\Delta Q\|_{H^s}\bigg)
\nonumber
\end{eqnarray}

\begin{eqnarray}
|\mathcal{I}_{16}|=|\sum_{|q'-q|\le 5} \big((S_{q'-1}Q_{\alpha\beta}-S_{q-1}Q_{\alpha\beta})\Delta_q\Delta_{q'} \textrm{tr}(Q\nabla u),\Delta\Delta_q Q_{\alpha\beta})|\nonumber\\
\le \sum_{|q-q'|\le 5}\|\Delta_q(Q\nabla u)\|_{L^\infty}\|S_{q'-1}Q-S_{q-1}Q\|_{L^2}\|\Delta_q \Delta Q\|_{L^2}\le \sum_{|q'-q|\le 5} 2^q\|Q\nabla u\|_{L^2}\|S_{q'-1}Q-S_{q-1}Q\|_{L^2}\|\Delta\Delta_q Q\|_{L^2}\nonumber\\
\le \sum_{|q'-q|\le 5}2^{-2qs}b_q(t)\|Q\|_{L^\infty}\|\nabla u\|_{L^2}\|\nabla Q\|_{H^s}\|\Delta Q\|_{H^s}
\nonumber
\end{eqnarray}

\begin{eqnarray}
|\mathcal{I}_{17}|=|\sum_{q'>q-5}\Big(\Delta_q(S_{q'+2}\textrm{tr}(Q\nabla u)\Delta_{q'}Q_{\alpha\beta}),\Delta_q \Delta Q_{\alpha\beta}\Big)|\le \sum_{q'>q-5}\|S_{q'+2}(Q\nabla u)\|_{L^\infty}\|\Delta_{q'} Q_{\alpha\beta}\|_{L^2}\|\Delta\Delta_q Q_{\alpha\beta}\|_{L^2}\nonumber\\
\le \sum_{q'>q-5} 2^{q'}\|Q\nabla u\|_{L^2} \|\Delta_{q'}Q_{\alpha\beta}\|_{L^2}\|\Delta\Delta_q Q_{\alpha\beta}\|_{L^2}\le \|Q\|_{L^\infty}\|\nabla u\|_{L^2} \sum_{q'>q-5} 2^{-(q+q')s}2^{q's}\|\Delta_{q'}\nabla Q_{\alpha\beta}\|_{L^2}2^{qs}\|\Delta\Delta_q Q_{\alpha\beta}\|_{L^2}\nonumber\\
\le 2^{-2qs}b_q(t)\|Q\|_{L^\infty} \|\nabla u\|_{L^2}\|\nabla Q\|_{H^s}\|\Delta Q\|_{H^s}
\nonumber
\end{eqnarray} where $b_q\stackrel{def}{=}\sum_{q'>q-5} 2^{-(q'-q)s}\tilde a_{q'}(t)a_q(t)\in l^1$ with $a_q(t),\tilde a_{q'}(t)\in l^2$.

      Using the commutator estimate (\ref{commutator}) with $p=2, q=\frac{2}{\varepsilon}, r=\frac{2}{1-\varepsilon}$ (where we  restrict  $0<\varepsilon<\frac{1}{2}$) we get:
\begin{eqnarray}
|\mathcal{I}_{18}|=|\big(S_{q-1}Q_{\alpha\beta}\sum_{|q'-q|\le 5} [\Delta_q, S_{q'-1} Q_{\gamma\delta}]\Delta_{q'} u_{\gamma,\delta},\Delta_q\Delta Q_{\alpha\beta}\big)|\le \|S_{q-1}Q\|_{L^\infty} \sum_{|q'-q|\le 5} \|[\Delta_q, S_{q'-1} Q]\Delta_{q'}\nabla u\|_{L^2}\|\Delta_q\Delta Q\|_{L^2}\nonumber\\
\le\|Q\|_{L^\infty} \sum_{|q'-q|\le 5} 2^{-q}\|S_{q'-1}\nabla Q\|_{L^{\frac{2}{\varepsilon}}}\|\Delta_{q'}\nabla u\|_{L^{\frac{2}{1-\varepsilon}}}\|\Delta\Delta_q Q\|_{L^2}
\nonumber
\end{eqnarray}

      We continue estimating exactly as in the proof of the estimates for the term $\mathcal{I}_{15}$ and obtain, for $0<\varepsilon<\frac 12$:

\begin{eqnarray}
|\mathcal{I}_{18}|\le 2^{-2qs}\bigg(C\Big(\|S_q\nabla Q\|_{L^{\frac{2}{\varepsilon}}}\|Q\|_{L^\infty}\Big)^{\frac{2}{1-\varepsilon}}\|u\|_{H^s}^2+\frac{\nu}{100}\|\nabla u\|_{H^s}^2+\frac{\Gamma L^2}{100}\|\Delta Q\|_{H^s}\bigg)
\nonumber
\end{eqnarray}

\begin{eqnarray}
|\mathcal{I}_{19}|=|\Big(S_{q-1}Q_{\alpha\beta}\big(\sum_{|q'-q|\le 5}(S_{q'-1} Q_{\gamma\delta}-S_{q-1}Q_{\gamma\delta})\Delta_q \Delta_{q'} u_{\gamma,\delta}\big),\Delta\Delta_q Q_{\alpha\beta}\Big)|\nonumber\\
\le \|Q\|_{L^\infty}\sum_{|q'-q|\le 5}\|S_{q'-1} Q-S_{q-1}Q\|_{L^2}\|\Delta_q \nabla u\|_{L^\infty}\|\Delta\Delta_q Q\|_{L^2}\nonumber\\\le\|Q\|_{L^\infty}\sum_{|q'-q|\le 5} 2^q \|S_{q'-1} Q-S_{q-1}Q\|_{L^2}\|\nabla u\|_{L^2}\|\Delta\Delta_q Q\|_{L^2}\le
\|Q\|_{L^\infty}\|\nabla u\|_{L^2}\sum_{|q'-q|\le 5}2^{-2qs}\|\nabla Q\|_{H^s}\|\Delta Q\|_{H^s}
\nonumber
\end{eqnarray}

\begin{eqnarray}
|\mathcal{I}_{20}|=|\big(S_{q-1} Q_{\alpha\beta}\sum_{q'>q-5} \Delta_q (S_{q'+2}u_{\gamma,\delta}\Delta_{q'}Q_{\gamma\delta}), \Delta_q \Delta Q_{\alpha\beta}\big)|\nonumber\\
\le \|Q\|_{L^\infty}\sum_{q'>q-5} \|S_{q'+2}u_{\gamma,\delta}\Delta_{q'}Q_{\gamma\delta}\|_{L^2}\|\Delta_q\Delta Q\|_{L^2}\le
\|Q\|_{L^\infty}\|\nabla u\|_{L^2}\sum_{q'>q-5} 2^{q'}\|\Delta_{q'}Q\|_{L^2}\|\Delta_q \Delta Q\|_{L^2}\nonumber\\
\le \|Q\|_{L^\infty}\|\nabla u\|_{L^2}\sum_{q'>q-5}\|\Delta_{q'}\nabla Q\|_{L^2}\|\Delta_q \Delta Q\|_{L^2}\le \|Q\|_{L^\infty}\|\nabla u\|_{L^2} \sum_{q'>q-5} 2^{-(q+q')s}2^{q's}\|\Delta_{q'}\nabla Q_{\alpha\beta}\|_{L^2}2^{qs}\|\Delta\Delta_q Q_{\alpha\beta}\|_{L^2}\nonumber\\
\le 2^{-2qs}b_q(t)\|Q\|_{L^\infty} \|\nabla u\|_{L^2}\|\nabla Q\|_{H^s}\|\Delta Q\|_{H^s}\nonumber
\end{eqnarray} where $b_q\stackrel{def}{=}\sum_{q'>q-5} 2^{-(q'-q)s}\tilde a_{q'}(t)a_q(t)\in l^1$ with $a_q(t),\tilde a_q(t)\in l^2$.

\begin{eqnarray}
|\mathcal{J}_1|=|(\Delta_q(u\nabla u),\Delta_q u)|=\underbrace{|\int S_{q-1}u\nabla\Delta_q u\cdot\Delta_q u|}_{\mathcal{J}_{1a}}+\underbrace{\sum_{|q'-q|\le 5} |\int [\Delta_q;S_{q'-1} u]\Delta_{q'}\nabla u \Delta_q u|}_{\mathcal{J}_{1b}}\nonumber\\
+\underbrace{\sum_{|q'-q|\le 5} |\int(S_{q'-1}u-S_{q-1}u)\Delta_q\Delta_{q'}\nabla u\Delta_q u|}_{\mathcal{J}_{1c}}+\underbrace{\sum_{q'> q-5}|\int \Delta_q (S_{q'+2}\nabla u\cdot\Delta_{q'}u)\Delta_q u|}_{\mathcal{J}_{1d}}\nonumber
 \end{eqnarray}

with
\begin{eqnarray}
|\mathcal{J}_{1a}|\le \|S_{q-1}u\|_{L^4}\|\Delta_q\nabla u\|_{L^2}\|\Delta_q u\|_{L^4}\le \|u\|_{L^2}^{\frac 12}\|\nabla u\|_{L^2}^{\frac 12} 2^{-2qs} b_q(t)\|\nabla u\|_{H^s}^{\frac 32}\|u\|_{H^s}^{\frac 12}\nonumber
\end{eqnarray}

\begin{eqnarray}
|\mathcal{J}_{1b}|=|\sum_{|q'-q|\le 5} \int [\Delta_q;S_{q'-1} u]\Delta_{q'}\nabla u \Delta_q u|\le C 2^{-q}\|S_{q-1}\nabla u\|_{L^4}\|\Delta_{q'}\nabla u\|_{L^2}\|\Delta_q u\|_{L^4}\nonumber\\
\le C \|u\|_{L^2}^{\frac 12}\|\nabla u\|_{L^2}^{\frac 12}2^{-2qs}b_q(t)\|\nabla u\|_{H^s}^{\frac 32}\|u\|_{H^s}^{\frac 12}\nonumber
\end{eqnarray}

\begin{eqnarray}
|\mathcal{J}_{1c}|\le \sum_{|q'-q|\le 5}\|(S_{q'-1}-S_{q-1})u\|_{L^4}\|\Delta_q\nabla u\|_{L^2}\|\Delta_q u\|_{L^4}\le C \|u\|_{L^2}^{\frac 12}\|\nabla u\|_{L^2}^{\frac 12} 2^{-2qs} b_q(t)\|\nabla u\|_{H^s}^{\frac 32}\|u\|_{H^s}^{\frac 12}\nonumber
\end{eqnarray}

\begin{eqnarray}
|\mathcal{J}_{1d}|=|\sum_{q'> q-5}\left(\Delta_q\left(S_{q'+2}\nabla u\Delta_{q'}u\right),\Delta_q u\right)|\le \sum_{q'> q-5}\|\Delta_q(S_{q'+2}\nabla u\Delta_{q'}u)\|_{L^{\frac 43}}\|\Delta_q u\|_{L^4}\nonumber\\
\le \sum_{q'> q-5}\|S_{q'+2}\nabla u\|_{L^2}\|\Delta_{q'}u\|_{L^4}\|\Delta_q u\|_{L^4}\le  \sum_{q'> q-5} C\|\nabla u\|_{L^2}\|\Delta_{q'}u\|_{L^2}^{\frac 12}\|\Delta_{q'}\nabla u\|_{L^2}^{\frac 12} \|\Delta_q u\|_{L^2}^{\frac 12}\|\Delta_q\nabla u\|_{L^2}^{\frac 12}\nonumber\\
\le C\|\nabla u\|_{L^2}2^{-qs}\sum_{q'> q-5} 2^{-q's}(2^{q's}\|\Delta_{q'}u\|_{L^2})^{\frac 12}(2^{q's}\|\Delta_{q'}\nabla u\|_{L^2})^{\frac 12}(2^{qs}\|\Delta_q u\|_{L^2})^{\frac 12}(2^{qs}\|\Delta_q\nabla u\|_{L^2})^{\frac 12}\nonumber\\
\le C\|\nabla u\|_{L^2}\|\nabla u\|_{H^s}\|u\|_{H^s}2^{-2qs}\bigg(2^{qs}\sum_{q'> q-5} c2^{-q's}a_{q'}(t)\bar a_{q'}(t)\Big)
\le C2^{-2qs} b_q(t)\|\nabla u\|_{L^2}\|\nabla u\|_{H^s}\|u\|_{H^s}\nonumber
\end{eqnarray} where $b_q(t)=\sum_{q'> q-5} 2^{-(q'-q)s}a_{q'}(t)\bar a_{q'}(t)\in l^1_q, \forall t\ge 0$.

\par We claim that we have:

\begin{eqnarray}
|\mathcal{J}_2|=|\int \Delta_q\left(   \partial_\alpha Q_{\gamma\delta}   \partial_\beta Q_{\gamma\delta}\right)\Delta_q u_{\alpha,\beta}|\le \|\Delta_q\left(   \partial_\alpha Q_{\gamma\delta}   \partial_\beta Q_{\gamma\delta}\right)\|_{L^2}\|\Delta_q\nabla u\|_{L^2}\nonumber\\
\le C2^{-2qs}b_q(t)\|   \partial_\alpha Q_{\gamma\delta}   \partial_\beta Q_{\gamma\delta}\|_{H^s}\|\nabla u\|_{H^s}
\le C2^{-2qs}b_q(t)\|\nabla Q\|_{L^2}^{\frac 12}\|\Delta Q\|_{L^2}^{\frac 12}\|\nabla Q\|_{H^s}^{\frac 12}\|\Delta Q\|_{H^s}^{\frac 12}\|\nabla u\|_{H^s}.
\label{est:j2}
\end{eqnarray} 

\par In order to prove the claim we write, using Bony's paraproduct decomposition (\ref{bonydecomp}):

\begin{eqnarray}
|\mathcal{J}_2|\le  \|\Delta_q \big(   \partial_\alpha Q_{\gamma\delta}   \partial_\beta Q_{\gamma\delta}\big)\|_{L^2}\|\Delta_q\nabla u\|_{L^2}\le \underbrace{\|S_{q-1} Q_{\gamma\delta,\alpha}\Delta_q Q_{\gamma\delta,\beta}\|_{L^2}\|\Delta_q\nabla u\|_{L^2}}_{\stackrel{\rm{def}}{=}A}\nonumber\\
+\underbrace{\sum_{|q'-q|\le 5}\|(S_{q'-1}Q_{\gamma\delta,\alpha}-S_{q-1}Q_{\gamma\delta,\alpha})\Delta_q\Delta_{q'}Q_{\gamma\delta,\beta}\|_{L^2}\|\Delta_q\nabla u\|_{L^2}}_{\stackrel{\rm{def}}{=}B}+\underbrace{\sum_{|q'-q|\le 5}\| [\Delta_q,S_{q'-1}Q_{\gamma\delta,\alpha}]\Delta_{q'}Q_{\gamma\delta,\beta}\|_{L^2}\|\Delta_q\nabla u\|_{L^2}}_{\stackrel{\rm{def}}{=}C}\nonumber\\
+\underbrace{ \sum_{q'>q-5}\|\Delta_q(S_{q'+2}Q_{\gamma\delta,\beta}\Delta_{q'}Q_{\gamma\delta,\alpha})\|_{L^2}\|\Delta_q\nabla u\|_{L^2}}_{\stackrel{\rm{def}}{=}D}
\end{eqnarray}
\par We estimate:
$$|A|\le \|\nabla Q\|_{L^4}\|\Delta_q\nabla Q\|_{L^4}\|\Delta_q\nabla u\|_{L^2}\le 2^{-2qs}b_q(t)\|\nabla Q\|_{L^2}^{\frac{1}{2}}\|\Delta Q\|_{L^2}^{\frac{1}{2}}\|\nabla Q\|_{H^s}^{\frac{1}{2}}\|\Delta Q\|_{H^s}^{\frac{1}{2}}\|\nabla u\|_{H^s}$$ and a similar estimate holds for $B$.
\begin{eqnarray}|C|\le \sum_{|q'-q|\le 5} 2^{-q}\|\Delta S_{q'-1} Q\|_{L^4}\|\Delta_{q'} \nabla Q\|_{L^4}\|\Delta_q\nabla u\|_{L^2}\le \|\nabla Q\|_{L^4}\|\Delta_q\nabla Q\|_{L^4}\|\Delta_q\nabla u\|_{L^2}\nonumber\\
\le 2^{-2qs}b_q(t)\|\nabla Q\|_{L^2}^{\frac{1}{2}}\|\Delta Q\|_{L^2}^{\frac{1}{2}}\|\nabla Q\|_{H^s}^{\frac{1}{2}}\|\Delta Q\|_{H^s}^{\frac{1}{2}}\|\nabla u\|_{H^s}\nonumber
\end{eqnarray}
\begin{eqnarray}
|D|\le \sum_{q'>q-5}\|\Delta_q(S_{q'+2}Q_{\gamma\delta,\beta}\Delta_{q'}Q_{\gamma\delta,\alpha})\|_{L^2}\|\Delta_q \nabla u\|_{L^2}\le \sum_{q'>q-5}\|\nabla Q\|_{L^4}\|\Delta_{q'}\nabla Q\|_{L^4}\|\Delta_q\nabla u\|_{L^2}\nonumber\\
\le 2^{-2qs}b_q(t)\|\nabla Q\|_{L^2}^{\frac{1}{2}}\|\Delta Q\|_{L^2}^{\frac{1}{2}}\|\nabla Q\|_{H^s}^{\frac{1}{2}}\|\Delta Q\|_{H^s}^{\frac{1}{2}}\|\nabla u\|_{H^s}\nonumber
\end{eqnarray} where $b_q(t)=\sum_{q'>q-5}2^{-(q'-q)s}a_{q'}(t)\bar a_q(t)$. 
\par The last three estimates imply the claimed estimate (\ref{est:j2}).

\begin{eqnarray}
|\mathcal{J}_3|=|\sum_{|q'-q|\le 5}\int[\Delta_q;S_{q'-1} Q_{\alpha\gamma}]\Delta_{q'}\Delta Q_{\gamma\beta}\Delta_q u_{\alpha,\beta}|\le \sum_{|q'-q|\le 5} \|[\Delta_q;S_{q'-1}Q_{\alpha\gamma}]\Delta_{q'}\Delta Q_{\gamma\beta}\|_{L^2}\|\Delta_q\nabla u\|_{L^2}\nonumber\\
\le C2^{-q}\|S_{q'-1}\nabla Q_{\alpha\gamma}\|_{L^\infty}\|\Delta_q\Delta Q_{\gamma\beta}\|_{L^2}\|\Delta_q\nabla u\|_{L^2}\le C2^{-\frac{q}{2}}\|S_{q'-1}\nabla Q\|_{L^4}\|\Delta_q\Delta Q\|_{L^2}\|\Delta_q\nabla u\|_{L^2}\nonumber\\
\le 2^{\frac q2}\|\nabla Q\|_{L^2}^{\frac 12}\|\Delta Q\|_{L^2}^{\frac 12}\|\Delta_q\Delta Q\|_{L^2}\|\Delta_q u\|_{L^2} \nonumber\\
\le C2^{-2qs}b_q(t)\|\nabla Q\|_{L^2}^{\frac 12}\|\Delta Q\|_{L^2}^{\frac 12}\|\Delta Q\|_{H^s}\|u\|_{H^s}^{\frac 12}\|\nabla u\|_{H^s}^{\frac 12}\nonumber
\end{eqnarray}

Concerning the term $\mathcal{J}_4$ we use that  $(S_{q'-1}Q_{\alpha\gamma}-S_{q-1}Q_{\alpha\gamma})$ is localized in a dyadic ring, so we have
$$\|S_{q'-1}Q_{\alpha\gamma}-S_{q-1}Q_{\alpha\gamma}\|_{L^\infty}\leq C2^{-\frac q2}\|\nabla Q\|_{L^4}\leq C 2^{-\frac q2}\|\nabla Q\|_{L^2}^{\frac 12}\|\Delta Q\|_{L^2}^{\frac 12},$$
and we obtain
$$|\mathcal{J}_4|=|\int\sum\limits_{|q'-q|\leq 5}(S_{q'-1}Q_{\alpha\gamma}-S_{q-1}Q_{\alpha\gamma})\Delta_q\Delta_{q'}\Delta Q_{\gamma\beta}\Delta_q u_{\alpha,\beta}|\leq C2^{-\frac q2}\|\nabla Q\|_{L^2}^{\frac 12}\|\Delta Q\|_{L^2}^{\frac 12}\|\Delta_q\Delta Q\|_{L^2}\|\Delta_q u_{\alpha,\beta}\|_{L^2}.$$

{     Using the fact that $\|\Delta_q u_{\alpha,\beta}\|_{L^2}\leq C 2^{\frac q2}2^{-qs} a_q^1(t)\|u\|_{H^s}^{\frac 12}\|\nabla u\|_{H^s}^{\frac 12}$ and $\|\Delta_q\Delta Q\|_{L^2}\leq C2^{-qs} a_q^2(t)\|\Delta Q\|_{H^s}$ and denoting $b_q(t)\stackrel{def}{=}a_q^1(t)a_q^2(t)$ we find
$$|\mathcal{J}_4|\leq C2^{-2qs} b_q(t)\|\nabla Q\|_{L^2}^{\frac 12}\|\Delta Q\|_{L^2}^{\frac 12}\|\Delta Q\|_{H^s}\|u\|_{H^s}^{\frac 12}\|\nabla u\|_{H^s}^{\frac 12}$$
The following term to estimate is $\mathcal{J}_5$. Using Bernstein inequalities $\|S_{q'+2}\Delta Q\|_{L^\infty}\leq C2^{q'}2^{\frac{q'}{2}}\|\nabla Q\|_{L^4}\leq C2^{\frac{3q'}{2}}\|\nabla Q\|_{L^2}^{\frac 12}\|\Delta Q\|_{L^2}^{\frac 12}$ and $\|\Delta_{q'} Q_{\alpha\gamma}\|_{L^2}\leq C2^{-\frac{3q'}{2}}\|\nabla \Delta_{q'} Q_{\alpha\gamma}\|_{L^2}^{\frac 12}\|\Delta_{q'}\Delta Q_{\alpha\gamma}\|_{L^2}^{\frac 12}$, we obtain
\begin{eqnarray}|\mathcal J_5|=|\sum_{q'>q-5}\int \Delta_q\big(S_{q'+2}\Delta Q_{\gamma\beta} \Delta_{q'}Q_{\alpha\gamma}\big)\Delta_q u_{\alpha,\beta}|\le |\sum_{q'>q-5} \|S_{q'+2}\Delta Q\|_{L^\infty}\|\Delta_{q'}Q\|_{L^2}\|\Delta_q\nabla u\|_{L^2}
\nonumber\\
\leq C\sum\limits_{q'>q-5}2^{\frac{3q'}{2}}\|\nabla Q\|_{L^2}^{\frac 12}\|\Delta Q\|_{L^2}^{\frac 12} 2^{-\frac{3q'}{2}}\|\Delta_{q'}\nabla Q\|_{L^2}^{\frac 12}\|\Delta_{q'}\Delta Q\|_{L^2}^{\frac 12}\|\Delta_q\nabla u\|_{L^2}\\
\le C\|\nabla Q\|_{L^2}^{\frac 12}\|\Delta Q\|_{L^2}^{\frac 12}\sum_{q'>q-5}2^{-q's}a_{q'}(t)\|\nabla Q\|_{H^s}^{\frac 12}\|\Delta Q\|_{H^s}^{\frac 12} 2^{-qs}\bar a_q(t)\|\nabla u\|_{H^s}\nonumber\\
\leq C2^{-2qs} b_q(t)\|\nabla Q\|_{L^2}^{\frac 12}\|\Delta Q\|_{L^2}^{\frac 12}\|\nabla Q\|_{H^s}^{\frac 12}\|\Delta Q\|_{H^s}^{\frac 12}\|\nabla u\|_{H^s}\nonumber
\end{eqnarray} where $b_q(t)=\sum_{q'>q-5} 2^{-(q'-q)s}a_{q'}(t)\bar a_q(t)$ .

\smallskip         The term $\mathcal{J}_k, k=6,7,8$ is estimated exactly as the term $\mathcal{J}_{k-3}$ that we have  already studied above. We also have that $\mathcal{J}_k=\mathcal{J}_{k-6}$ for $k\in\{ 9,\dots, 14\}$.

\smallskip      For $\mathcal{J}_{15}$ we apply Schwartz inequality together with the commutator estimate (\ref{commutator}) with $0<\varepsilon<\frac{1}{2}$ and $p=\frac{2}{1+\varepsilon}$, $r=\frac{2}{\varepsilon}$ and $s=2$  to obtain:

\begin{eqnarray}
|\mathcal{J}_{15}|=|\big(\sum_{|q'-q|\le 5}[\Delta_q, S_{q'-1}Q_{\alpha\beta}]\Delta_{q'}\textrm{tr}(Q\Delta Q),\Delta_q u_{\alpha,\beta}\big)|\le \sum_{|q'-q|\le 5} \|[\Delta_q,S_{q'-1}Q]\Delta_{q'}(Q\Delta Q)\|_{L^{\frac{2}{1+\varepsilon}}}\|\Delta_{q'}\nabla u\|_{L^{\frac{2}{1-\varepsilon}}}\nonumber\\
\le \sum_{|q'-q|\le 5}  2^{-q}\|S_{q'-1}\nabla Q\|_{L^{\frac{2}{\varepsilon}}}\|\Delta_{q'}(Q\Delta Q)\|_{L^2}\|\Delta_{q'}\nabla u\|_{L^{\frac{2}{1-\varepsilon}}}\nonumber\\
\le \sum_{|q'-q|\le 5}  2^{-q}\|S_{q'-1}\nabla Q\|_{L^{\frac{2}{\varepsilon}}}\|Q\|_{L^\infty}\|\Delta_{q'}\Delta Q\|_{L^2}\|\Delta_{q'}\nabla u\|_{L^{\frac{2}{1-\varepsilon}}}\nonumber
\end{eqnarray}

       We continue estimating exactly as in the proof of the estimates for the term $\mathcal{I}_{15}$ and obtain:

\begin{eqnarray}
|\mathcal{J}_{15}|\le 2^{-2qs}\bigg(C\Big(\|S_q\nabla Q\|_{L^{\frac{2}{\varepsilon}}}\|Q\|_{L^\infty}\Big)^{\frac{2}{1-\varepsilon}}\|u\|_{H^s}^2+\frac{\nu}{100}\|\nabla u\|_{H^s}^2+\frac{\Gamma L^2}{100}\|\Delta Q\|_{H^s}\bigg)\nonumber
\end{eqnarray}

\begin{eqnarray}
|\mathcal{J}_{16}|\le \sum_{|q'-q|\le 5}|\Big( (S_{q'-1}Q_{\alpha\beta}-S_q Q_{\alpha\beta})\Delta_q\Delta_{q'}\textrm{tr}(Q\Delta Q),\Delta_q u_{\alpha,\beta}\Big)|\nonumber\\
\le \sum_{|q'-q|\le 5}\|Q\|_{L^\infty}\|S_{q'-1}Q-S_q Q\|_{L^\infty}\|\Delta_q\Delta Q\|_{L^2}\|\Delta_q \nabla u\|_{L^2}\le \|Q\|_{L^\infty}\|\nabla u\|_{L^2}\sum_{|q'-q|\le 5} 2^q\|S_{q'-1}Q-S_q Q\|_{L^2}\|\Delta_q\Delta Q\|_{L^2}\nonumber\\
\le 2^{-2qs}b_q(t)\|Q\|_{L^\infty}\|\nabla u\|_{L^2}\|\nabla Q\|_{H^s}\|\Delta Q\|_{H^s}\nonumber
\end{eqnarray}

\begin{eqnarray}
|\mathcal{J}_{17}|\le \sum_{q'>q-5}|(\Delta_q (S_{q'+2}\textrm{tr}(Q\Delta Q)\Delta_{q'}Q_{\alpha\beta}),\Delta_q u_{\alpha,\beta})|\le \sum_{q'>q-5}\|\Delta_q\textrm{tr}(Q\Delta Q)\|_{L^2}\|\Delta_{q'}Q\|_{L^\infty}\|\Delta_q \nabla u\|_{L^2}\nonumber\\
\le \|Q\|_{L^\infty}\|\nabla u\|_{L^2} \sum_{q'>q-5}\|\Delta_q' \nabla Q\|_{L^2}\|\Delta_q \Delta Q\|_{L^2}\le 2^{-2qs}b_q(t)\|Q\|_{L^\infty}\|\nabla u\|_{L^2}\|\nabla Q\|_{H^s}\|\Delta Q\|_{H^s}\nonumber
\end{eqnarray}

 For $\mathcal{J}_{18}$ we apply Schwartz inequality together with the commutator estimate (\ref{commutator}) with $0<\varepsilon<\frac{1}{2}$ and $p=\frac{2}{1+\varepsilon}$, $r=\frac{2}{\varepsilon}$ and $s=2$  to obtain:

\begin{eqnarray}
|\mathcal{J}_{18}|\le |\big(S_{q-1}Q_{\alpha\beta}\sum_{|q'-q|\le 5} [\Delta_q, S_{q'-1}Q_{\gamma\delta}]\Delta_{q'}\Delta Q_{\gamma\delta},\Delta_q u_{\alpha,\beta}\big)|\le \|Q\|_{L^\infty}\sum_{|q'-q|\le 5}\|[\Delta_q, S_{q'-1} Q]\Delta_{q'}\Delta Q\|_{L^{\frac{2}{1+\varepsilon}}}\|\Delta_q \nabla u\|_{L^{\frac{2}{1-\varepsilon}}}\nonumber\\
\le\|Q\|_{L^\infty}\sum_{|q'-q|\le 5} 2^{-q}\|S_{q'-1}\nabla Q\|_{L^{\frac{2}{\varepsilon}}}\|\Delta_{q'}\Delta Q\|_{L^2}\|\Delta_q \nabla u\|_{L^{\frac{2}{1-\varepsilon}}}\nonumber
\end{eqnarray} 

       We continue estimating exactly as in the proof of the estimates for the term $\mathcal{I}_{15}$ and obtain:

\begin{eqnarray}
|\mathcal{J}_{18}|\le 2^{-2qs}\bigg(C\Big(\|S_q\nabla Q\|_{L^{\frac{2}{\varepsilon}}}\|Q\|_{L^\infty}\Big)^{\frac{2}{1-\varepsilon}}\|u\|_{H^s}^2+\frac{\nu}{100}\|\nabla u\|_{H^s}^2+\frac{\Gamma L^2}{100}\|\Delta Q\|_{H^s}\bigg)\nonumber
\end{eqnarray}

\begin{eqnarray}
|\mathcal{J}_{19}|\le |\Big(S_{q-1}Q_{\alpha\beta}\sum_{|q'-q|\le 5} (S_{q'-1} Q_{\gamma\delta}-S_{q-1} Q_{\gamma\delta})\Delta_q \Delta_{q'}\Delta Q_{\gamma\delta},\Delta_q u_{\alpha,\beta}\Big)|\nonumber\\
\le \|Q\|_{L^\infty}\sum_{|q'-q|\le 5} \|S_{q'-1}Q-S_{q-1} Q\|_{L^\infty}
\|\Delta_q \Delta Q\|_{L^2}\|\Delta_q \nabla u\|_{L^2}\le \|Q\|_{L^\infty}\|\nabla u\|_{L^2}\sum_{|q'-q|\le 5}\|\Delta_q\nabla Q\|_{L^2}\|\Delta_q\Delta Q\|_{L^2}\nonumber\\
\le 2^{-2qs}b_q(t)\|Q\|_{L^\infty}\|\nabla u\|_{L^2}\|\nabla Q\|_{H^s}\|\Delta Q\|_{H^s}\nonumber
\end{eqnarray}

\begin{eqnarray}
|\mathcal{J}_{20}|\le |\Big(S_{q-1}Q_{\alpha\beta}\sum_{q'>q-5}\Delta_q (S_{q'+2}\Delta Q_{\gamma\delta}\Delta_{q'}Q_{\gamma\delta}),\Delta_q u_{\alpha,\beta}\Big)|\nonumber\\
\le \|Q\|_{L^\infty}\sum_{q'>q-5}\|\Delta_{q'}Q\|_{L^\infty}\|\Delta_q \Delta Q\|_{L^2}\|\Delta_q \nabla u\|_{L^2}\le \|Q\|_{L^\infty}\|\nabla u\|_{L^2}\sum_{q'>q-5}\|\Delta_{q'}\nabla Q\|_{L^2}\|\Delta_q\Delta Q\|_{L^2}\nonumber\\
\le \|Q\|_{L^\infty}\|\nabla u\|_{L^2}\sum_{q'>q-5}2^{-(q'+q)s}2^{q's}\|\Delta_{q'}\nabla Q\|_{L^2}2^{qs}\|\Delta_q\Delta Q\|_{L^2}\le
2^{-2qs}b_q(t)\|Q\|_{L^\infty}\|\nabla u\|_{L^2}\|\nabla Q\|_{H^s}\|\Delta Q\|_{H^s}\nonumber
\end{eqnarray} where $b_q(t)=\sum_{q'>q-5}2^{-(q'-q)s}a_{q'}(t)\bar a_q(t)$.

      Finally, we claim that:
\begin{eqnarray}
|\mathcal{J}_{21}|+|\mathcal{J}_{22}|+|\mathcal{J}_{23}|\le 2^{-2qs}b_q(t) \bigg[C\big(\sum_{j=2}^5 \|Q\|_{L^{2(j-1)}}^{j-1}\big)^2\|\nabla Q\|_{H^s}^2+\| u\|_{H^s}^2+\frac{\nu}{100}\|\nabla u\|_{H^s}^2\big]
\label{lastJs}
\end{eqnarray}

      In order to prove the above estimate, we observe that the simplest terms are those of the form $(\Delta_q Q_{\alpha\beta},\Delta_q u_{\alpha,\beta})$ that can be easily estimated:

\begin{equation}
|(\Delta_q Q_{\alpha\beta},\Delta_q u_{\alpha,\beta})|\le C \|\Delta_q \nabla  Q\|_{L^2}\| \Delta_q u\|_{L^2}\le  2^{-2qs}b_q(t)\|\nabla Q\|_{H^s}\|u\|_{H^s}\nonumber
\end{equation}

      For the rest of the terms we just consider a generic term from $\mathcal{J}_{21},\mathcal{J}_{22},\mathcal{J}_{23}$, namely $(\Delta_q(Q_{11}^j),\Delta_q u_{\alpha,\beta})$ where $2\le j\le 5$ and use Lemma ~\ref{lemma:deltaqpowers} to obtained the claimed estimate (\ref{lastJs}).

      Putting together the estimates for all terms, multiplying by $2^{2qs}$ and taking the sum in $q$, observing that we can write any sequence $b_q\in l^1_q$ as $b_q=a_q\cdot\bar a_q$ with $a_q,\bar a_q\in l^2_q$, using $ab\leq C\varepsilon^{-1} a^2+\varepsilon b^2$, with appropriate $\varepsilon$, we obtain the claimed estimate (\ref{longest}).

\end{document}